\documentclass[12pt,a4paper]{article}
\usepackage[utf8]{inputenc}
\usepackage{amsmath}
\usepackage{amsfonts}
\usepackage{dsfont}
\usepackage{mathrsfs}
\usepackage{fullpage,amssymb,amsmath,amsfonts,hyperref,amsthm,paralist,graphicx,color,float, mathtools,tikz,caption,subcaption}
\usepackage{here}
\usepackage[normalem]{ulem}
\usepackage{centernot}
\usepackage{ifthen,xkeyval, tikz, calc, graphicx}
\usepackage{xcolor}
\usepackage{amssymb}
\usepackage{ragged2e}
\usepackage{graphicx}
\usepackage[margin=0.85 in]{geometry}
\usepackage{enumitem}

\newtheorem{theorem}{Theorem}[section]

\newtheorem{corollary}[theorem]{Corollary}

\newtheorem{definition}[theorem]{Definition}

\newtheorem{lemma}[theorem]{Lemma}

\newtheorem{proposition}[theorem]{Proposition}

\newtheorem{remark}[theorem]{Remark}

\def\acknowledgementsname{Acknowledgments}
\newenvironment{acks}[1][\acknowledgementsname]{\noindent\textbf{#1.}\space\ignorespaces}{\par}

\title{Critical Exponents for Marked Random Connection Models}

\author{%
  Alejandro Caicedo\footnote{Ludwig-Maximilians-Universität München, Mathematisches Institut, Theresienstraße 39, 80333 München, Germany;
    Email: caicedo@math.lmu.de
    \includegraphics[height=1em]{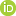}~\url{https://orcid.org/0009-0009-4299-6356}}
  \and 
  Matthew Dickson\footnote{University of British Columbia, Department of Mathematics, Vancouver, BC, Canada, V6T 1Z2;
  Email: dickson@math.ubc.ca
  \includegraphics[height=1em]{ORCIDiD_icon16x16.png}~\url{https://orcid.org/0000-0002-8629-4796}}}

\definecolor{darkorange}{RGB}{255,165,0}
\definecolor{altviolet}{RGB}{139,0,139}
\definecolor{turquoise}{RGB}{64,224,208}

\usepackage{xcolor}
\definecolor{nicegreen}{rgb}{0.0, 0.5, 0.0}
\definecolor{niceblue}{rgb}{0.36, 0.54, 0.66}
\usepackage[normalem]{ulem}

\newcommand{\p}{\mathbb P}
\newcommand{\pl}{\mathbb P_\lambda}
\newcommand{\plg}{\mathbb P_{\lambda,\gamma}}
\newcommand{\el}{\mathbb E_\lambda}
\newcommand{\elg}{\mathbb E_{\lambda,\gamma}}
\newcommand{\E}{\mathbb E}
\newcommand{\R}{\mathbb R}
\newcommand{\Rd}{\mathbb R^d}

\newcommand{\X}{\mathbb X}
\newcommand{\N}{\mathbb N}
\newcommand{\Complex}{\mathbb C}

\newcommand{\Ical}{\mathcal I}
\newcommand{\Jcal}{\mathcal J}
\newcommand{\Ncal}{\mathcal N}

\newcommand{\G}{\mathcal G}
\newcommand{\dd}{\, \mathrm{d}}
\newcommand{\C}{\mathscr {C}}
\newcommand{\Leb}{\mathrm{Leb}}
\newcommand{\thinn}[1]{\langle #1 \rangle}
\newcommand{\piv}[1]{{\mathrm {Piv}}\left(#1\right)}
\newcommand{\con}[2]{#1 \longleftrightarrow #2}
\newcommand{\ncon}[2]{#1 \centernot\longleftrightarrow #2}
\newcommand{\dcon}[2]{#1 \Longleftrightarrow #2}
\newcommand{\conn}[3]{#1 \longleftrightarrow #2\textrm { in } #3}
\newcommand{\nconn}[3]{#1 \centernot\longleftrightarrow #2\textrm { in } #3}
\newcommand{\dconn}[3]{#1 \Longleftrightarrow #2 \textrm { in } #3}
\newcommand{\ndconn}[3]{#1 \centernot\Longleftrightarrow #2\textrm { in } #3}

\newcommand{\clx}[2]{\mathscr{C}(#1,#2)}
\newcommand{\cls}[2]{\mathscr{C}^#2(#1)}
\newcommand{\sqarrow}{\leftrightsquigarrow}
\DeclareMathOperator*{\argmax}{arg\,max}

\newcommand{\xconn}[4]{#1 \xleftrightarrow{\,#4\,} #2\textrm { in } #3}

\newcommand{\inconn}[4]{#1 \overset{#2}{\longleftrightarrow} #3\textrm { in } #4 }
\newcommand{\sqconn}[3]{#1 \leftrightsquigarrow #2\textrm { in } #3 }

\newcommand{\thinning}[2]{#1_{\langle #2 \rangle}}
\newcommand{\adjconn}[3]{#1 \sim #2\textrm { in } #3}
\newcommand{\tlam}{\tau_\lambda}
\newcommand{\tlamo}{\tau_\lambda^\circ}
\newcommand{\tlg}{\tau_{\lambda,\gamma}}

\newcommand{\trilam}{\triangle_\lambda}

\newcommand{\Vol}[1]{{\rm Vol}_{#1}}

\newcommand{\e}{\text{e}}

\newcommand{\connf}{\varphi}

\newcommand{\Z}{\mathbb{Z}}
\newcommand{\origin}[1]{{o_{#1}}}
\newcommand{\Pcal}{\mathcal P}
\newcommand{\Kcal}{\mathcal K}

\newcommand{\Pmin}{{\mathcal P}_{\mathrm{min}} }
\DeclareMathOperator*{\esssup}{ess\,sup}
\DeclareMathOperator*{\essinf}{ess\,inf}
\newcommand{\Ecal}{\mathcal E}
\newcommand{\xbar}{\overline{x}}
\newcommand{\ybar}{\overline{y}}
\newcommand{\ubar}{\overline{u}}
\newcommand{\zerobar}{\overline{0}}
\DeclarePairedDelimiter\abs{\lvert}{\rvert}
\DeclarePairedDelimiter\norm{\lVert}{\rVert}
\DeclarePairedDelimiterX{\inner}[2]{\langle}{\rangle}{#1, #2}
\newcommand{\pla}{\mathbb P_{\lambda}}
\newcommand{\plaga}{\mathbb P_{\lambda,\gamma}}

\newcommand{\OpNorm}[1]{\left\lVert#1\right\rVert_{\rm op}}
\newcommand{\OneNorm}[1]{\left\lVert#1\right\rVert_{1,\infty}}

\newcommand{\InfNorm}[1]{\left\lVert#1\right\rVert_{\infty,\infty}}

\newcommand{\Id}{\mathds 1}
\newcommand{\Optlam}{\mathcal{T}_\lambda}
\newcommand{\OptlamT}{\mathcal{T}_{\lambda_O}}

\newcommand{\Optnlam}{\mathcal{T}^{(n)}_\lambda}

\newcommand{\Opconnf}{\mathcal{D}}

\newcommand{\ConstantTriangle}{C_\Delta}
\newcommand{\cbar}[1]{\overline{c}_{#1}}
\newcommand{\MagSupTilde}{\widetilde{M}_{\sup}}

\newcommand{\rom}[1]{\uppercase\expandafter{\romannumeral #1\relax}}

\numberwithin{equation}{section}

\begin{document}

\maketitle.                                  

\begin{abstract}
    Here we prove critical exponents for Random Connections Models (RCMs) with random marks. The vertices are given by a marked Poisson point process on $\mathbb{R}^d$ and an edge exists between any pair of vertices independently with a probability depending upon their spatial displacement and on their respective marks. Given conditions on the edge probabilities, we prove mean-field lower bounds for the susceptibility and percolation functions. In particular, we prove the equality of the susceptibility and percolation critical intensities. If we assume that a form of the triangle condition holds, then we also prove that the susceptibility, percolation and cluster tail critical exponents exist and take their mean-field values. Our proof approach adapts the differential inequality and magnetization function approaches that have been previously applied to discrete homogeneous settings to our continuum marked setting. This includes a proof of the analyticity of the magnetization function in the required parameter regime.
\end{abstract}



\section{Introduction}
\label{sec:Introduction}
\paragraph{A Preview on Lattice Percolation.}
Bond percolation on $\Z^d$ is the most classical model of percolation. In this model each pair of nearest-neighbour vertices in $\Z^d$ forms an edge independently with probability $p$. Then there exists a critical probability $p_c\in\left[0,1\right]$ such that $p>p_c$ implies that there almost surely exists an infinite cluster (percolation does occur), and $p<p_c$ implies that there almost surely is no infinite cluster (percolation does not occur). It is known that for $d\geq 2$ the critical probability $p_c\in\left(0,1\right)$ and the phase transition is said to be non-trivial. There are a number of properties that undergo a transition at this critical value $p_c$, including the percolation probability (the probability that the origin is part of an infinite cluster) and the susceptibility (the expected size of the origin's cluster). Since these seem to undergo a transition at the critical probability, it is natural to ask questions about their behaviour as $p\nearrow p_c$ and as $p\searrow p_c$. This behaviour is often characterised by a so-called \textit{critical exponent}. There are two general methods that have been used in this bond percolation model and similar models to study such near-critical properties. On one hand for $d=2$ there are particular special properties and techniques such as planarity and complex analysis tools that allow the study of behaviour near the critical point, see for example \cite{SmirnWerne01, Werne09, BolRio06}. On the other hand when $d$ is large (at least $d \geq 6$) Hara and Slade \cite{Har90} have proven the triangle condition holds which allows a proof that various critical exponents exist and take on their mean-field behaviour, see for example \cite{AizNew84,AizBar91,Ngu87,KozNac09,KozNac11,HeyHof17,hutchcroft2022derivation}.

\paragraph{The Marked Random Connection Model.}
Our focus in this paper is more in the spirit of the second approach. We consider a much richer family of models of percolation called the Marked Random Connection Model (MRCM), and for the percolation and susceptibility critical exponents we derive a mean-field bound that holds under fairly weak conditions on the MRCM. Then we show that if a version of the triangle condition holds for this model then this bound is sharp in the sense that the critical exponents indeed do take this value - we show mean-field behaviour. The connection between the triangle condition and large values of $d$ for the MRCM has been made in \cite{DicHey2022triangle}, in which they show that for suitable families of MRCMs the triangle condition holds in sufficiently high spatial dimensions. This had been performed previously in \cite{HeyHofLasMat19} for non-marked plain Random Connection Models. The conditions in these papers are not required in our analysis, since we simply assume that the triangle condition holds.

We now give an informal description of the MRCM. A more formal definition can be found in \cite{DicHey2022triangle}. We consider a random graph embedded in $\Rd$, where the vertices are assigned by a homogeneous Poisson point process with a given intensity $\lambda\geq 0$. Then each vertex is randomly and independently assigned a ``mark" from a Polish space $\Ecal$ with probability distribution $\Pcal$ (and the Borel $\sigma$-algebra). Equivalently - and the approach we pursue in this paper - we can denote $\X=\R^d\times \Ecal$ and our vertex set $\eta$ is given by a Poisson point process on $\X$ with intensity $\lambda\nu$, where $\nu = \Leb\otimes \Pcal$ ($\Leb$ denotes the Lebesgue measure on $\Rd$). The edges of our graph then occur randomly and independently with a probability that depends upon the spatial displacement of the two concerned vertices and upon their two marks. This probability is given by a measurable and symmetric adjacency function $\connf\colon \X^2\to\left[0,1\right]$. As mentioned above, we ask that $\connf$ depends only upon the two marks and the spatial displacement of the two vertices. This means that if $x=\left(\xbar,a\right)$ and $y=\left(\ybar,b\right)$ for $\xbar,\ybar\in\Rd$ and $a,b\in\Ecal$ then $\connf(x,y)=\connf(\xbar-\ybar;a,b)$ only. We then denote a realization of the full random graph by $\xi$, which means that $\xi$ is a realization of not only the marked point process but also includes the realization of the edges between them. Note that since both the distribution of the vertex set $\eta$ and the adjacency function $\connf$ are spatially translation invariant, the whole random graph $\xi$ is also spatially translation invariant.

\paragraph{Augmentation by Fixed Points.}
We also want to consider augmented versions of $\xi$. Given a vertex set $\eta$ and a point $x\in\X$, we denote by $\eta^x$ the vertex set where we have inserted one vertex at $x$. Given a graph $\xi$ we get $\xi^x$ by augmenting the vertex set $\eta$ to get $\eta^x$ as above and $\xi^x$ then inherits all the edges that were present in $\xi$. The graph $\xi^x$ then gains new edges by having the new vertex $x$ form edges with the vertices of $\eta$ independently of everything else with probability described by the same adjacency function $\connf$. This can be easily generalised to get $\eta^{x_1,\dots,x_m}$ and $\xi^{x_1,\ldots,x_m}$ for any $m\in\N$.

\paragraph{Percolation Quantities.}
Given $x,y\in\eta$, we say that $x$ and $y$ are adjacent, or $x\sim y$, in $\xi$ if there exists an edge with end vertices $x$ and $y$. We say that $x$ and $y$ are connected in $\xi$, or $\conn{x}{y}{\xi}$, if there exists a finite sequence of adjacent vertices in $\xi$ starting at $x$ and ending at $y$. That is, either $x\sim y$ in $\xi$ or there exist $n\in\N$ and $x_1,\ldots,x_{n}\in\eta$ such that $x\sim x_1$, $y\sim x_{n}$, and $x_i\sim x_{i+1}$ in $\xi$ for all $i\in\left\{1,\ldots,n-1\right\}$. We can now ask about the probability that two vertices are connected. Furthermore given $x\in\eta$ and a (possibly random) set, $B\subset\X$, we say that $\conn{x}{B}{\xi}$ if there exists $y\in \eta\cap B$ such that $\conn{x}{y}{\xi}$. Since the intensity measure $\nu$ is non-atomic, any given point $x\in\X$ is almost surely not a vertex of $\xi$. We therefore use the augmented graph $\xi^{x,y}$ in the following definition.

\begin{definition}(Two-point function)
For $\lambda\geq 0$, define the two-point function $\tlam\colon \X^2 \to [0,1]$ by 
\begin{align}
\tlam(x,y) := \pla( \conn x y \xi^{x,y}).
\end{align}
\end{definition}

It is also be useful to define the set of all point connected to some root vertex. We call this the cluster of that point.
\begin{definition}(Cluster of a point)
For a given point $x \in \X$ define the cluster of $x$ in $\xi$ as
\begin{align}
    \C(x,\xi) := \{y \in \eta \mid \conn x y \xi\}.
\end{align}
Since a fixed point $x \in \X$ is almost surely not in $\xi$ it does not generally make sense to consider $\C(x,\xi)$. Instead we usually consider clusters of the type $\C(x,\xi^x)$ and for the sake of not making the notation too heavy we denote $\C(x) := \C\left(x,\xi^x\right)$.
\end{definition}

Here we define two functions which in turn can be used to define a critical transition for the MRCM. In the following definitions we denote $\origin{a} = \left(\zerobar,a\right)\in\X$, where $\zerobar\in\Rd$ is the spatial origin and $a\in\Ecal$. The spatial translation invariance of the MRCM ensures that the choice of $\zerobar$ plays no special role, but the potential inhomogeneity of the mark space $\Ecal$ allows for the possibility that the choice of $a$ could be impactful.
\begin{definition} (Percolation function)
For $\lambda\geq 0$, we define the percolation function $\theta_\lambda\colon\Ecal\to \left[0,1\right]$ as
\begin{align}
    \theta_\lambda(a) &:= \pla\left(\abs*{\C(\origin{a})}=\infty\right).
\end{align}
\end{definition}

\begin{definition} (Susceptibility function)
For $\lambda\geq 0$, we define the susceptibility function $\chi_\lambda\colon\Ecal\to \left[0,\infty\right]$ as
\begin{align}
    \chi_\lambda(a) &:= \E_\lambda\left[\abs*{\C(\origin{a})}\right].
\end{align}
\end{definition}

\paragraph{Critical Behaviour.}
Note that for all $a\in\Ecal$, both $\theta_\lambda(a)$ and $\chi_\lambda(a)$ are non-decreasing in $\lambda$. It is then natural to define one critical intensity as the point at which $\theta_\lambda$ becomes positive, and another critical intensity as the point at which $\chi_\lambda$ becomes infinite. In the following definition we use the $L^p$-norms to determine the size of these functions.

\begin{definition}(Critical intensities)
We define the critical intensities
\begin{align}
    \lambda_c &:= \sup\left\{\lambda\geq 0\colon \theta_\lambda(a)=0 \text{ for $\Pcal$-a.e. }a\in\Ecal\right\},\label{eqn:percolationCritical}\\
    \lambda_T^{(p)} &:= \inf\left\{\lambda\geq 0\colon \norm{\chi_\lambda}_p =\infty\right\},\qquad\text{for all }p\in\left[1,+\infty\right].\label{eqn:susceptibilityCritical}
\end{align}

\end{definition}

First we observe that the picture becomes far simpler in the case where $\Ecal$ is a singleton - that is for the plain (single-mark) RCM. In that case it is clear that the variable $p$ plays no role and we only have two critical intensities $\lambda_c$ and $\lambda_T$. It was also proven in \cite{Mee95} that $\lambda_c=\lambda_T$ for the plain RCM.

Returning to the MRCM, there are still some simplifications that can easily be made. Note that for every $a\in\Ecal$, we have $\theta_\lambda(a)>0$ implies $\chi_\lambda(a) = \infty$. Therefore $\lambda_T^{(p)} \leq \lambda_c$ for all $p$. From Jensen's inequality it is also clear that $\lambda_T^{(p)}$ is non-increasing in $p$. Therefore we have the partial ordering
\begin{equation}
    0\leq \lambda_T^{(\infty)} \leq \lambda_T^{(p_2)} \leq \lambda_T^{(p_1)} \leq \lambda_T^{(1)} \leq \lambda_c
\end{equation}
for all $1\leq p_1\leq p_2 \leq \infty$.

Given the condition $\esssup_{a,b\in\Ecal}\int\connf(\xbar;a,b)\dd \xbar <\infty$ (we call this \ref{Assump:BoundExpectedDegree} below), \cite[Proposition~2.1]{DicHey2022triangle} proves that $\lambda^{(p)}_T$ is also $p$-independent. When we are working under this condition, we write $\lambda_T$ for this value. In the specific context of the Boolean disc model with bounded random radii, \cite[Proposition~2.2]{DicHey2022triangle} further shows that $\lambda_c=\lambda_T$. However, the question remained open in any greater generality. We show in Corollary~\ref{thm:Sharpness of Phase Transition} that this equality of the critical intensities holds under Assumptions \ref{Assump:BoundExpectedDegree} and \ref{Assump:AllReachablebySome}, generalising their result.

In Lemma~\ref{lem:nontriviallambdaO} we prove that $\lambda_c<\infty$, and this same lemma can be used to show that $\lambda_T>0$ (defined if \ref{Assump:BoundExpectedDegree} holds). This means that the the critical intensities $\lambda_c$ and $\lambda_T$ are non-trivial, and we are now interested in deriving critical exponents. These are numbers that allows us to better understand the behaviour of a given model near or at criticality. In this paper we are going to derive three different exponents, namely $\gamma$, $\beta$ and $\delta$. In the following definition and hereafter we use the notation $\left(X\right)_+=\max\left\{X,0\right\}$ for $X\in\R$.

\begin{definition}(Critical Exponents)
\label{def_crit_exp}

\begin{itemize}
    \item \textbf{$\gamma$ Exponent}. Suppose Assumption~\ref{Assump:BoundExpectedDegree} holds so that $\lambda_T$ is defined. We then say that $\gamma$ exists in the bounded ratio sense if there exists a constant $\gamma \in \left(0,\infty\right)$ as well as two other constants $C, C' \in \left(0,\infty\right)$ such that 
    \begin{equation}
        C\left(\lambda_T-\lambda\right)^{-\gamma} \leq \norm*{\chi_\lambda}_p \leq C'\left(\lambda_T-\lambda\right)^{-\gamma}
    \end{equation}
    for $\lambda<\lambda_T$ and for all $p\in\left[1,\infty\right]$.
    
    \item \textbf{$\beta$ Exponent.} We say that $\beta$ exists in the bounded ratio sense if there exists a constant $\beta \in \left(0,\infty\right)$ as well as three other constants $\epsilon, C, C' \in \left(0,\infty\right)$ such that 
\begin{equation}
    C\left(\lambda - \lambda_c\right)^\beta_+ \leq \norm*{\theta_\lambda}_p \leq C'\left(\lambda - \lambda_c\right)^\beta_+
\end{equation}
for $\lambda < \lambda_c + \epsilon$ and for all $p\in\left[1,\infty\right]$.

\item \textbf{$\delta$ Exponent.}
We say that $\delta$ exists in the bounded ratio sense if there exists a constant $\delta \in \left(0,\infty\right)$ as well as two other constants $C, C' \in \left(0,\infty\right)$ such that $\forall n \in \N$ and $\Pcal$-almost every $a\in\Ecal$
\begin{equation}
    Cn^{-\frac{1}{\delta}} \leq \mathbb{P}_{\lambda_c}\left(\abs*{\C\left(\origin{a}\right)}\geq n\right) \leq C'n^{-\frac{1}{\delta}}.
\end{equation}
\end{itemize}
\end{definition}

There are many other critical exponents that one can consider to better understand the behaviour of a model, but in this paper we focus our study on these three. 

\subsection{Results}

Our results require certain subsets of three assumptions. For the first two we introduce the following notation. Given $a,b\in\Ecal$ and $k\geq 1$, let us define
\begin{align}
    D(a,b) &:= \int_{\Rd}\connf\left(\xbar;a,b\right)\dd \xbar,\label{eqn:DegreeFunction}\\
    D^{(k)}(a,b) &:= \int_{\Ecal^{k-1}}\left(\prod^k_{j=1} D(c_{j-1},c_j)\right)\Pcal^{\otimes (k-1)}\left(\dd \Vec{c}_{\left[1,\ldots,k-1\right]}\right),\label{eqn:kDegreeFunction}
\end{align}
where $c_0=a$ and $c_k=b$, and the measure $\Pcal^{\otimes (k-1)}\left(\dd \Vec{c}_{\left[1,\ldots,k-1\right]}\right)$ denotes the product measure $\prod^{k-1}_{i=1}\Pcal\left(\dd c_i\right)$ on $\Ecal^{k-1}$. We can interpret these probabilistically by considering $\xi^{\origin{a}}$ for a given MRCM with intensity $\lambda$. For a measurable set $B\subset\Ecal$, Mecke's formula \eqref{eq:prelim:mecke_n} then tells us that $\lambda\int_B D(a,b)\Pcal(\dd b)$ gives the expected number of neighbours of $\origin{a}$ with a mark in $B$. Similarly, Mecke's formula tells us that $\lambda^k\int_B D^{(k)}(a,b)\Pcal(\dd b)$ gives the expected number of paths of length $k$ starting at $\origin{a}$ and ending at a vertex with a mark in $B$. 

\paragraph{Assumption D.} We make the following assumptions about $D(a,b)$ and $D^{(k)}(a,b)$:
\begin{enumerate}[label=\textbf{(D.\arabic*)}]
\item \label{Assump:BoundExpectedDegree} \emph{``Every mark has bounded expected degree with every other mark"}
    \begin{equation}
        \esssup_{a,b\in\Ecal}D(a,b) < \infty, \label{eqn:supsupbound}
    \end{equation}

\item \label{Assump:AllReachablebySome} \emph{``Some mark can be connected to every other mark in exactly $k$ steps for some $k$"}
    \begin{equation}
        \esssup_{a\in\Ecal}\sup_{k\geq 1}\essinf_{b\in\Ecal}D^{(k)}(a,b) >0.\label{eqn:supinfbound}
    \end{equation}
\end{enumerate}

Assumption~\ref{Assump:BoundExpectedDegree} is required by each of our theorems. It is proven in \cite[Proposition~2.1]{DicHey2022triangle} that \ref{Assump:BoundExpectedDegree} implies $\lambda^{(p)}_T = \lambda^{(1)}_T$ for all $p\in\left[1,\infty\right]$. Therefore, when we call upon \ref{Assump:BoundExpectedDegree} we are also be able to use $\lambda_T$ to denote this value. We can view \ref{Assump:AllReachablebySome} as a strong form of irreducibility. \ref{Assump:AllReachablebySome} is equivalent to there existing some integer $k\geq 1$ and some $\Pcal$-positive set of marks that has a positive expected number of paths of some length $k$ to $\Pcal$-almost every mark, and so there is a positive probability that such a path exists. This differs from actually being irreducible (for example like in \cite{bollobas2007phase}) in that under \ref{Assump:AllReachablebySome} there is a single $k$ that must connect to all other marks.

At first glance these conditions may appear to be very strong: \ref{Assump:BoundExpectedDegree} in particular excludes many interesting models. Nevertheless it is worth noting that even the equality of different `critical intensities' is not guaranteed without imposing some conditions. In \cite{JahnelLuchtrath2023existence} the percolation critical intensity $\lambda_c$ is compared to an annulus-crossing critical intensity, denoted $\widehat{\lambda}_c$ (introduced in \cite{gouere2008subcritical}). They show that for a class of `soft Boolean' models $\widehat{\lambda}_c=0$, while \cite{GraLucMor2021Percolation} shows that for these models $\lambda_c>0$. Note that soft Boolean models do not satisfy \ref{Assump:BoundExpectedDegree}, and therefore fall outside the scope of this paper.

The two assumptions above are sufficient to derive the mean-field bounds on $\gamma$ and $\beta$, but to show that the mean-field bounds are in fact attained in the mean-field regime we also make use of a much stronger assumption called the \emph{Triangle Condition}.  The \emph{triangle diagram} is defined using the three-fold convolution of this two-point function: given $\lambda\geq 0$ let
\begin{equation}
\label{eqn:triangle}
    \triangle_\lambda := \lambda^2 \esssup_{x,y\in\X}\int \tlam(x,u)\tlam(u,v)\tlam(v,y)\nu^{\otimes 2}\left(\dd u, \dd v\right).
\end{equation}
To formulate the following assumption, we introduce the following notation:
\begin{align}
    \Ical_{\lambda,a} &:= \left(\sup_{k\geq 1}\left(\frac{\lambda}{1+\lambda}\right)^k\essinf_{b\in\Ecal}D^{(k)}(a,b)\right)^{-1},\\
    \Jcal_{\lambda,a} &:= \left(\sup_{k\geq 1}\left(\frac{\lambda}{2\left(1+\lambda\right)}\right)^k\essinf_{b\in\Ecal}D^{(k)}(a,b)\right)^{-1},\\
    \cbar{\lambda} &:= 1+\lambda\esssup_{a,b,c\in\Ecal}D(a,b)\Jcal_{\lambda,c},\\
    \ConstantTriangle &:= \min\left\{\frac{1}{\left(1+\lambda_T\esssup_{a,b,c\in\Ecal}D(a,b)\Ical_{\lambda_T,c}\right)^2},\right.\nonumber\\
    &\hspace{5cm}\left.\frac{1}{\cbar{\lambda_T}}\frac{\lambda_T^2\left(\essinf_{a\in\Ecal} \int D(a,b)\Pcal(\dd b)\right)^2}{1 + 2\lambda_T\esssup_{a\in\Ecal}\int D(a,b)\Pcal(\dd b)}\right\}.
\end{align}
If $\sup_{k\geq 1}\left(\frac{\lambda}{1+\lambda}\right)^k\essinf_{b\in\Ecal}D^{(k)}(a,b)=0$ then we define $\Ical_{\lambda,a}=\infty$, and similarly for $\Jcal_{\lambda,a}$.

\paragraph{Triangle Condition}
\begin{enumerate}[label=\textbf{(T)}]
    \item The triangle diagram is bounded above at criticality:
    \begin{equation}
    \triangle_{\lambda_T} < \ConstantTriangle.
\end{equation}\label{TriangleCondition_Assumption}
\end{enumerate}

This assumption is a critical ingredient for proofs of the mean-field behaviours, without it only one-sided bounds can be found for the critical exponents (and none for the critical exponent $\delta$). We do not expect the exact value of $\ConstantTriangle$ to be of special importance. It is simply a bound sufficient for our proof technique to work. The notion of triangle diagram and triangle condition was originally introduced in Percolation by Aizenman and Newman in \cite{AizNew84} as an unverified condition. It came as a analogous quantity of the ``bubble diagram" introduced also by Aizenman but for the Ising model in \cite{aizenising}. In \cite[Theorem~2.5]{DicHey2022triangle} it is shown for a wide class of families of MRCMs that $\triangle_{\lambda_T}\to 0$ as $d \to \infty$, while \cite{HeyHofLasMat19} previously showed this for plain  Random Connection Models (where $\Ecal$ is a singleton). We can then interpret \ref{TriangleCondition_Assumption} as saying ``for $d$ sufficiently large" for these families of models.

It is worth remarking that having a large dimension is not the only mechanism by which the triangle condition can be made to hold. In \cite{HeyHofLasMat19} they also show that for the plain RCM in dimensions $d>6$ the triangle condition also holds for sufficiently \emph{spread-out} models. While \cite{DicHey2022triangle} only poses the question in the ``$d$ sufficiently large" case for MRCMs, the argument transfers across to the ``sufficiently spread-out" case in essentially the same way as it did for \cite{HeyHofLasMat19}.

It is important to know that in this paper we don't contribute in any way to Assumption \ref{TriangleCondition_Assumption}, we just make use of it in order to prove critical behaviours.

Note that when we call upon $\ConstantTriangle$ (that is, when we use Assumption \ref{TriangleCondition_Assumption}), we are also assuming \ref{Assump:BoundExpectedDegree} and therefore the intensity $\lambda_T$ is well-defined.

Observe that \ref{TriangleCondition_Assumption} implies \ref{Assump:AllReachablebySome}, since otherwise $\cbar{\lambda}\equiv\infty$ and $\ConstantTriangle=0$. Nevertheless, for the sake of clarity we list both \ref{Assump:AllReachablebySome} and \ref{TriangleCondition_Assumption} in our theorems because we use the former at an intermediate point of the argument independently of the latter.

Our first two theorems relate to the behaviour of the susceptibility function, $\chi_\lambda$. First we prove a lower bound that holds in relatively high generality.
\begin{theorem}[Susceptibility Mean-Field Bound]\label{thm:Susceptibility Mean-Field Bound}
If Assumption \ref{Assump:BoundExpectedDegree} holds, then there exists $C>0$ such that
\begin{equation}
    \norm*{\chi_\lambda}_p \geq C\left(\lambda_T-\lambda\right)^{-1}
\end{equation}
for $\lambda<\lambda_T$ and for all $p\in\left[1,\infty\right]$.
\end{theorem}

Our second result for the susceptibility function shows that under additional conditions there is a matching upper bound. This proves that the $\gamma$ exponent exists and takes its mean-field value in this regime.
\begin{theorem}[Susceptibility Mean-Field Behaviour]\label{thm:Susceptibility Mean-Field Behaviour}
If Assumptions \ref{Assump:BoundExpectedDegree}, \ref{Assump:AllReachablebySome}, and \ref{TriangleCondition_Assumption} all hold, then there exist $0<C\leq C'<\infty$ such that
\begin{equation}
    C\left(\lambda_T-\lambda\right)^{-1} \leq \norm*{\chi_\lambda}_p \leq C'\left(\lambda_T-\lambda\right)^{-1}
\end{equation}
for $\lambda<\lambda_T$ and for all $p\in\left[1,\infty\right]$. That is, the critical exponent $\gamma =1$.
\end{theorem}

We also pursue analogous results for the percolation function, $\theta_\lambda$. First we prove a lower bound that holds in relatively high generality - albeit not a quite as high as in Theorem~\ref{thm:Susceptibility Mean-Field Bound}.
\begin{theorem}[Percolation Mean-Field Bound]\label{thm:Percolation Mean-Field Bound}
If Assumptions \ref{Assump:BoundExpectedDegree} and \ref{Assump:AllReachablebySome} both hold, then there exist $\varepsilon>0$ and $C>0$ such that
\begin{equation}
    \norm*{\theta_\lambda}_p \geq C\left(\lambda - \lambda_T\right)_+
\end{equation}
for $\lambda < \lambda_T + \epsilon$ and for all $p\in\left[1,\infty\right]$.
\end{theorem}
This result is particularly important because it allows us to immediately show the equality of the critical intensities associated with the susceptibility and percolation functions.
\begin{corollary}[Sharpness of Phase Transition]\label{thm:Sharpness of Phase Transition}
If Assumptions \ref{Assump:BoundExpectedDegree} and \ref{Assump:AllReachablebySome} both hold, then $\lambda_c = \lambda_T$.
\end{corollary}

\begin{proof}
If $\theta_\lambda(a)>0$, then clearly $\chi_\lambda(a)=\infty$. This implies that $\lambda_c \geq \lambda_T^{(p)} = \lambda_T$ for all $p\in\left[1,\infty\right]$. Theorem~\ref{thm:Percolation Mean-Field Bound} then implies $\lambda_c \leq \lambda_T$, proving the result.
\end{proof}

The following result complements Theorem~\ref{thm:Percolation Mean-Field Bound} by proving a matching upper bound under an additional condition. In particular this proves that the $\beta$ exponent exists and takes its mean-field value in this regime.
\begin{theorem}[Percolation Mean-Field Behaviour]\label{thm:Percolation Mean-Field Behaviour}
If Assumptions \ref{Assump:BoundExpectedDegree}, \ref{Assump:AllReachablebySome}, and \ref{TriangleCondition_Assumption} all hold, then there exist $\epsilon > 0$ and $0<C\leq C'<\infty$ such that
\begin{equation}
    C\left(\lambda - \lambda_c\right)_+ \leq \norm*{\theta_\lambda}_p \leq C'\left(\lambda - \lambda_c\right)_+
\end{equation}
for $\lambda < \lambda_c + \epsilon$ and for all $p\in\left[1,\infty\right]$. That is, the critical exponent $\beta = 1$.
\end{theorem}

The following result gives conditions under which the exponent $\delta$ exists and takes its mean-field value. In contrast to the previous results, our proof method requires all the conditions to hold to prove either of the bounds. For specific classes of MRCM the lower bound has been proven to hold in a higher generality (see \cite{DewMui}), but our result is more general in the sense that it can be applied to a larger class of MRCMs.
\begin{theorem}[Cluster Tail Mean-Field Behaviour]\label{thm:Cluster Tail Mean-Field Behaviour}
If Assumptions \ref{Assump:BoundExpectedDegree}, \ref{Assump:AllReachablebySome}, and \ref{TriangleCondition_Assumption} all hold, then there exist $0<C\leq C'<\infty$ such that for all $n\in\N$ and $\Pcal$-almost every $a\in\Ecal$,
\begin{equation}
    Cn^{-\frac{1}{2}} \leq \mathbb{P}_{\lambda_c}\left(\abs*{\C\left(\origin{a}\right)}\geq n\right) \leq C'n^{-\frac{1}{2}}.
\end{equation}
That is, the critical exponent $\delta=2$.
\end{theorem}

\paragraph{Discussion.}
The theorems above have been labelled as `bound' and `behaviour' results. The exponents in the `bound' results are expected to hold in very high generality. Indeed, for Boolean disc models the bound of Theorem~\ref{thm:Percolation Mean-Field Bound} has already been proven beyond the Assumption~\ref{Assump:BoundExpectedDegree} that we require (for example see \cite{DumRaoTas20,DewMui}). The novelty here is that our results can be applied to many models for which the presence of an edge has its own randomness and is not merely a function of the point configuration (examples of such models are presented in Section~\ref{sec:Examples}). While these Boolean disc models do satisfy Assumption~\ref{Assump:AllReachablebySome}, we also don't expect this assumption to be necessary. 

The `behaviour' results state that the `bound' exponents are the correct exponents and are expected to hold in less generality. For the plain RCM (where $\Ecal$ is a singleton) \cite{HeyHofLasMat19} uses the finiteness of a `triangle' to prove the exponent $\gamma=1$ and proves that it is indeed finite for sufficiently high dimensions. Here we prove for MRCMs that if a triangle condition \ref{TriangleCondition_Assumption} holds (for which \cite{DicHey2022triangle} provides a similar `high dimensional' interpretation), then the critical exponents $\gamma$, $\beta$ and $\delta$ all take their mean-field values. In contrast to \ref{Assump:BoundExpectedDegree} and \ref{Assump:AllReachablebySome}, we expect that some sort of triangle condition is indeed necessary for this mean-field behaviour. That said, we certainly don't give significance to the value $\ConstantTriangle$, and the exact form of $\trilam$ will cause issues if one wanted to consider models outside the scope of \ref{Assump:BoundExpectedDegree}. 

It is worth remarking on how our triangle condition \ref{TriangleCondition_Assumption} differs from other triangle conditions for other models. While we require that the triangle diagram is smaller than $\ConstantTriangle$ at criticality, for Bernoulli bond percolation on certain graphs it is usually only necessary to require that (the analogue of) $\triangle_{\lambda_T}$ is finite. In \cite{AizBar91} \emph{ultraviolet regularization} is used to make the jump from finiteness to ``smallness," while \cite{kozma2011triangle} uses an alternative operator theory argument. The ultraviolet regularization approach becomes troublesome for MRCMs because each pair of marks could have a different scale on which connections occur. If $\trilam<\infty$ then we may be able to find a suitable rescaling for each pair of marks $a,b$, but we are not guaranteed that this would be suitable for all pairs of marks. A fundamental ingredient in the operator theory approach is that the two-point function on vertex-transitive graphs is positive - that is, when treated as a symmetric infinite matrix the spectrum is non-negative. In \cite{DicHey2022triangle} the two-point function of the MRCM, $\tlam$, is treated as a self-adjoint linear operator on a Hilbert space - we outline part of this in Section~\ref{sec:susceptibilityProofs}. There exist quite reasonable examples of MRCM models (even non-marked RCM models) for which this operator is \emph{not} positive. For example, in the case of the Boolean disc model with fixed radii (also called continuum percolation, Gilbert model, etc.) the two-point operator is clearly not positive at $\lambda=0$ (where it is explicit and we can evaluate the Fourier transform to evaluate its spectrum) and it is differentiable in a neighbourhood of $\lambda=0$. Therefore there exists $\lambda>0$ for which the two-point operator is not positive in this model. The lace expansion argument (see \cite{HeyHofLasMat19,DicHey2022triangle}) derives an Ornstein-Zernike equation which can also be used to show that this non-positivity extends to the critical intensity if the spatial dimension is sufficiently high.

The proof techniques employed in this paper largely build upon techniques previously used in \cite{AizNew84,AizBar91} for bond percolation. Here we adapt them for the continuum by making heavy use of Mecke's formula. For the percolation function and the critical cluster tail probabilities, this includes using a magnetization function argument. For the susceptibility arguments, the step to a continuum had already been made in \cite{HeyHofLasMat19} and their treatment of the plain RCM. The presence of marks in our model adds complications though, and in our proof of the lower bound in Section~\ref{sec:susceptibilitylowerbound} we utilise an operator formulation like that used in \cite{DicHey2022triangle} to manage this extra difficulty. It is not clear that a similar approach works for the upper bound, and so in Section~\ref{sec:susceptibilityupperbound} we instead work more directly with the functions. 

One natural concern with taking the arguments for classical models and trying to apply them in our context could be that vertices with different marks could a-priori behave \emph{very} differently. In our arguments we mitigate this by bounding the essential supremum (over $\Ecal$) of the susceptibility and magnetization functions by a factor times their pointwise values - see lemmas \ref{thm:inf-pointwise bound susceptibility} and \ref{lem:inf-pointwise_bound}. With Assumptions~\ref{Assump:BoundExpectedDegree} and \ref{Assump:AllReachablebySome} this allows us to show that up to a constant factor the essential supremum and essential infimum of these functions behave the same. This is sufficient control for our arguments to work, and is the primary reason the Assumption \ref{Assump:AllReachablebySome} appears.

A powerful technique that we do not use in this paper is the manipulation of randomised algorithms and the OSSS inequality. This was central to the proof in \cite{DumRaoTas20} to show the sharpness of the phase transition for some Boolean disc models with unbounded random radii, and in \cite{DewMui} to show mean-field bounds for these models. Furthermore, \cite{faggionato2019connection} uses these techniques to prove sharpness for some random connection models with bounded edge lengths (a stronger condition than \ref{Assump:BoundExpectedDegree}). Note that these all provide mean-field `bound' results. We are not aware of the randomised algorithm approach having been used to derive a matching bound and show a `behaviour' result. Since we wanted to arrive at results proving mean-field critical exponents, we chose to use the older techniques as a skeleton for our argument.

\subsection{Examples}
\label{sec:Examples}

Here we present some examples of models that satisfy Assumptions \ref{Assump:BoundExpectedDegree}, \ref{Assump:AllReachablebySome}, and \ref{TriangleCondition_Assumption}.

\paragraph{Boolean Disc Models.}

For $0\leq R_{\min} \leq R_{\max} <\infty$ with $R_{\max}>0$, let $\Ecal=\left[R_{\min},R_{\max}\right]$ and $\Pcal$ be any probability measure on $\Ecal$ that is not supported solely on $\left\{0\right\}$. Then let the adjacency function take the form
\begin{equation}
    \connf(\xbar;a,b) = \Id_{\left\{\abs*{\xbar}<a+b\right\}}.
\end{equation}
Observe that
\begin{equation}
     D(a,b) = \frac{\pi^{\frac{d}{2}}}{\Gamma\left(\frac{d}{2}+1\right)}\left(a+b\right)^d.
\end{equation}
This model satisfies \ref{Assump:BoundExpectedDegree} because $R_{\max}<\infty$. Since the ``radii" $a$ and $b$ are not $\Pcal$-almost surely equal to $0$, we have $\esssup_{a\in\Ecal}\essinf_{b\in\Ecal}D(a,b)>0$ and therefore \ref{Assump:AllReachablebySome} is also satisfied. Let us now think of a sequence of these models indexed by the spatial dimension $d$. In particular, now $R_{\max}= R_{\max}(d)$ and $R_{\min}= R_{\min}(d)$. Let $\Vol{d}(r)$ denote the $d$-volume of the $d$-ball with radius $r\geq0$. If $\Vol{d}(2R_{\max}(d))$ is bounded in $d$ and there exists $c>0$ such that $R_{\min}(d) > c\sqrt{d}$, then \cite{DicHey2022triangle} shows that $\triangle_{\lambda_c}\to0$ as $d\to\infty$. Therefore for each such sequence of models, the condition \ref{TriangleCondition_Assumption} holds for $d$ sufficiently large.

In the specific context of Boolean disc models, some of our results have been proven in greater generality. If the radius distribution has finite $(5d-3)$-moment, then \cite{DumRaoTas20} proves that the percolation mean-field bound holds, and \cite{DewMui} uses this to prove the susceptibility mean-field bound and the lower bound on the cluster tail probability. The sharpness of the phase transition has been proven in $d=2$ under the minimal assumption that the radius distribution has finite second moment in \cite{ahlberg2018sharpness}. The planarity of $\R^2$ was crucial for their argument which relied on Russo–Seymour–Welsh theory. In general spatial dimensions, \cite{dembin2022almost} proves `subcritical sharpness' for models where the radius distribution has tail density $r^{-d-1-\delta}$ with $\delta>0$ (except for at most countably many values of $\delta$).

\paragraph{Marked Multivariate Gaussian Model.}
Let $\Sigma\colon\Ecal^2\to\R^{d\times d}$ be a measurable map where for every $a,b\in\Ecal$, $\Sigma\left(a,b\right)$ is itself a symmetric positive definite covariance matrix. We further require that there exist $\Sigma_{\max}<\infty$,  $\Sigma_{\min} > 0$, and $\mathcal{A}>0$ such that the set of eigenvalues $\sigma(\Sigma(a,b)) \subset \left[\Sigma_{\min},\Sigma_{\max}\right]$ and $\mathcal{A}^2 \leq \left(2\pi\right)^d\det \Sigma(a,b)$ for all $a,b\in\Ecal$. Then let the adjacency function be given by
\begin{equation}
    \connf(\xbar;a,b) = \mathcal{A}\left(2\pi\right)^{-d/2}\left(\det \Sigma(a,b)\right)^{-1/2}\exp\left(-\frac{1}{2}\xbar^{\intercal}\Sigma(a,b)^{-1}\xbar\right).
\end{equation}
For this adjacency function, $ D(a,b) = \mathcal{A}$ and is independent of the marks, and therefore \ref{Assump:BoundExpectedDegree} and \ref{Assump:AllReachablebySome} both hold. Let us once again think of a sequence on these models indexed by the spatial dimension $d$. If $\Sigma_{\max}=\Sigma_{\max}(d)$ is uniformly bounded above and $\limsup_{d\to\infty} \mathcal{A}(d)^2\left(4\pi \Sigma_{\min}(d)\right)^{-d/2} = 0$, then \cite{DicHey2022triangle} proves that $\triangle_{\lambda_c}\to0$ as $d\to\infty$. Therefore for each such sequence of models, the condition \ref{TriangleCondition_Assumption} holds for $d$ sufficiently large.

\paragraph{Factorisable Model.}
Let $\psi\colon \Rd\to \left[0,1\right]$ be a measurable function such that $\psi(\xbar)=\psi(-\xbar)$ and $\int \psi(\xbar)\dd \xbar\in\left(0,\infty\right)$. Also let $K\colon\Ecal^2\to\left[0,1\right]$ be measurable, symmetric, and satisfy $\esssup_{c_0\in\Ecal}\sup_{k\geq 1}\essinf_{c_k\in\Ecal}\int \left(\prod^k_{j=1}K(c_{j-1},c_j)\right)\Pcal^{\otimes (k-1)}\left(\dd \vec{c}_{\left[1,\ldots,k-1\right]}\right)>0$. Then let the adjacency function be given by
\begin{equation}
    \connf(\xbar;a,b) = \psi(\xbar)K(a,b).
\end{equation}
We clearly have
\begin{equation}
    \esssup_{a,b\in\Ecal}D(a,b) \leq \int \psi(\xbar)\dd \xbar<\infty
\end{equation}
because $K(a,b)\in\left[0,1\right]$, and so \ref{Assump:BoundExpectedDegree} holds. Note that
\begin{multline}
    \esssup_{a\in\Ecal}\sup_{k\geq 1}\essinf_{b\in\Ecal}D^{(k)}(a,b)\\= \esssup_{c_0\in\Ecal}\sup_{k\geq 1}\left(\int \psi(\xbar)\dd \xbar\right)^k\essinf_{c_k\in\Ecal}\int \left(\prod^k_{j=1}K(c_{j-1},c_j)\right)\Pcal^{\otimes(k-1)}\left(\dd \vec{c}_{\left[1,\ldots,k-1\right]}\right).
\end{multline}
Our last condition on $K$ gives that this is positive (the factors of $\int \psi(\xbar)\dd \xbar$ are all positive and so cannot change the positivity of the whole) and so \ref{Assump:AllReachablebySome} holds. Let us index $\psi=\psi_d$ by the spatial dimension but let $K$ be $d$-independent. Let us also construct the self-adjoint integral operator $\Kcal\colon L^2(\Ecal) \to L^2(\Ecal)$ by
\begin{equation}
    \left(\Kcal f\right)(a) = \int K(a,b)f(b)\Pcal\left(\dd b\right)
\end{equation}
for all $f\in L^2(\Ecal)$ and $\Pcal$-almost every $a\in\Ecal$, and have $\sigma(\Kcal)\subset\R$ denote the spectrum of $\Kcal$. If $\sup \sigma(\Kcal) \geq \abs*{\inf\sigma(\Kcal)}$ and $\psi_d$ satisfies the conditions of a finite-variance model given in \cite{HeyHofLasMat19}, then \cite{DicHey2022triangle} proves that $\triangle_{\lambda_c}\to0$ as $d\to\infty$. Therefore for each such sequence of models, the condition \ref{TriangleCondition_Assumption} holds for $d$ sufficiently large.

To perhaps make the above clearer, let us consider the concrete example where $\Ecal=\left\{1,2,3\right\}$, $\Pcal$ gives mass $1/3$ to each of the three singletons, and
\begin{equation}
    \left\{K(i,j)\right\}_{i,j=1,2,3} = 
    \begin{pmatrix}
        1 & 1 & 0\\
        1 & 0 & 1\\
        0 & 1 & 0
    \end{pmatrix}. \label{eqn:Kmatrix}
\end{equation}
This model can then be interpreted as a plain RCM with adjacency function $\psi$, where we independently assign each vertex a mark $i\in\left\{1,2,3\right\}$ with probability $1/3$, and then delete edges that are between two $2$ vertices, between two $3$ vertices, and between $1$ and $3$ vertices. As in the more general case \ref{Assump:BoundExpectedDegree} clearly holds, and we can show \ref{Assump:AllReachablebySome} holds by considering $k=4$:
\begin{equation}
    \left\{\frac{1}{3^3}\sum^3_{c_1,c_2,c_3=1} \left(\prod^4_{j=1}K(c_{j-1},c_j)\right)\right\}_{c_0,c_4=1,2,3} = \frac{1}{27}
    \begin{pmatrix}
        6 & 4 & 3\\
        4 & 5 & 1\\
        3 & 1 & 2
    \end{pmatrix}.
\end{equation}
Also note that the operator $\mathcal{K}$ can be represented by the matrix $\left\{\frac{1}{3}K(i,j)\right\}_{i,j=1,2,3}$ acting on vectors in $\R^3$, and so the Perron-Frobenius Theorem guarantees that $\sup \sigma(\Kcal) \geq \abs*{\inf\sigma(\Kcal)}$. If $\psi_d$ is sufficiently nicely behaved, then we have $\triangle_{\lambda_c}\to 0$ as $d\to\infty$.

\subsection{Overview of the Paper}
We conclude Section~\ref{sec:Introduction} by identifying some probabilistic tools that we use frequently in our arguments. Then Section~\ref{sec:susceptibilityProofs} contains the proofs of our results regarding the susceptibility function (Theorems~\ref{thm:Susceptibility Mean-Field Bound} and \ref{thm:Susceptibility Mean-Field Behaviour}). In Section~\ref{sec:criticalintensities} we also present some results on the different critical intensities that we talk about in this paper. Sections~\ref{sec:susceptibilitylowerbound} and \ref{sec:susceptibilityupperbound} then respectively prove the lower and upper bounds on the susceptibility function by various uses of the Margulis-Russo formula. For the percolation and cluster tail results we make use of the magnetization function. In Section~\ref{sec:MagnetizationBounds} we derive upper and lower bounds on this magnetization function by construction and manipulations of partial differential equations satisfied by the magnetization. We use Mecke's formula and the Weierstrass M-test to prove that the magnetization function is analytic on the required domain and therefore that the required partial derivatives exist. In Section~\ref{sec:PercolationProofs} we first use the behaviour of the magnetisation (as the magnetization parameter approaches zero) to derive a pointwise lower bound for the percolation function which we then relate to the $L^p$-norms stated in Theorem~\ref{thm:Percolation Mean-Field Bound}. In Section~\ref{sec:percolation_upper_bound} we use our magnetization upper bound and the extrapolation principle to complete the proof of Theorem~\ref{thm:Percolation Mean-Field Behaviour}. Finally in Section~\ref{sec:ClusterTailProof} we use our magnetization bounds to derive the bounds on the cluster tail probabilities required of Theorem~\ref{thm:Cluster Tail Mean-Field Behaviour}.

\subsection{Preliminaries}\label{sec:Prelims}

Before we embark on proving the main results of this paper let us introduce some useful probability tools and theorems that will reappear multiple times. For a given set $E$ we denote $\mathbf{N}(E)$ the set of all at most countably infinite subsets of $E$. In our case in particular the space of configurations is $\mathbf{N} := \mathbf{N}(\X^{[2]} \times [0,1])$, where $\X^{[2]}$ are the subsets of $\X$ containing exactly two elements (these represents the possible edges), the interval $[0,1]$ is what contains the information about the presence or not of the edge in question. Notice that it is also possible to replace $\X^{[2]}$ by $(\X \times [0,1]^\N)^{[2]}$ if one wants to includes thinned configurations, and all the results below will still hold. For a better understanding on how a configuration is encoded as well as the $\sigma$-algebra associated the reader is referred to the formal construction of the MRCM in \cite[Section 3.1]{DicHey2022triangle}. For more details about the probabilistic lemmas below the reader is referred to \cite[Section 3.2]{DicHey2022triangle}

\paragraph{Mecke's Formula.}
Given $m\in\N$ and a measurable function $f\colon\mathbf{N} \times \X^m \to \R_{\geq 0}$, the Mecke equation for $\xi$ states that 
\begin{equation}
    \E_\lambda \left[ \sum_{\vec x \in \eta^{(m)}} f(\xi, \vec{x})\right] = \lambda^m \int
				\E_\lambda\left[ f\left(\xi^{x_1, \ldots, x_m}, \vec x\right)\right] \nu^{\otimes m}\left(\dd \vec{x}\right),  \label{eq:prelim:mecke_n}
\end{equation}
where $\vec x=(x_1,\ldots,x_m)$, $\eta^{(m)}=\{(x_1,\ldots,x_m)\colon x_i \in \eta, x_i \neq x_j \text{ for } i \neq j\}$, and $\nu^{\otimes m}$ is the $m$-product measure of $\nu$ on $\X^m$. A proof and discussion of Mecke's formula can be found in \cite[Chapter~4]{LasPen17}.

    When we introduced the notation $D(a,b)$ and $D^{(k)}(a,b)$ in \eqref{eqn:DegreeFunction} and \eqref{eqn:kDegreeFunction} we remarked that Mecke's formula allowed us to interpret them as the expected number of $b$-neighbours of $\origin{a}$, and as the expected number of length $k$ paths starting at $\origin{a}$ and ending at some $b$-vertex respectively. Specifically, the interpretation of $D(a,b)$ follows from applying Mecke's formula with $m=1$ and $f(\xi,x)=\Id_{\left\{x\sim \origin{a}\right\}}\Id_{\left\{x\in\Rd\times B\right\}}$ for any measurable $B\subset\Ecal$:
    \begin{equation}
        \E_\lambda \left[ \sum_{ x \in \eta} \Id_{\left\{x\sim \origin{a}\right\}}\Id_{\left\{x\in\Rd\times B\right\}}\right] = \lambda\int_{B}\left(\int_{\Rd}\connf\left(\xbar;a,b\right)\dd \xbar\right)\Pcal\left(\dd b\right) = \lambda\int_B D(a,b)\Pcal(\dd b).
    \end{equation}
    Similarly, the interpretation of $D^{(k)}(a,b)$ follows from applying Mecke's formula with $m=k$ and $f\left(\xi,\left(x_1,\ldots,x_k\right)\right)=\Id_{\left\{\origin{a}\sim x_1 \sim \ldots \sim x_k\right\}}\Id_{\left\{x_k\in\Rd\times B\right\}}$ for any measurable $B\subset\Ecal$:
    \begin{align}
        &\E_\lambda \left[ \sum_{ x \in \eta^{(k)}} \mathbb{1}_{\left\{\origin{a}\sim x_1 \sim \ldots \sim x_k\right\}}\Id_{\left\{x_k\in\Rd\times B\right\}}\right]\nonumber\\
        &\quad= \lambda^k \int_B \int_{\Ecal^{k-1}}\left(\int_{\left(\Rd\right)^k}\prod^k_{j=1}\connf(\xbar_j-\xbar_{j-1};c_j,c_{j-1})\dd \vec{\xbar}_{\left[1,\ldots,k\right]}\right) \nonumber\\
        &\hspace{8cm}\Pcal^{\otimes(k-1)}\left(\dd \vec{c}_{\left[1,\ldots,k-1\right]}\right)\Pcal\left(\dd c_k\right)\nonumber\\
        &\quad= \lambda^k \int_B\left( \int_{\Ecal^{k-1}} \prod^k_{j=1}D(c_j,c_{j-1})\Pcal^{\otimes(k-1)}\left(\dd \vec{c}_{\left[1,\ldots,k-1\right]}\right) \right)\Pcal\left(\dd c_k\right)\nonumber\\
        &\quad= \lambda^k\int_B D^{(k)}(a,c_k)\Pcal(\dd c_k),
    \end{align}
    where $\left(\xbar_0,c_0\right)=\origin{a}$.

    Another particular application of this is to take $m=1$ and let $f\left(\xi,x\right) = \Id_{\left\{\conn{\origin{a}}{x}{\xi^\origin{a}}\right\}}$. This tells us that
    \begin{equation}
        \chi_\lambda(a) = \E_\lambda\left[\abs*{\C\left(\origin{a}\right)}\right] = 1 + \E_\lambda\left[\sum_{x\in\eta}\Id_{\left\{\conn{\origin{a}}{x}{\xi^\origin{a}}\right\}}\right]\\
        =1 + \lambda\int\tlam\left(\origin{a},x\right)\nu\left(\dd x\right).
    \end{equation}
    Here the term $1$ accounts for $\origin{a}\in\C\left(\origin{a}\right)$, while the second term accounts for the vertices in $\xi$ that are connected to $\origin{a}$.

\paragraph{Margulis-Russo formula.} A classical result in percolation theory is the Margulis-Russo formula which allows us to get the derivative of functions of the form $\lambda \mapsto \el[f(\xi)]$ for a fixed configuration $\xi$ and function $f$. More precisely let $\Lambda \subset \X$ be $\nu$-finite, $\zeta \in \mathbf{N}$ and define 
\begin{align}
\label{restriction_conf}
    \zeta_\Lambda := \{(\{(x,v),(y,w)\}, u) \in \zeta \colon \{x,y\} \subseteq \Lambda \}
\end{align}
we call $\zeta_\Lambda$ the restriction of $\zeta$ to $\Lambda$.  We say that $f \colon \mathbf{N} \rightarrow \mathbb{R}$ \textit{lives} on $\Lambda$ if $f(\zeta) = f(\zeta_\Lambda)$ for every $\zeta \in \mathbf{N}$. Assume that there exists a $\nu$-finite $\Lambda \subset \X$ such that $f$ lives on $\Lambda$. Moreover, assume that there exists $\lambda_1 > 0$ such that $\mathbb{E}_{\lambda_1}[|f(\xi)|] < \infty$. Note that the $\nu$-finiteness of $\Lambda$ implies that $\mathbb{P}_{t}\left(\eta_\Lambda=\emptyset\right)>0$ for all $t>0$, and therefore the superposition theorem of Poisson point processes implies that $\mathbb{E}_{\lambda}[|f(\xi)|] < \infty$ for all $\lambda\leq\lambda_1$ (see \cite[Exercise~3.8]{LasPen17}). Then the Margulis-Russo's formula states that, for all $\lambda \leq \lambda_1$, 
\begin{align}
\frac{\partial}{\partial \lambda}\el[f(\xi)]= \int_\Lambda \el[f(\xi^x) - f(\xi)] \nu(\mathrm{d}x).
\end{align}
More details on this result can be found in \cite{chebunin2024uniqueness}, and the $\X=\Rd$ case can be found in \cite[Theorem~3.2]{LasZie17}.

\paragraph{BK inequality.} Before introducing the BK inequality, another  classical result in percolation that still holds for our model, we need to introduce the notion of increasing sets, increasing events, events living on a set and disjoint occurrence. We denote the $\sigma$-algebra on $\mathbf{N}$ by $\mathcal{N}$. We call a set $E \subset \mathbf{N}$ increasing if $\mu \in E$ implies $\nu \in E$ for each $\nu \in \mathbf{N}$ with $\mu \subseteq \nu$. Let $(\mathbb{Y}_1,\mathcal{Y}_1), (\mathbb{Y}_2,\mathcal{Y}_2)$ be two measurable spaces. We say that a set $E_i \subset \mathbf{N} \times \mathbb{Y}_i$ is \textit{increasing} if $E_i^z := \{\mu \in \mathbf{N} \colon (\mu,z) \in E_i\}$ is increasing for each $z \in \mathbb{Y}_i$.

Given a  Borel  set $\Lambda \in  \mathcal{B}(\Rd)$ and $\mu \in \mathbf{N}$, we define $\mu_\Lambda$, the restriction of $\mu$ to all edges completely contained in $(\Lambda \times \Ecal)$, analogously to \ref{restriction_conf}. Furthermore let $\mathcal{R}$ denote the ring of all finite unions of half-open $d$-dimensional rectangles in $\Rd$ with rational coordinates. Then for increasing $E_i \in \mathbf{N} \otimes \mathbb{Y}_i$, we define
\begin{multline}
E_1 \circ E_2 := \{(\mu,z_1,z_2) \in \mathbf{N} \times \mathbb{Y}_1 \times \mathbb{Y}_2 \\\colon \exists K_1,K_2 \in \mathcal{R} \text { s.t. }
 K_1 \cap K_2 = \emptyset,(\mu_{K_1} ,z_1) \in E_1,(\mu_{K_2} ,z_2) \in E_2\}
\end{multline}
and say that $E_1$ and $E_2$ \textit{occurs disjointly}, basically the vertices and edges that make $E_1$ happen are disjoint of the ones that make $E_2$ happen.
A set $E_i \in \mathbf{N} \otimes \mathbb{Y}_i$ \textit{lives} on $\Lambda$ if $\mathds{1}_{E_i} (\mu,z) = \mathds{1}_{E_i} (\mu_\Lambda,z)$ for each $(\mu,z) \in \mathbf{N} \times \mathbb{Y}_i$. We consider random elements $W_1, W_2$ of $\mathbb{Y}_1$ and $\mathbb{Y}_2$, respectively, and assume that $\xi$, $W_1$ and $W_2$ are independent. The following theorem is proven in \cite{HeyHofLasMat19} for $\X=\Rd$, but generalising to $\X=\Rd\times \Ecal$ adds no further difficulty.

\begin{theorem}(BK inequality)
Let $E_1 \in \mathbf{N} \otimes \mathbb{Y}_1$ and $E_2 \in \mathbf{N} \otimes \mathbb{Y}_2$ be increasing events that live on $\Lambda \times \Ecal$ for some bounded set $\Lambda \in \mathcal{B}(\Rd)$. Then
\begin{align}
    \pl((\xi,W_1,W_2) \in E_1 \circ E_2) \leq \pl((\xi,W_1) \in E_1)\pl((\xi,W_2) \in E_2)
\end{align}
\end{theorem}

\paragraph{FKG inequality.}  The last classical result in percolation that we use is the FKG inequality, which in contrast with the BK inequality give us a lower bound to some events instead of and upper bound. Given two increasing and integrable function $f, g$ on $\mathbf{N}$, we have
\begin{align}
    \el[f(\xi)g(\xi)] \geq \el[\el[f(\xi)g(\xi)\mid \eta]]  \geq \el[\el[f(\xi)\mid\eta]\el[g(\xi)\mid\eta]] \geq \el[f(\xi)]\el[g(\xi)]
\end{align}

We saw previously it is very useful to consider augmentation of a given realisation, like for example to consider the cluster of a point. In the same way sometimes it is also useful to delete some vertices and their respective edges. For that we introduce the notion of thinnings.

\paragraph{Thinnings.}
At various points of our arguments we use thinning events. These are used to describe connectivity properties by asking whether a connection still exists if vertices are removed in a particular way. We give a brief description here, but for a full and complete description we refer the reader to \cite{DicHey2022triangle,HeyHofLasMat19}.

Let $u,x\in\X$ and $A\subset\X$ be locally finite.
\begin{itemize}
    \item Set
\begin{equation}
    \overline{\connf}(A,x) := \prod_{y\in A}\left(1-\connf(y,x)\right).
\end{equation}
We define an $A$-thinning of some point configuration $\eta$ by taking each $u\in\eta$ and retaining $u$ with probability $\overline{\connf}(A,u)$ independently of all other points of $\eta$. We use $\eta_{\thinn{A}}$ to denote this $A$-thinning of $\eta$.

\item Let $\left\{\xconn{u}{x}{\xi}{A}\right\}$ denote the event that both $u,x\in\eta$ and that $\conn{u}{x}{\xi}$, but that this connection does not survive an $A$-thinning of $\eta\setminus\left\{u\right\}$.

\item Let $\tlam^{A}(u,x)$ denote the probability that there exists an open path between $u$ and $x$ on an $A$-thinning of $\eta^{x}$, where $u$ is conditioned to be present. That is,
\begin{equation}
    \tlam^{A}(u,x) = \pla \left(\conn{u}{x}{\xi^{u,x}[\thinning{\eta^x}{A}\cup\{u\}]} \right). \label{eq:def:LE:offconn}
\end{equation}

\item Given $x,y \in \X$ and edge-marking $\xi$, we say $u \in \X$ is pivotal and write $u \in \text{Piv}(x,y,\xi)$ if $\{\conn{x}{y}{\xi^{x,y}}\}$ but $\{\nconn{x}{y}{\xi[\eta \setminus u]}\}$. That is, every path on $\xi^{x,y}$ connecting $x$ and $y$ uses the vertex $u$. Note that the end points $x$ and $y$ are never pivotal.

\item Given $x \in \X$ and a (possibly random) set $B \subset \X$, we extend the last definition, we say that $u \in \X$ is pivotal for the connection of $u$ to $B$ and write $u \in \text{Piv}(x,B,{\xi})$ if $ \conn{x}{B}{\xi^{x,u}}$ but $\nconn{x}{B}{\xi^{x}}$. Notice that in particular $u$ and $x$ must be connected in $\xi^{x,u}$.

\end{itemize}

\section{Susceptibility Proofs}
\label{sec:susceptibilityProofs}

In this section we prove the bounds on the susceptibility function $\chi_\lambda$. We first address some potential ambiguities in the idea of a ``critical intensity" in Section~\ref{sec:criticalintensities}. In \eqref{eqn:percolationCritical} and \eqref{eqn:susceptibilityCritical} we introduced critical intensities arising from the percolation function and the susceptibility function respectively. In the later arguments the two-point operator is used and the natural associated critical intensity for that is $\lambda_O$ defined in \eqref{eqn:operatorCritical}. We provide necessary and sufficient conditions for $\lambda_O$ to be non-trivial (that is, finite and strictly positive), and show that under Assumption~\ref{Assump:BoundExpectedDegree} the operator-based critical intensity equals the susceptibility-based critical intensity. In Section~\ref{sec:susceptibilitylowerbound} we use differential inequalities to arrive at a lower bound on the operator norm of the two-point operator for $\lambda<\lambda_O$. This holds in very high generality, and when we add in the Assumption \ref{Assump:BoundExpectedDegree} we prove the lower bound on $\norm*{\chi_\lambda}_p$ in Theorem~\ref{thm:Susceptibility Mean-Field Bound}. In Section~\ref{sec:susceptibilityupperbound} we adapt an argument from \cite{HeyHofLasMat19} to prove the mean-field behaviour in Theorem~\ref{thm:Susceptibility Mean-Field Behaviour}.

In our arguments we are making use of the following notation. Given measurable $h\colon \Ecal^2\to\R\cup\left\{\pm\infty\right\}$ and $p_1,p_2\in\left[1,\infty\right)$, define the norms
\begin{align}
    \norm*{h}_{p_1,p_2} :=& \left(\int \left(\int \abs*{h\left(a,b\right)}^{p_1} \Pcal\left(\dd b\right)\right)^\frac{p_2}{p_1} \Pcal\left(\dd a\right)\right)^\frac{1}{p_2},\\
    \label{eqn:normPtoInfinity}\norm*{h}_{p_1,\infty} :=& \esssup_{a\in\Ecal}\left(\int \abs*{h\left(a,b\right)}^{p_1} \Pcal\left(\dd b\right)\right)^\frac{1}{p_1},\\
    \norm*{h}_{\infty,\infty} :=& \esssup_{a,b\in\Ecal} \abs*{h\left(a,b\right)}.
\end{align}
We apply these to symmetric $h$, and therefore the question of which \emph{argument} of $h$ each $L^p$-norm is being applied to first does not matter.

\subsection{Critical Intensities}
\label{sec:criticalintensities}

Here we address potential ambiguities in the idea of a critical intensity.

The following lemma is a refinement of \cite[Proposition~2.1]{DicHey2022triangle}. In particular, it also produces the result that if $\norm*{ D}_{\infty,\infty}<\infty$ then $\lambda^{(p)}_T\equiv \lambda_T$ for all $p\in\left[1,\infty\right]$.

\begin{lemma}
\label{thm:equality of lambda_T}
Let $1\leq p_1 \leq p_2 \leq \infty$ and $q_1:= \frac{p_1}{p_1-1}\in\left[1,\infty\right]$. Then
\begin{equation}
    \norm*{ D}_{q_1,p_2}<\infty \implies \lambda^{(p_1)}_T = \lambda^{(p_2)}_T.
\end{equation}
\end{lemma}

\begin{proof}
First note that Jensen's inequality implies that $\norm*{\chi_\lambda}_{p_2}\geq \norm*{\chi_\lambda}_{p_1}$, and therefore $\lambda^{(p_1)}_T \geq \lambda^{(p_2)}_T$. We only need to then show the reverse inequality.

By considering the graph neighbours of $\origin{a}$ and applying Mecke's formula, we find
\begin{align}
    \E_\lambda\left[\abs*{\C(\origin{a})}\right] &\leq 1 + \lambda\int \E_\lambda\left[\abs*{\C(x)}\right]\connf(x,\origin{a})\nu\left(\dd x\right) \nonumber\\
    & = 1 + \lambda\int \E_\lambda\left[\abs*{\C(\origin{b})}\right] D(b,a)\Pcal\left(\dd b\right) \nonumber \\
    &\leq 1 + \lambda \norm*{\chi_\lambda}_{p_1} \left(\int  D(b,a)^{q_1}\Pcal\left(\dd b\right)\right)^{\frac{1}{q_1}}. 
\end{align}
In this last inequality we used H{\"o}lder's inequality. Then taking the $p_2$-norm over $a\in\Ecal$ and using the triangle inequality gives
\begin{equation}
\label{eqn:norm_relation}
    \norm*{\chi_\lambda}_{p_2} \leq 1 + \lambda \norm*{\chi_\lambda}_{p_1} \norm*{ D}_{q_1,p_2}.
\end{equation}
This inequality means that if $\norm*{ D}_{q_1,p_2}<\infty$, then $\norm*{\chi_\lambda}_{p_2}=\infty\implies \norm*{\chi_\lambda}_{p_1}=\infty$. Thus $\lambda^{(p_1)}_T \leq \lambda^{(p_2)}_T$ and the result is proven.
\end{proof}

To prove Theorem~\ref{thm:Susceptibility Mean-Field Bound}, we make use an operator formalism also used by \cite{DicHey2022triangle}. Let $L^2\left(\Ecal\right)$ denote the space of square-integrable functions on $\Ecal$ (with respect to the measure $\Pcal$). When accompanied by the inner product $\inner{f}{g} := \int \overline{f(a)}g(a)\Pcal(\dd a)$, where $f,g\in L^2(\Ecal)$ and $\overline{f(a)}$ is the complex conjugate of $f(a)$, and the corresponding norm $\norm*{f}_2 := \inner{f}{f}^{\frac{1}{2}}$, the space $L^2\left(\Ecal\right)$ constitutes a Hilbert space (in particular, a Banach space). The space of bounded linear operators on a Banach space is itself a Banach space when augmented with the \emph{operator norm} which we define by
\begin{equation}
    \OpNorm{H} := \sup_{f\in L^2(\Ecal):f\ne 0}\frac{\norm*{Hf}_2}{\norm*{f}_2},
\end{equation}
for a linear operator $H\colon L^2(\Ecal)\to L^2(\Ecal)$.

Integral operators are of particular importance to us. Let us denote
\begin{equation}
     D(a,b) = \int\connf(\xbar;a,b)\dd \xbar,\qquad  T_\lambda(a,b) = \int \tlam(\xbar;a,b)\dd \xbar.
\end{equation}
Given $ D(a,b)$ and $ T_\lambda(a,b)$, we define the integral operators $\Opconnf,\Optlam\colon L^2(\Ecal)\to L^2(\Ecal)$ as those operators that act as
\begin{equation}
    \left(\Opconnf f\right)(a) = \int  D(a,b)f(b)\Pcal(\dd b), \qquad \left(\Optlam f\right)(a) = \int  T_\lambda(a,b)f(b)\Pcal(\dd b),
\end{equation}
for $f\in L^2(\Ecal)$ and $\Pcal$-almost every $a\in\Ecal$.

For some of our calculations it is convenient to work with the critical intensity
\begin{equation}
\label{eqn:operatorCritical}
    \lambda_O = \inf\left\{\lambda>0\colon \OpNorm{\Optlam }=\infty\right\}.
\end{equation}
Observe that $ T_\lambda(a,b)$ is symmetric in $a$ and $b$ and for $\lambda<\lambda_O$ the operator $\Optlam $ is a bounded operator. Therefore $\Optlam $ a self-adjoint operator in this regime. Since $ D(a,b)$ is also symmetric, if $\Opconnf $ is a bounded operator then it is also self-adjoint. Lemma~\ref{lem:nontriviallambdaO} below gives very general conditions under which $\lambda_O$ is non-trivial (i.e. $\lambda_O\in\left(0,\infty\right)$), and Lemma~\ref{lem:lambdaO_lambdaT} then shows that under some other conditions we can equate this $\lambda_O$ with other critical intensities that we have already defined. The relation of an operator norm to a critical threshold is not new. In the context of inhomogeneous random graphs, \cite{bollobas2007phase} prove that the emergence of giant components corresponds to the operator norm of the edge-probability integral operator being strictly larger than $1$. In the context of bond percolation on connected, locally finite, transitive graphs, \cite{Hutch20l2boundedness} discusses consequences of the conjecture that the two-point matrix has bounded operator norm at the critical edge probability.

\begin{lemma}
\label{lem:nontriviallambdaO}
$\OpNorm{\Opconnf }<\infty$ if and only if $\lambda_O >0$, and $\OpNorm{\Opconnf }>0$ if $\lambda_O <\infty$. If $d\geq 2$ and $\OpNorm{\Opconnf }>0$, then $\lambda_O <\infty$.
\end{lemma}

\begin{proof}
We first show $\OpNorm{\Opconnf }=\infty$ implies $\lambda_O =0$. Clearly two vertices are connected if they are adjacent. Thus $\connf(\xbar;a,b) \leq \tlam(\xbar;a,b)$ for all $\lambda\leq 0$, $\xbar\in\Rd$, and $a,b\in\Ecal$ (in fact, equality holds at $\lambda=0$), and hence $0\leq  D(a,b) \leq  T_\lambda(a,b)$ for all $\lambda\leq 0$, and $a,b\in\Ecal$. It is shown in \cite[Lemma~3.3]{DicHey2022triangle} that if two kernel functions are ordered in this way then the operator norms of their respective integral operators inherit this ordering. That is, we have $\OpNorm{\Opconnf }\leq \OpNorm{\Optlam }$ for all $\lambda\geq 0$. If $\OpNorm{\Opconnf }=\infty$, then $\OpNorm{\Optlam }=\infty$ for all $\lambda\geq0$ and $\lambda_O=0$.

To show $\OpNorm{\Opconnf }<\infty$ implies $\lambda_O >0$, we use a `method of generations' approach and compare the cluster of a vertex with the tree produced by a sub-critical branching process starting at that vertex. For $a\in\Ecal$, let $\origin{a}$ be the root of our cluster (generation $0$). Those vertices adjacent to $\origin{a}$ (generation $1$) are distributed as a Poisson point process with intensity density $\lambda\connf(\cdot;a,\cdot)$ with respect to $\Leb\times\Pcal$. For integer $k\geq 1$, generation $k+1$ then consists of vertices adjacent to vertices in generation $k$ but no vertices from earlier generations. We also assign an arbitrary ordering to each generation. \cite[Corollary~3.7]{LasPen17} shows how a measurable ordering can be given to the whole Poisson point process, and this can be used to give an ordering to each generation because determining if a vertex in in generation $k$ is a finite procedure and there are almost surely only countably many vertices. Given some vertex $(\xbar,b)$ in generation $k$, we then call those vertices in generation $k+1$ that are neighbours of $(\xbar,b)$ but not neighbours of any earlier vertex in generation $k$ (in our ordering) to be the `immediate descendants' of $(\xbar,b)$. The immediate descendants of some vertex $(\xbar,b)$ are then distributed according to a thinned Poisson point process. This point process is also a Poisson point process, whose intensity density we can bound above by $\lambda\connf(\cdot-\xbar;b,\cdot)$. We can therefore bound $\tlam$ from above using the kernel of a `spatial' branching process. For all $\xbar\in\Rd$, and $a,b\in \Ecal$ we get
\begin{equation}
    \tlam(\xbar;a,b) \leq \connf(\xbar;a,b) + \sum^\infty_{k=1}\lambda^k\int\prod^{k+1}_{i=1}\connf(\ubar_{i}-\ubar_{i-1};c_{i-1},c_i)\dd \ubar_{[1:k]}\Pcal^{\otimes k}(\dd c_{[1:k]}),
\end{equation}
where for each integral $\ubar_0=\zerobar$, $\ubar_{k+1}=\xbar$, $c_0=a$, and $c_{k+1}=b$. Integrating $\xbar$ over $\Rd$ (and using Tonelli's Theorem to swap an infinite sum and integral) then gives us
\begin{equation}
     T_\lambda(a,b) \leq  D(a,b) + \sum^\infty_{k=1}\lambda^k\int\prod^{k+1}_{i=1} D(c_{i-1},c_i)\Pcal^{\otimes k}\left(\dd c_{[1:k]}\right),
\end{equation}
for all $a,b\in \Ecal$. It can now be seen that the right hand side is the kernel function of the integral operator $\sum^\infty_{k=1}\lambda^{k-1}\Opconnf ^k$. As above, this inequality on the kernel functions passes to an inequality on the operator norm of the associated operators. We can then use the triangle inequality and the sub-multiplicativity of the operator norm to get
\begin{equation}
    \OpNorm{\Optlam }\leq \sum^\infty_{k=1}\lambda^{k-1}\OpNorm{\Opconnf }^k.
\end{equation}
If $\OpNorm{\Opconnf }<\infty$, then it is possible to choose $\lambda>0$ sufficiently small that $\lambda\OpNorm{\Opconnf }<1$. For such a $\lambda$ the sum then converges and we have
\begin{equation}
    \OpNorm{\Optlam }\leq \OpNorm{\Opconnf }\left(1-\lambda\OpNorm{\Opconnf }\right)^{-1}<\infty.
\end{equation}
Hence $\lambda_O\geq \OpNorm{\Opconnf }^{-1}>0$.

Note that the above argument also shows that $\OpNorm{\Opconnf }=0$ implies $\lambda_O=\infty$.

We now aim to show that $\OpNorm{\Opconnf }>0$ implies $\lambda_O<\infty$ for $d\geq 2$. In \cite[Proposition~2.1]{DicHey2022triangle} it is proven that $\lambda_c\geq \lambda_O$ - this is done by comparing both with $\lambda_T^{(1)}$. Here we adapt a coarse-graining argument from \cite{Pen91} that shows that $\lambda_c<\infty$, the extra complication being the possible inhomogeneity arising from the marks. The idea is to couple the cluster in the marked RCM with a cluster in $\Z^2$ arising from the independent Bernoulli bond percolation model. A higher intensity $\lambda$ corresponds to a higher edge probability in the bond percolation model, and so when $\lambda$ is sufficiently high we are in a super-critical regime.

Suppose for contradiction that $\connf(\xbar;a,b)=0$ for $\Leb\times\Pcal^2$-almost every $\left(\xbar,a,b\right)\in\Rd\times\Ecal^2$. Then $ D(a,b)=0$ for $\Pcal^2$-almost every $\left(a,b\right)\in\Ecal^2$, and $\Opconnf f\equiv0$ for all $f\in L^2(\Ecal)$. This contradicts $\OpNorm{\Opconnf }>0$. Therefore there exist measurable sets $\Lambda\subset \Rd$ and $F,G\subset \Ecal$ such that $\Leb\left(\Lambda\right)>0$, $\Pcal\left(F\right)>0$, $\Pcal\left(G\right)>0$, and 
\begin{equation}
    \essinf_{\xbar\in\Lambda,a\in F,b\in G}\connf(\xbar;a,b) >0.
\end{equation}
Let us denote $\Pmin := \min\left\{\Pcal(F),\Pcal(G)\right\}>0$.
Since $\Leb\left(\Lambda\right)>0$ and $d\geq 2$, there exist linearly independent $\xbar_1,\xbar_2\in \Lambda$ and $\delta>0$ small enough that 
\begin{equation}
    \varepsilon:= \min_{i=1,2}\essinf_{\substack{a\in F,b\in G,\\ \ybar\in\xbar_i + B_\delta}}\connf\left(\ybar;a,b\right) >0,
\end{equation}
where $B_\delta = \left[-\delta/2,\delta/2\right]^d$. We also assume that $\delta>0$ is sufficiently small that the sets of the form $m\xbar_1 + n\xbar_2 + B_\delta$ with $m,n\in\Z$ are disjoint.

By using Mecke's formula, given $\left(\xbar,a\right)\in B_\delta\times F$ and $i=1,2$, the number of vertices in $\left(\xbar_i+B_\delta\right)\times G$ that are adjacent to $\left(\xbar,a\right)$ is a Poisson random variable with mean:
\begin{align}
    \E_{\lambda}\left[\#\left\{y\in\eta\cap\left(\left(\xbar_i+B_\delta\right)\times G\right) \colon y \sim \left(\xbar,a\right)\right\}\right] &= \lambda \int_{G} \int_{\xbar_i+B_\delta+\xbar} \connf\left(\ybar;a,b\right)\dd \ybar \Pcal\left(\dd b\right)\nonumber\\
    &\geq \lambda \int_{G} \int_{\xbar_i+\left(B_\delta+\xbar\right)\cap B_\delta} \connf\left(\ybar;a,b\right)\dd \ybar \Pcal\left(\dd b\right)\nonumber\\
    &\geq \lambda \Pcal\left(G\right)\left(\frac{\delta}{2}\right)^d\varepsilon \geq \lambda \Pmin\left(\frac{\delta}{2}\right)^d\varepsilon.
\end{align}
Here we used that $\xbar_i+\left(B_\delta+\xbar\right)\cap B_\delta\subset \xbar_i+ B_\delta$, and that $\Leb\left(\left(B_\delta+\xbar\right)\cap B_\delta\right) \geq 2^{-d}\Leb\left(B_\delta\right)$ uniformly in $\xbar\in B_\delta$. It then follows that the probability that there exists such a vertex in $\left(\xbar_i+B_\delta\right)\times G$ has the lower bound
\begin{multline}
    \pla\left(\exists x\in\eta\cap\left(\left(\xbar_i+B_\delta\right)\times G\right) \colon x \sim \origin{a}\right) \\= 1 - \exp\left(-\E_{\lambda}\left[\#\left\{x\in\eta\cap\left(\left(\xbar_i+B_\delta\right)\times G\right) \colon x \sim \origin{a}\right\}\right]\right) \geq 1- \e^{-\lambda \Pmin\left(\frac{\delta}{2}\right)^d\varepsilon}.
\end{multline}
Since $\connf(\ybar;a,b) = \connf(\ybar;b,a)$, this calculation can be repeated for $a\in G$ to get
\begin{equation}
    \pla\left(\exists x\in\eta\cap\left(\left(\xbar_i+B_\delta\right)\times F\right) \colon x \sim \origin{a}\right) \geq 1- \e^{-\lambda \Pmin\left(\frac{\delta}{2}\right)^d\varepsilon}.
\end{equation}

Given an instance $\xi^\origin{a}$ of the marked RCM with $a\in F$, we construct a discrete tree $\mathfrak{D}(0)\subset \Z^2$ by the following algorithm.
\begin{enumerate}[label=Step \arabic*)]
\item \label{Step 1}Place $\zerobar\in \mathfrak{D}(0)$. Denote $X_{0,0}=\origin{a}\in B_\delta\times F$.

\item \label{Step 2}Consider some edge $e$ of $\Z^2$ connecting nearest neighbours $(m,n),(m',n')\in\Z^2$ such that $(m,n)\in \mathfrak{D}(0)$, $(m',n')\notin \mathfrak{D}(0)$, and $e$ has not previously been considered. If no such $e$ exists, stop.

\item \label{Step 3}If $m+n$ is even and there exists a vertex in $\left(m'\xbar_1 + n'\xbar_2 + B_\delta\right)\times G$ that is adjacent to $X_{m,n}\in \left(m\xbar_1 + n\xbar_2 + B_\delta\right)\times F$ in $\xi^\origin{a}$, add $(m',n')$ to the set $\mathfrak{D}(0)$ and denote this new vertex $X_{m',n'}$. 

If $m+n$ is odd and there exists a vertex in $\left(m'\xbar_1 + n'\xbar_2 + B_\delta\right)\times F$ that is adjacent to $X_{m,n}\in \left(m\xbar_1 + n\xbar_2 + B_\delta\right)\times G$ in $\xi^\origin{a}$, add $(m',n')$ to the set $\mathfrak{D}(0)$ and denote this new vertex $X_{m',n'}$.

\item \label{Step 4}Return to \ref{Step 2}.
\end{enumerate}

If the cluster $\C\left(\origin{a},\xi^{\origin{a}}\right)$ is finite, then so is $\mathfrak{D}(0)$. Also note that if $m+n$ is even, then any $(m',n')\in\Z^2$ that is a nearest neighbour has $m'+n'$ odd, and vice versa. Therefore for every $(m,n)\in \mathfrak{D}(0)$, $X_{m,n}\in \left(m\xbar_1 + n\xbar_2 + B_\delta\right)\times F$ if $n+m$ is even and $X_{m,n}\in \left(m\xbar_1 + n\xbar_2 + B_\delta\right)\times G$ if $m+n$ is odd. 

For each edge we consider, the probability of adding a new point of $\Z^2$ to $\mathfrak{D}(0)$ is independent of other edges and is bounded below by $1-\exp\left(-\lambda \Pmin\left(\frac{\delta}{2}\right)^d\varepsilon\right)$. The probability that $\mathfrak{D}(0)$ is finite is then less than or equal to the probability that the cluster of $\zerobar$ is finite in the Bernoulli bond model on $\Z^2$ with edge probability $1-\exp\left(-\lambda \Pmin\left(\frac{\delta}{2}\right)^d\varepsilon\right)$. It was famously proven in \cite{kesten1980critical} that the critical edge probability for that model is $p_{\rm c}=1/2$, and so for $\lambda> \tfrac{\log 2}{\Pmin\varepsilon}\left(\frac{2}{\delta}\right)^d$ the tree $\mathfrak{D}(0)$ is infinite with some positive probability. Hence $\lambda> \tfrac{\log 2}{\Pmin\varepsilon}\left(\frac{2}{\delta}\right)^d$ implies that $\theta_\lambda(a) > 0$. Since this holds for $\Pcal$-almost all $a\in F$ and $\Pcal(F)>0$, we have $\norm*{\theta_\lambda}_\infty>0$ for $\lambda> \tfrac{\log 2}{\Pmin\varepsilon}\left(\frac{2}{\delta}\right)^d$. Hence $\lambda_O \leq \lambda_c \leq \tfrac{\log 2}{\Pmin\varepsilon}\left(\frac{2}{\delta}\right)^d <\infty$.
\end{proof}

The following lemma gives sufficient conditions for the equivalence of the operator-based critical intensity and the susceptibility-based critical intensities.

\begin{lemma}
\label{lem:lambdaO_lambdaT}
The critical intensities satisfy
\begin{equation}
    \lambda^{(\infty)}_T \leq \lambda_O \leq \lambda^{(2)}_T.
\end{equation}
Therefore $\norm*{ D}_{2,\infty}<\infty$ implies that $\lambda_O = \lambda^{(p)}_T$ for all $p\in\left[2,\infty\right]$.
\end{lemma}

\begin{proof}
From \cite[Lemma~3.2]{DicHey2022triangle}, we know $\OpNorm{\Optlam }\leq \norm*{ T_\lambda}_{1,\infty}$ (this is proven by Schur's test). Therefore $\lambda^{(\infty)}_T \leq \lambda_O$ holds. Also note that $\chi_\lambda(a) = 1 + \lambda\left(\Optlam 1\right)(a)$. Therefore $\OpNorm{\Optlam }<\infty$ and the triangle inequality imply that $\norm*{\chi_\lambda}_2<\infty$, and hence we have $\lambda_O \leq \lambda^{(2)}_T$.

Lastly, Lemma~\ref{thm:equality of lambda_T} shows that $\lambda^{(\infty)}_T =\lambda^{(2)}_T$ if $\norm*{ D}_{2,\infty}<\infty$. This sandwiches all the relevant critical intensities and therefore proves the result.
\end{proof}

\subsection{Susceptibility Lower Bound}
\label{sec:susceptibilitylowerbound}

To get the desired lower bound on the susceptibility function, we first derive a similar bound on the two-point operator in very high generality. This uses a result on the derivatives of the two-point function and two-point operator, for which we draw on \cite{DicHey2022triangle} (which in turn was influenced by \cite[Lemma~2.2 \& Lemma~2.3]{HeyHofLasMat19}). Once we have our general lower bound on the two-point operator, we can use \ref{Assump:BoundExpectedDegree} to relate this to the susceptibility function and prove Theorem~\ref{thm:Susceptibility Mean-Field Bound}.

We now recall a result on the differentiability of the two-point function and two-point operator from \cite{DicHey2022triangle}.

\begin{proposition}
\label{thm:DifferentiateTlam}
Let $x,y\in\X$. Then $\tlam(x,y)$ is differentiable w.r.t.~$\lambda$ on $\left[0,\lambda_O\right)$ and 
\begin{equation}
\label{eqn:functionDerivative}
    \frac{\dd}{\dd\lambda}\tlam(x,y) = \int \pla(\conn{y}{x}{\xi^{x,y,u}}, \nconn{y}{x}{\xi^{x,y}}) \nu\left(\dd u\right).
\end{equation}
Furthermore, the operator $\Optlam $ is differentiable w.r.t.~$\lambda$ on $\left[0,\lambda_O\right)$ and 
\begin{equation}
    \OpNorm{\frac{\dd}{\dd \lambda}\Optlam} \leq \OpNorm{\Optlam}^2.
\end{equation}
\end{proposition}

\begin{proof}
This is proven in Lemma~3.7 and Corollary~3.8 of \cite{DicHey2022triangle}. The proof of \eqref{eqn:functionDerivative} uses a truncation of $\tlam(x,y)$ upon which the Margulis-Russo formula can be applied. The operator bound for the derivative then follows from applying the BK inequality to the probability on the right hand side of \eqref{eqn:functionDerivative} and then using the definition of the operator norm.
\end{proof}

\begin{theorem}
\label{thm:OperatorLowerBound}
For $\lambda\in\left(0,\lambda_O\right)$,
\begin{equation}
\label{eqn:OpNormMeanFieldBound}
    \OpNorm{\Optlam } \geq \frac{1}{\lambda_O-\lambda}.
\end{equation}
Now let $q\in\left[2,\infty\right]$ and suppose $\norm*{ D}_{q,\infty}<\infty$. Then for $p\in\left[\tfrac{q}{q-1},\infty\right]$ we have
\begin{equation}
    \norm{\chi_\lambda}_p \geq \frac{1}{\norm*{ D}_{\frac{p}{p-1},\infty}}\frac{1}{\lambda^{(p)}_T-\lambda}.
\end{equation}
\end{theorem}

\begin{proof}
We prove \eqref{eqn:OpNormMeanFieldBound} by first deducing $\OpNorm{\OptlamT }=\infty$ and then extrapolating the bound by a differential inequality.

We evaluate $\OpNorm{\OptlamT }$ by considering a sequence of truncated operators. Let us introduce $\Lambda_n=\left[-n,n\right]^d\times\Ecal\subset \X$, and define $\tau^{(n)}_\lambda\colon \X^2\to\left[0,1\right]$ as
\begin{equation}
    \tau^{(n)}_\lambda(x,y) := \pla\left(\conn{x}{y}{\xi^{x,y}_{\Lambda_n}}\right).
\end{equation}
We also define
\begin{equation}
    t^{(n)}_\lambda(a,b) := \esssup_{\xbar\in\left[-n,n\right]^d}\int \tau^{(n)}_\lambda\left(\left(\xbar,a\right),\left(\ybar,b\right)\right) \dd \ybar.
\end{equation}
Since $\tau^{(n)}_\lambda(x,y)\nearrow \tlam(x,y)$ monotonically for all $x,y\in\X$ as $n\to\infty$, we have $t^{(n)}_\lambda(a,b) \nearrow  T_\lambda(a,b)$ monotonically for all $a,b\in\Ecal$ as $n\to\infty$. If we define $\Optnlam :L^2(\Ecal)\to L^2(\Ecal)$ by
\begin{equation}
    \left(\Optnlam  f\right)\left(a\right) = \int t^{(n)}_\lambda(a,b) f(b) \Pcal\left(\dd b\right)
\end{equation}
for all $f\in L^2(\Ecal)$ and $\Pcal$-almost every $a\in\Ecal$, then $\OpNorm{\Optnlam }\nearrow \OpNorm{\Optlam }$ monotonically and $1/\OpNorm{\Optnlam } \searrow 1/ \OpNorm{\Optlam }$ monotonically as $n\to\infty$. We now want to show that the functions $\lambda\mapsto 1/ \OpNorm{\Optnlam }$ are equicontinuous, in order to show that the monotone pointwise limit of these continuous functions is also continuous.

To show equicontinuity, we aim to uniformly (in $\lambda$) bound the rate of change of $1/\OpNorm{\Optnlam }$. By the monotone increasing in $\lambda$ property of $\OpNorm{\Optnlam }$, and by triangle inequality, we have
\begin{align}
    \liminf_{\varepsilon\to0}\frac{1}{\varepsilon}\left(\OpNorm{\mathcal{T}^{(n)}_{\lambda+\varepsilon}} - \OpNorm{\Optnlam }\right)& \geq 0\\
    \limsup_{\varepsilon\to0}\frac{1}{\varepsilon}\left(\OpNorm{\mathcal{T}^{(n)}_{\lambda+\varepsilon}} - \OpNorm{\Optnlam }\right)& \leq \limsup_{\varepsilon\to0}\OpNorm{\frac{1}{\varepsilon}\left(\mathcal{T}^{(n)}_{\lambda+\varepsilon} - \Optnlam \right)}.
\end{align}
An adaptation of the proof of Proposition~\ref{thm:DifferentiateTlam} then tells us that $\Optnlam $ is differentiable with respect to $\lambda$ for $\lambda\in\left(0,\infty\right)$, with the bound $\OpNorm{\frac{\dd}{\dd \lambda}\Optnlam }\leq \OpNorm{\Optnlam }^2$. Note that while the results for $\Optlam $ in Proposition~\ref{thm:DifferentiateTlam} hold only for $\lambda\in\left(0,\lambda_O\right)$, this adapted result holds for $\lambda\in\left(0,\infty\right)$. The requirement that $\lambda<\lambda_O$ in the former arises from the need that certain integrals are finite and so certain limits can be exchanged. For our adapted result, the truncation parameter $n$ ensures that all the required integrals are finite for any $\lambda$ we wish to choose. Hence the result holds for the larger range of $\lambda$. In summary, for any $n\in\N$ and $\lambda\in\left(0,\infty\right)$,
\begin{equation}
    \limsup_{\varepsilon\to0}\frac{1}{\varepsilon}\left(\OpNorm{\mathcal{T}^{(n)}_{\lambda+\varepsilon}} - \OpNorm{\Optnlam }\right) \leq \lim_{\varepsilon\to0}\OpNorm{\frac{1}{\varepsilon}\left(\mathcal{T}^{(n)}_{\lambda+\varepsilon} - \Optnlam \right)} \leq \OpNorm{\Optnlam }^2.
\end{equation}

Now since $\OpNorm{\Optnlam }$ is non-decreasing in $\lambda$, and by a variant of the chain rule respectively, we get the following bounds for any $n\in\N$ and $\lambda\in\left(0,\infty\right)$:
\begin{align}
    \limsup_{\varepsilon\to0}\frac{1}{\varepsilon}\left(\frac{1}{\OpNorm{\mathcal{T}^{(n)}_{\lambda+\varepsilon}}} - \frac{1}{\OpNorm{\Optnlam }}\right)& \leq 0,\\
    \liminf_{\varepsilon\to0}\frac{1}{\varepsilon}\left(\frac{1}{\OpNorm{\mathcal{T}^{(n)}_{\lambda+\varepsilon}}} - \frac{1}{\OpNorm{\Optnlam }}\right)& \geq -1.
\end{align}
These uniform bounds on the rate of change of $1/\OpNorm{\Optnlam }$ then prove that the functions $\lambda\mapsto 1/ \OpNorm{\Optnlam }$ are equicontinuous. This in turn proves that the pointwise limit $\lambda\mapsto 1/ \OpNorm{\Optlam }$ is continuous on $\left(0,\infty\right)$. By monotonicity and the definition of $\lambda_O$ we have $\OpNorm{\Optlam }=\infty$ for $\lambda>\lambda_O$, and therefore we now have $\OpNorm{\OptlamT }=\infty$. 

From \cite[Corollary~3.8]{DicHey2022triangle}, we know that $\lambda\mapsto\Optlam $ is differentiable and
\begin{equation}
    \OpNorm{\frac{\dd}{\dd \lambda}\Optlam } \leq \OpNorm{\Optlam }^2
\end{equation}
for $\lambda\in\left(0,\lambda_O\right)$. As for the truncated version above, this inequality and the triangle inequality then implies that
\begin{equation}
    \limsup_{\varepsilon \to 0}\frac{1}{\varepsilon}\left(\OpNorm{\mathcal{T}_{\lambda+\varepsilon}}- \OpNorm{\Optlam }\right) \leq \OpNorm{\Optlam }^2.
\end{equation}
Therefore, by a variation of the chain rule,
\begin{equation}
    \liminf_{\varepsilon\to 0}\frac{1}{\varepsilon}\left(\frac{1}{\OpNorm{\mathcal{T}_{\lambda+\varepsilon}}}- \frac{1}{\OpNorm{\Optlam }}\right) \geq -1.
\end{equation}
Then `integrating' from $\lambda$ to $\lambda_O$ gives
\begin{equation}
    -\frac{1}{\OpNorm{\Optlam }} = \frac{1}{\OpNorm{\OptlamT }} - \frac{1}{\OpNorm{\Optlam }} \geq -\left(\lambda_O-\lambda\right).
\end{equation}
Rearranging this inequality gives \eqref{eqn:OpNormMeanFieldBound} as required.

Having $\norm*{ D}_{q,\infty}<\infty$ for some $q\in\left[2,\infty\right]$ implies $\norm*{ D}_{2,\infty}<\infty$ (by Jensen's inequality), and therefore Lemma~\ref{lem:lambdaO_lambdaT} implies $\lambda_O = \lambda^{(p)}_T$ for all $p\in\left[2,\infty\right]$. We also have $\norm*{ T_\lambda}_{1,\infty} \geq \OpNorm{\Optlam }$ from Schur's test (see \cite[Lemma~3.2]{DicHey2022triangle}), and we get
\begin{equation}
    \norm*{ T_\lambda}_{1,\infty} \geq \OpNorm{\Optlam } \geq \frac{1}{\lambda_O-\lambda} = \frac{1}{\lambda^{(\infty)}_T-\lambda},
\end{equation}
and
\begin{equation}
    \norm*{\chi_\lambda}_\infty \geq 1 + \frac{\lambda}{\lambda^{(\infty)}_T-\lambda}
\end{equation}
for all $\lambda\in\left[0,\lambda^{(\infty)}_T\right)$. We then relate $\norm*{\chi_\lambda}_\infty$ and $\norm*{\chi_\lambda}_p$ using \eqref{eqn:norm_relation} to get
\begin{equation}
    \norm*{\chi_\lambda}_p \geq \frac{1}{\norm*{ D}_{\frac{p}{p-1},\infty}}\frac{1}{\lambda^{(\infty)}_T-\lambda}
\end{equation}
for all $\lambda\in\left[0,\lambda^{(\infty)}_T\right)$. Finally, note that Lemma~\ref{thm:equality of lambda_T} and $\norm*{ D}_{q,\infty}<\infty$ imply that $\lambda^{(p)}_T = \lambda^{(\infty)}_T$ for all $p\in\left[\frac{q}{q-1},\infty\right]$, and the result is proven.
\end{proof}

\begin{proof}[Proof of Theorem~\ref{thm:Susceptibility Mean-Field Bound}]
The condition \ref{Assump:BoundExpectedDegree} means we can apply Theorem~\ref{thm:OperatorLowerBound} with $q=\infty$ and $p=1$. With Lemma~\ref{thm:equality of lambda_T}, \ref{Assump:BoundExpectedDegree} also implies $\lambda^{(p)}_T = \lambda_T$ for all $p\in\left[1,\infty\right]$, and therefore the bound on $\norm*{\chi_\lambda}_1$ from Theorem~\ref{thm:OperatorLowerBound} gives us our result for all $\norm*{\chi_\lambda}_p$ by Jensen's inequality.
\end{proof}

\subsection{Susceptibility Upper Bound}
\label{sec:susceptibilityupperbound}

To prove our upper bound on the susceptibility function, we adapt the proof from \cite{HeyHofLasMat19} for the non-marked RCM. The major innovation here is the use of Assumptions \ref{Assump:BoundExpectedDegree} and \ref{Assump:AllReachablebySome} to control the mark-dependence of the susceptibility. Central to controlling this variation, we recall the following notation:
\begin{equation}
    \Ical_{\lambda,a} := \left(\sup_{k\geq 1}\left(\frac{\lambda}{1+\lambda}\right)^k\essinf_{b\in\Ecal}D^{(k)}(a,b)\right)^{-1}.
\end{equation}
We want to make some comments about $\Ical_{\lambda,a}$. 
\begin{itemize}
    \item From the definition, $\Ical_{0,a} = \infty$ for all $a\in\Ecal$.

    \item Since $\lambda\mapsto\frac{\lambda}{1+\lambda}$ is non-decreasing, $\lambda\mapsto\Ical_{\lambda,a}$ is non-increasing for all $a\in\Ecal$. Therefore $\lambda\mapsto\esssup_{a\in\Ecal}\Ical_{\lambda,a}$ is also non-increasing.
\end{itemize}
Before we make some more comments on $\Ical_{\lambda,a}$, we need to make an observation about Assumption~\ref{Assump:AllReachablebySome}.
\begin{remark}
\label{obs:D2equivalent}
    The condition in \ref{Assump:AllReachablebySome} says that there exist $\varepsilon>0$, $k\geq 1$, and a $\Pcal$-positive set $A\subset\Ecal$ such that for all $a\in A$ and $\Pcal$-almost every $b\in\Ecal$ we have $D^{(k)}(a,b) > \varepsilon$. From the symmetry of $D^{(k)}$ in its two arguments and an application of Mecke's formula \eqref{eq:prelim:mecke_n}, we can see
    \begin{equation}
        D^{(2k)}(b,c) = \lambda\int D^{(k)}(b,a)D^{(k)}(a,c)\Pcal(\dd a) > \lambda\varepsilon^2 \Pcal\left(A\right) >0,
    \end{equation}
    for $\Pcal$-almost every $b,c\in\Ecal$. In particular this lower bound is independent of $b$ and $c$. Therefore we have $\essinf_{a\in\Ecal}\sup_{k\geq 1}\essinf_{b\in\Ecal}D^{(k)}(a,b)>0$. Since the essential infimum is bounded above by the essential supremum, we therefore have
    \begin{equation}
        \ref{Assump:AllReachablebySome} \iff \essinf_{a\in\Ecal}\sup_{k\geq 1}\essinf_{b\in\Ecal}D^{(k)}(a,b)>0.
    \end{equation}
    In an imprecise sense, we can then view \ref{Assump:AllReachablebySome} as saying that \emph{``every mark can be connected to every other mark in exactly $k$ steps for some $k$"}.
\end{remark}
We can now deduce the following comments.
\begin{itemize}
    \item As seen above, $\text{\ref{Assump:AllReachablebySome}}\implies\essinf_{a\in\Ecal}\sup_{k\geq 1}\essinf_{b\in\Ecal}D^{(k)}(a,b)>0$. If $\lambda>0$, then this implies that $\esssup_{a\in\Ecal}\Ical_{\lambda,a}<\infty$.

    \item Since $\lambda\mapsto\frac{\lambda}{1+\lambda}$ is continuous for $\lambda\geq 0$, $\lambda\mapsto\sup_{k\geq 1}\left(\frac{\lambda}{1+\lambda}\right)^k\essinf_{b\in\Ecal}D^{(k)}(a,b)$ is continuous for all $a\in\Ecal$ and $\lambda\geq 0$. Since \ref{Assump:AllReachablebySome} implies $\esssup_{a\in\Ecal}\Ical_{\lambda,a}<\infty$, we then have $\lambda\mapsto\esssup_{a\in\Ecal}\Ical_{\lambda,a}$ is also continuous for $\lambda> 0$.
\end{itemize}

We now see how we can make use of $\Ical_{\lambda,a}$.
\begin{lemma}
\label{thm:inf-pointwise bound susceptibility}
For all $\lambda>0$ and $a\in\Ecal$ such that $\int D\left(a,b\right)\Pcal\left(\dd b\right)<\infty$,
\begin{equation}
\label{eqn:inf-pointwise bound susceptibility}
    \norm*{\chi_\lambda}_\infty \leq \left(1+\lambda\InfNorm{D}\Ical_{\lambda,a}\right) \chi_\lambda(a).
\end{equation}
\end{lemma}

\begin{proof}
We begin by bounding $\norm*{\chi_\lambda}_\infty$ from above using $\norm*{\chi_\lambda}_1$. From \eqref{eqn:norm_relation} in the proof of Lemma~\ref{thm:equality of lambda_T}, we have
\begin{equation}
\label{eqn:SuspInf_vs_Susp1}
    \norm*{\chi_\lambda}_{\infty} \leq 1 + \lambda \norm*{\chi_\lambda}_{1} \norm*{D}_{\infty,\infty}.
\end{equation}

Now we aim to get a lower bound on $\chi_\lambda\left(a\right)$ using $\norm*{\chi_\lambda}_{1}$. Given a (possibly augmented) configuration $\xi$, we construct a new configuration $\left[\xi\right]^\origin{a}$. First, let $\left\{u_x\right\}_{x\in\eta}$ be a sequence of independent and identically distributed $\textrm{Uniform}\left(0,1\right)$ random variables that is also independent of $\xi$. Now let $\Ncal_\origin{a}:=\left\{y\in\eta\colon \adjconn{\origin{a}}{y}{\xi^\origin{a}}\right\}$ be the set of neighbours of $\origin{a}$ in $\xi^\origin{a}$ and observe that $\abs*{\Ncal_{\origin{a}}}$ is a Poisson random variable with mean $\lambda\int D\left(a,b\right)\Pcal\left(\dd b\right)<\infty$. If $\abs*{\Ncal_{\origin{a}}}=0$, then we set $\left[\xi\right]^\origin{a}=\xi^\origin{a}$. Otherwise the finiteness of $\E_\lambda\left[\abs*{\Ncal_{\origin{a}}}\right]$ implies that $\abs*{\Ncal_{\origin{a}}}$ is almost surely finite, and so 
\begin{equation}
    \left[x\right] := \argmax\left\{u_x\colon \adjconn{\origin{a}}{x}{\xi^\origin{a}}\right\}
\end{equation}
is almost surely defined given $\abs*{\Ncal_{\origin{a}}}\geq 1$. In particular, $\left[x\right]$ is uniformly chosen from $\Ncal_{\origin{a}}$. Then $\left[\xi\right]^\origin{a}$ is constructed from $\xi^\origin{a}$ by removing all of the edges of $\origin{a}$ except the edge between $\origin{a}$ and $\left[x\right]$.

Since every vertex and edge in $\left[\xi\right]^\origin{a}$ is also present in $\xi^\origin{a}$, we clearly have $\chi_\lambda(a) \geq \E_\lambda\left[\abs*{\C\left(\origin{a},\left[\xi\right]^\origin{a}\right)}\right]$. Since $\abs*{\Ncal_{\origin{a}}}$ is a Poisson random variable with mean $\lambda\int D\left(a,b\right)\Pcal\left(\dd b\right)$, we have
\begin{multline}
    \pla\left(\left[x\right]\text{ is defined}\right) = \pla\left(\#\left\{y\in\eta\colon \adjconn{\origin{a}}{y}{\xi^\origin{a}}\right\}\geq 1\right) \\= 1-\exp\left(-\lambda\int D\left(a,b\right)\Pcal\left(\dd b\right)\right).
\end{multline}
Then given any measurable set $B\subset \X$, we have
\begin{equation}
\label{eqn:NhbrDensity}
    \pla\left(\left[x\right]\in B\;\middle|\; \left[x\right]\text{ is defined} \right) = \frac{\int_B\connf\left(u,\origin{a}\right)\nu\left(\dd u\right)}{\int_\X\connf\left(v,\origin{a}\right)\nu\left(\dd v\right)}.
\end{equation}
For the moment we assume that $\pla\left(\abs*{\C\left(\origin{a},\left[\xi\right]^\origin{a}\right)}=\infty\;\middle|\; \left[x\right]=x, \left[x\right]\text{ is defined}\right)=0$. Then we can expand the following density:
\begin{align}
    &\E_\lambda\left[\abs*{\C\left(\origin{a},\left[\xi\right]^\origin{a}\right)}\;\middle|\; \left[x\right]=x, \left[x\right]\text{ is defined}\right]  \nonumber\\
    &\hspace{3cm}=\sum^\infty_{n=1}n\pla\left(\abs*{\C\left(\origin{a},\left[\xi\right]^\origin{a}\right)}=n\;\middle|\; \left[x\right]=x, \left[x\right]\text{ is defined}\right)\nonumber\\
    &\hspace{3cm}=\sum^\infty_{n=1}n\pla\left(\abs*{\C\left(x,\xi\right)}=n-1\;\middle|\; \left[x\right]=x, \left[x\right]\text{ is defined}\right).
\end{align}
This second equality holds because having $\left[x\right]=x$ implies that $x\in\eta$ and because every vertex connected to $\origin{a}$ in $\left[\xi\right]^\origin{a}$ must be connected to $\left[x\right]$ in $\xi$ from the construction of $\left[\xi\right]^\origin{a}$. Now let us condition on there being $k\geq 1$ neighbours of $\origin{a}$ in $\xi^\origin{a}$. These $k$ vertices are independent and identically distributed with density given by \eqref{eqn:NhbrDensity}. Since $\left[x\right]$ is chosen uniformly from these $k$ vertices, the distribution of the other $k-1$ neighbours of $\origin{a}$ in $\xi^\origin{a}$ (and the other vertices in $\xi$) are independent of the position of $\left[x\right]$. The conditioning event that $\left\{\left[x\right]=x\right\}\cap\left\{\left[x\right]\text{ is defined}\right\}$ does imply that $x\in\xi$ though, and so
\begin{equation}
    \pla\left(\abs*{\C\left(x,\xi\right)}=n-1\;\middle|\; \left[x\right]=x, \left[x\right]\text{ is defined}\right) = \pla\left(\abs*{\C\left(x,\xi^x\right)}=n-1\right).
\end{equation}
For the moment assume that for $x=\left(\xbar,b\right)\in\X$ we have $\pla\left(\abs*{\C\left(x,\xi^x\right)}=\infty\right)=0$. Then by using the fact that
\begin{multline}
    \pla\left(\abs*{\C\left(\origin{a},\left[\xi\right]^\origin{a}\right)}=\infty\;\middle|\; \left[x\right]=x, \left[x\right]\text{ is defined}\right)\\
    = \pla\left(\abs*{\C\left(x,\xi\right)}=\infty\;\middle|\; \left[x\right]=x, \left[x\right]\text{ is defined}\right) = \pla\left(\abs*{\C\left(x,\xi^x\right)}=\infty\right),
\end{multline}
we get that
\begin{multline}
    \E_\lambda\left[\abs*{\C\left(\origin{a},\left[\xi\right]^\origin{a}\right)}\;\middle|\; \left[x\right]=\left(\xbar,b\right), \left[x\right]\text{ is defined}\right] \\= \sum^\infty_{n=1}n\pla\left(\abs*{\C\left(\left(\xbar,b\right),\xi^{\left(\xbar,b\right)}\right)}=n-1\right) = 1 + \chi_\lambda\left(b\right).
\end{multline}
On the other hand, if $\pla\left(\abs*{\C\left(x,\xi^x\right)}=\infty\right)>0$, then
\begin{equation}
    \E_\lambda\left[\abs*{\C\left(\origin{a},\left[\xi\right]^\origin{a}\right)}\;\middle|\; \left[x\right]=\left(\xbar,b\right), \left[x\right]\text{ is defined}\right] = \infty = 1 + \chi_\lambda\left(b\right)
\end{equation}
as well.

In summary, we can bound
\begin{align}
    \chi_\lambda\left(a\right) &\geq \E_\lambda\left[\abs*{\C\left(\origin{a},\left[\xi\right]^\origin{a}\right)}\right]\nonumber\\
    & \geq \E_\lambda\left[\abs*{\C\left(\origin{a},\left[\xi\right]^\origin{a}\right)}\Id_{\left[x\right]\text{ is defined}}\right]\nonumber\\
    & = \int \left(1 + \chi_\lambda\left(b\right)\right)\frac{\connf\left(\left(\ubar,b\right),\origin{a}\right)}{\int\connf(v,\origin{a})\nu(\dd v)}\left(1- \e^{-\lambda\int D(a,b)\Pcal(\dd b)}\right)\dd \ubar\Pcal\left(\dd b\right)\nonumber\\
    & \geq \frac{1- \e^{-\lambda\int D(a,b)\Pcal(\dd b)}}{\int D(a,b)\Pcal(\dd b)}\int D(a,b)\chi_\lambda\left(b\right) \Pcal\left(\dd b\right) \nonumber\\
    & \geq \frac{\lambda}{1+\lambda}\int D(a,b)\chi_\lambda\left(b\right) \Pcal\left(\dd b\right).  \label{eqn:Susp_vs_Susp1}
\end{align}
For this last inequality we start from the inequality $\e^z\leq \frac{1}{1-z}$ for all $z<1$, and therefore $1-e^{-y}\geq \frac{y}{1+y}$ for all $y>-1$. Letting $y=\lambda\int D(a,b)\Pcal(\dd b)$ gives the inequality written above.

Now let $k\geq 1$ be a fixed integer. By iteratively applying \eqref{eqn:Susp_vs_Susp1} $k$ times, we arrive at
\begin{equation}
    \chi_\lambda\left(a\right) \geq \left(\frac{\lambda}{1+\lambda}\right)^k\int D^{(k)}(a,b) \chi_\lambda(b)\Pcal\left(\dd b\right).
\end{equation}
Then by taking an essential infimum bound over $b$, we arrive at
\begin{equation}
    \chi_\lambda\left(a\right) \geq \norm*{\chi_\lambda}_1\left(\frac{\lambda}{1+\lambda}\right)^k\essinf_{b\in\Ecal} D^{(k)}(a,b).\label{eqn:Susp_vs_SuspGivenk}
\end{equation}
This holds for all $k\geq 1$, and so we can take the supremum over $k$ to get
\begin{equation}
    \chi_\lambda\left(a\right) \geq \norm*{\chi_\lambda}_1\sup_{k\geq 1}\left(\frac{\lambda}{1+\lambda}\right)^k\essinf_{b\in\Ecal} D^{(k)}(a,b) = \frac{\norm*{\chi_\lambda}_1}{\Ical_{\lambda,a}}. \label{eqn:Susp_vs_Suspk}
\end{equation}

Bringing together \eqref{eqn:SuspInf_vs_Susp1} and \eqref{eqn:Susp_vs_Suspk} then gives
\begin{equation}
    \norm*{\chi_\lambda}_\infty \leq  1 + \lambda\InfNorm{D}\Ical_{\lambda,a}\chi_\lambda(a).
\end{equation}
Finally using the bound $\chi_\lambda(a) \geq \chi_0(a) = 1$ gives the result.
\end{proof}

\begin{remark}
\label{obs:tauinftyboundpointwise}
From Mecke's formula we have
\begin{equation}
    \chi_\lambda(a) = 1 + \lambda\int T_\lambda(a,b)\Pcal(\dd b),
\end{equation}
and therefore Lemma~\ref{thm:inf-pointwise bound susceptibility} implies that for all $a\in\Ecal$
\begin{equation}
\label{eqn:tauinftyboundpointwise}
    \norm*{T_\lambda}_{1,\infty} \leq \left(1+\lambda\InfNorm{D}\Ical_{\lambda,a}\right) \int T_\lambda(a,b)\Pcal(\dd b).
\end{equation}
Taking the essential supremum over $a\in\Ecal$ for the $\lambda\InfNorm{D}\Ical_{\lambda,a}$ term preserves the inequality, and then taking the essential infimum over $a$ everywhere gives
\begin{equation}
    \norm*{T_\lambda}_{1,\infty} \leq \left(1+\lambda\InfNorm{D}\esssup_{a\in\Ecal}\Ical_{\lambda,a}\right) \essinf_{a\in\Ecal}\int T_\lambda(a,b)\Pcal(\dd b).
\end{equation}
\end{remark}

The following argument completes the proof of Theorem~\ref{thm:Susceptibility Mean-Field Behaviour} by proving the complementary bound to Theorem~\ref{thm:Susceptibility Mean-Field Bound}. We adapt the corresponding proof in \cite{HeyHofLasMat19}, now including the non-uniqueness of the mark. This replicates some of the calculations from the diagrammatic bounds in the lace expansion described there. 

\begin{definition}
Here we use some notation that builds upon the thinning events described in Section~\ref{sec:Prelims}. Let $\xi_1, \xi_2$ be two independent edge-markings with locally finite vertex sets $\eta_1,\eta_2$.
\begin{itemize}
    \item Let $\left\{\sqconn{u}{x}{(\xi_1, \xi_2)}\right\}$ denote the event that $u\in\eta_1$ and $x\in\eta_2$, but that $x$ does not survive a $\C\left(u, \xi^{u}_1\right)$-thinning of $\eta_2$. 
    \item Let $m\in\N$ and $\vec x, \vec y \in \X^m$. We define $\bigcirc_m^\leftrightarrow((x_j, y_j)_{1 \leq j \leq m}; \xi)$ as the event that $\{\conn{x_j}{y_j}{\xi}\}$ occurs for every $1 \leq j \leq m$ with the additional requirement that every point in $\eta$ is the interior vertex of at most one of the $m$ paths, and none of the $m$ paths contains an interior vertex in the set $\{x_j\colon j\in[m]\} \cup \{y_j\colon j\in [m]\}$.
    \item Let $\bigcirc_m^\sqarrow( (x_j,y_j)_{1 \leq j \leq m}; (\xi_1,\xi_2))$ be the intersection of the following two events. Firstly, that $\bigcirc_{m-1}^\leftrightarrow((x_j,y_j)_{1 \leq j <m};\xi_1)$ occurs but no path uses $x_{m}$ or $y_{m}$ as an interior vertex. Secondly, that $\{\sqconn{x_{m}}{y_{m}}{(\xi_1[\eta_1\setminus  \{x_i, y_i\}_{1 \leq i <m}],\xi_2)}\}$ occurs in such a way that at least one point $z$ in $\xi_1$ that is responsible for thinning out $y_m$ is connected to $x_m$ by a path $\gamma$ so that $z$ as well as all interior vertices of $\gamma$ are not contained in any path of the $\bigcirc_{m-1}^\leftrightarrow((x_j,y_j)_{1 \leq j <m};\xi_1)$ event.
\end{itemize}
\end{definition}

Before embarking on the proof of Theorem~\ref{thm:Susceptibility Mean-Field Behaviour}, we introduce the following two lemmas that describe complicated objects in more tangible objects. The first lemma allows us to translate between these $\bigcirc_m^\sqarrow( (x_j,y_j)_{1 \leq j \leq m}; (\xi_1,\xi_2))$ events and the simpler $\bigcirc_m^\leftrightarrow((x_j, y_j)_{1 \leq j \leq m}; \xi_1)$ events.
\begin{lemma}[Relating $\bigcirc_m^\sqarrow$ and $\bigcirc_m^\leftrightarrow$] \label{lem:DB:squigarrow_tlam_equality}
Let $m \in\N$ and $\vec x, \vec y \in \X^m$. Let $\xi_1, \xi_2$ be two independent edge-markings. Then
\begin{equation}
    \pla\left( \bigcirc_m^\sqarrow\left(\left(x_j,y_j\right)_{1 \leq j \leq m}; \left(\xi_1^{\vec x_{[1,m]}, \vec y_{[1,m-1]}}, \xi_2^{y_m}\right) \right) \right) 
					= \pla \left( \bigcirc_m^\leftrightarrow \left(\left(x_j,y_j\right)_{1 \leq j \leq m}; \xi_1^{\left(\vec x, \vec y\right)_{[1,m]}} \right) \right).
\end{equation}
\end{lemma}

\begin{proof}
    This is \cite[Lemma~7.5]{HeyHofLasMat19}, with $\X$ taking the role of $\Rd$. This change makes no difference to the proof.
\end{proof}

The following lemma allows us to relate the event that a vertex is pivotal to connection and thinned connection events.
\begin{lemma}
\label{lem:pivot_to_expectation}
    Let $\lambda\geq 0$, and $u,x\in\X$ be distinct. Then
    \begin{equation}
        \pla\left(u \in \piv{\origin{a},x;\xi^{\origin{a},x}}\right) = \E_\lambda\left[\mathds 1_{\{\conn{\origin{a}}{u}{\xi^{\origin{a},u}}\}} \tlam^{\C(\origin{a},\xi^{\origin{a}})}(u,x) \right].
    \end{equation}
\end{lemma}

\begin{proof}
    This follows in the same manner as the proof of the Cutting-Point Lemma in \cite[Lemma~3.6]{HeyHofLasMat19}. Once again, replacing $\Rd$ with the space $\X$ makes no difference to the proof.
\end{proof}

\begin{proof}[Proof of Theorem~\ref{thm:Susceptibility Mean-Field Behaviour}]
Applying Proposition~\ref{thm:DifferentiateTlam} and Lemma~\ref{lem:pivot_to_expectation} gives
\begin{align}
    \frac{\dd}{\dd\lambda} \int\tlam(\origin{a},x)\nu\left(\dd x\right) &= \int\frac{\dd}{\dd\lambda} \tlam(\origin{a},x)\nu\left(\dd x\right)\nonumber\\
    &= \int \pla\left(u \in \piv{\origin{a},x;\xi^{\origin{a},u,x}}\right) \nu^{\otimes 2}\left(\dd u, \dd x\right)\nonumber\\
    &= \int \E_\lambda\left[\mathds 1_{\{\conn{\origin{a}}{u}{\xi^{\origin{a},u}}\}} \tlam^{\C(\origin{a},\xi^{\origin{a}})}(u,x) \right]\nu^{\otimes 2}\left(\dd u, \dd x\right).
\end{align}
Noting that $\tlam^A(u,x) = \tlam(u,x) - \pla\left(\xconn{u}{x}{\xi^{u,x}}{A}\right)$ then produces

\begin{multline}
    \frac{\dd}{\dd\lambda} \int\tlam(\origin{a},x)\nu\left(\dd x\right) =  \int \E_\lambda\left[ \mathds 1_{\left\{\conn{\origin{a}}{u}{\xi^{\origin{a},u}}\right\}} \tlam(x,u)\right] \nu^{\otimes 2}\left(\dd u, \dd x\right)\\
	- \int \E_\lambda\left[ \mathds 1_{\left\{\conn{\origin{a}}{u}{\xi_0^{\origin{a},u}}\right\}} \mathds 1_{\left\{\xconn{u}{x}{\xi_1^{u,x}}{\C(\origin{a},\xi_0^{\origin{a}})}\right\}}\right] \nu^{\otimes2}\left(\dd u, \dd x\right). 
					\label{eq:MT:ftlam_identity}
\end{multline}
The first integral on the right-hand side of~\eqref{eq:MT:ftlam_identity} can be easily simplified to get a $\nu$-convolution of two $\tlam$ functions, since $\tlam(x,u)$ has a deterministic value. For the second integral, we use the `$\sqarrow$' notation described above. We also use the notation $\tlamo(x,y)=\delta_{x,y} + \lambda\tlam(x,y)$ (where $\delta_{x,y}$ is the Dirac delta function) for compactness and readability. The second integrand on the r.h.s. of \eqref{eq:MT:ftlam_identity} can be bounded by
\begin{align}
    &\E_\lambda \left[ \mathds 1_{\{\conn{\origin{a}}{u}{\xi_0^{\origin{a},u}}\}} \sum_{y \in \eta_1^x} \mathds 1_{\{\sqconn{\origin{a}}{y}{(\xi_0^{\origin{a}}, \xi_1^{x})} \}}
					\mathds 1_{\{\conn{u}{y}{\xi_1^u}\}\circ \{\conn{y}{x}{\xi_1^x}\}} \right] \nonumber\\
		&\qquad= \E_\lambda \left[ \mathds 1_{\{\conn{\origin{a}}{u}{\xi_0^{\origin{a},u}}\}} \mathds 1_{\{\sqconn{\origin{a}}{x}{(\xi_0^{\origin{a}}, \xi_1^{x})}\}} \right] \tlam(x,u) \nonumber\\
		& \hspace{2cm} + \lambda \int \E_\lambda \left[ \mathds 1_{\{\conn{\origin{a}}{u}{\xi_0^{\origin{a},u}}\}} \mathds 1_{\{\sqconn{\origin{a}}{y}{(\xi_0^\origin{a}, \xi_1^{x,y})}\}} \right]\nonumber\\
		&\hspace{6cm}\times\pla \left( \{ \conn{u}{y}{\xi^{u,y}} \} \circ \{ \conn{y}{x}{\xi^{y,x}}\} \right) \nu\left(\dd y\right) \nonumber\\
		&\qquad \leq \int \E_\lambda \left[ \mathds 1_{\{\conn{\origin{a}}{u}{\xi_0^{\origin{a},u}}\}} \mathds 1_{\{\sqconn{\origin{a}}{y}{(\xi_0^\origin{a}, \xi_1^{y})}\}} \right]
						 \tlam(y,u) \tlamo(x,y) \nu\left(\dd y\right)\label{eq:MT:through_conn_bound}
\end{align}
where we split the sum according to whether the vertex we are summing over, namely $y$, is equal to $x$ or not. This then gets recombined using the $\tlamo$ notation. Note that
\begin{equation}
    \mathds 1_{\{\conn{\origin{a}}{u}{\xi_0^{\origin{a},u}}\}} \mathds 1_{\{\sqconn{\origin{a}}{y}{(\xi_0^\origin{a}, \xi_1^{y})}\}} \leq \sum_{v \in \eta_0^\origin{a}} 
					\mathds 1_{\bigcirc^\sqarrow_3 \left((\origin{a}, v), (v,u), (v,y); (\xi_0^{\origin{a},u}, \xi_1^y)\right)}. \label{eq:MT:sconn_eta1_treegraph_bound}
\end{equation}
We now plug~\eqref{eq:MT:through_conn_bound} back into~\eqref{eq:MT:ftlam_identity} and apply~\eqref{eq:MT:sconn_eta1_treegraph_bound}, with the intent to use Lemma~\ref{lem:DB:squigarrow_tlam_equality}. The second integral on the r.h.s.~of~\eqref{eq:MT:ftlam_identity} is hence bounded by
\begin{align}
    &\int \left( \delta_{v,\origin{a}} \pla \left( \bigcirc^\sqarrow_2 ((\origin{a},u), (\origin{a},y); (\xi_0^{\origin{a},u}, \xi_1^y))\right)\right. \nonumber\\
    &\hspace{1cm}\left.+ \lambda \pla \left( \bigcirc^\sqarrow_3 ((\origin{a}, v), (v,u), (v,y); (\xi_0^{\origin{a},v,u}, \xi_1^y)) \right) \right) \tlam(y,u) \tlamo(x,y) \nu^{\otimes 4}\left(\dd v, \dd u, \dd x, \dd y\right) \nonumber\\
	&\qquad\leq \int \tlamo(v,\origin{a}) \tlam(u,v) \tlam(y,v) \tlam(y,u) \tlamo(x,y) \nu^{\otimes 4}\left(\dd v, \dd u, \dd x, \dd y\right) \nonumber\\
	& \qquad\leq\lambda^{-2} \trilam \chi_\lambda(a)\norm*{\chi_\lambda}_\infty.
\end{align}
The above estimate is achieved by first applying Lemma~\ref{lem:DB:squigarrow_tlam_equality} and then the BK inequality. In summary, we have
\begin{multline}
    \frac{\dd}{\dd\lambda} \int\tlam(\origin{a},x)\nu\left(\dd x\right) \geq \int \tlam(\origin{a},u)\tlam(u,x)\nu^{\otimes 2}\left(\dd u, \dd x\right) - \lambda^{-2} \trilam \chi_\lambda(a)\norm*{\chi_\lambda}_\infty \\
    \geq \left(\int T_\lambda(a,b)\Pcal\left(\dd b\right)\right)\left(\essinf_{b\in\Ecal}\int T_\lambda(b,c)\Pcal\left(\dd c\right)\right) - \lambda^{-2} \trilam \chi_\lambda(a)\norm*{\chi_\lambda}_\infty.
\end{multline}
We now want to pass this derivative out through a suprema. Given $\lambda,h>0$, let $c_{\lambda,h}\in\Ecal$ be such that
\begin{equation}
    \int T_\lambda(c_{\lambda,h},b)\Pcal\left(\dd b\right) \geq \esssup_{a\in\Ecal}\int T_\lambda(a,b)\Pcal\left(\dd b\right) - h^2 = \OneNorm{ T_\lambda} - h^2.
\end{equation}
Then
\begin{multline}
    \OneNorm{ T_{\lambda+h}} - \OneNorm{ T_\lambda} \geq \int\left(T_{\lambda+h}(c_{\lambda,h},b) - T_\lambda(c_{\lambda,h},b)\right)\Pcal\left(\dd b\right) - h^2\\ \geq \essinf_{a\in\Ecal}\int\left(\tau_{\lambda+h}(\xbar;a,b) - \tlam(\xbar;a,b)\right)\dd \xbar\Pcal\left(\dd b\right) - h^2,
\end{multline}
and dividing by $h$ and taking $h\to0$ gives
\begin{equation}
    \liminf_{h\downarrow 0}\frac{1}{h}\left(\norm*{ T_{\lambda+h}}_{1,\infty} - \norm*{ T_\lambda}_{1,\infty}\right) \geq \essinf_{a\in\Ecal} \frac{\dd}{\dd\lambda} \int\tlam(\origin{a},x)\nu\left(\dd x\right).
\end{equation}
Therefore we get
\begin{equation}
    \liminf_{h\downarrow 0}\frac{1}{h}\left(\norm*{ T_{\lambda+h}}_{1,\infty} - \norm*{ T_\lambda}_{1,\infty}\right) \geq \left(\essinf_{a\in\Ecal}\int T_\lambda(a,b)\Pcal(\dd b)\right)^2 - \lambda^{-2} \trilam \norm*{\chi_\lambda}_\infty^2,
\end{equation}
for $\lambda<\lambda_T$. Rearranging then gives
\begin{align}
    &\limsup_{h\downarrow 0}\frac{1}{h}\left(\frac{1}{\norm*{ T_{\lambda+h}}_{1,\infty}} - \frac{1}{\norm*{ T_\lambda}_{1,\infty}}\right) \nonumber\\
    &\hspace{5cm}\leq -\left(\frac{\essinf_{a\in\Ecal}\int T_\lambda(a,b)\Pcal(\dd b)}{\norm*{ T_\lambda}_{1,\infty}}\right)^2 + \trilam \frac{\norm*{\chi_\lambda}_\infty^2}{\lambda^2\norm*{ T_\lambda}_{1,\infty}^2} \nonumber\\
    &\hspace{5cm}\leq -\left(\frac{1}{1+\lambda\InfNorm{D}\esssup_{a\in\Ecal}\Ical_{\lambda,a}}\right)^2 +  \trilam \frac{\norm*{\chi_\lambda}_\infty^2}{\lambda^2\norm*{ T_\lambda}_{1,\infty}^2}.
\end{align}
Note that $\chi_\lambda(a) = 1 + \lambda\int T_\lambda(a,b)\Pcal(\dd b)$ by Mecke's formula, and therefore $\norm*{\chi_\lambda}_\infty = 1 + \lambda\OneNorm{ T_\lambda}$. Furthermore, Theorem~\ref{thm:Susceptibility Mean-Field Bound} proves that $\norm*{\chi_\lambda}_\infty\to\infty$ as $\lambda\uparrow \lambda_T$. Therefore for all $\varepsilon_1>0$ there exists $\varepsilon_2>0$ such that $\frac{\norm*{\chi_\lambda}_\infty^2}{\lambda^2\norm*{ T_\lambda}_{1,\infty}^2} \leq  1+\varepsilon_1$ for $\lambda\in\left[\lambda_T-\varepsilon_2,\lambda_T\right]$. In particular for $\lambda\in\left[\lambda_T-\varepsilon_2,\lambda_T\right]$ we have
\begin{align}
    &\limsup_{h\downarrow 0}\frac{1}{h}\left(\frac{1}{\norm*{ T_{\lambda+h}}_{1,\infty}} - \frac{1}{\norm*{ T_\lambda}_{1,\infty}}\right)
    \leq -\left(\frac{1}{1+\lambda\InfNorm{D}\esssup_{a\in\Ecal}\Ical_{\lambda,a}}\right)^2 \nonumber\\
    & \hspace{4cm}\times\left(1 - \left(1+\lambda\InfNorm{D}\esssup_{a\in\Ecal}\Ical_{\lambda,a}\right)^2\left(1+\varepsilon_1\right)\trilam\right).
\end{align}
Recall that \ref{Assump:BoundExpectedDegree} and \ref{Assump:AllReachablebySome} imply that $\InfNorm{D}\esssup_{a\in\Ecal}\Ical_{\lambda,a} <\infty$. Then since $\trilam$ is non-decreasing and $\lambda\esssup_{a\in\Ecal}\Ical_{\lambda,a}$ is continuous in $\lambda$, Assumption \ref{TriangleCondition_Assumption} implies that there exist $\varepsilon,\varepsilon'>0$ such that 
\begin{equation}
    \limsup_{h\downarrow 0}\frac{1}{h}\left(\frac{1}{\norm*{ T_{\lambda+h}}_{1,\infty}} - \frac{1}{\norm*{ T_\lambda}_{1,\infty}}\right) \leq -\varepsilon
\end{equation}
uniformly for $\lambda\in\left[\lambda_T-\varepsilon',\lambda_T\right)$. Integrating between $\lambda=\lambda_T$ (where $1/\norm*{T_{\lambda_T}}_{1,\infty}=0$) and $\lambda\in\left[\lambda_T-\varepsilon',\lambda_T\right)$, and rearranging then gives 
\begin{equation}
    \norm*{T_{\lambda}}_{1,\infty} \leq \frac{1}{\varepsilon}\frac{1}{\lambda_T-\lambda}
\end{equation}
for all $\lambda\in\left[\lambda_T-\varepsilon',\lambda_T\right]$. Relating this back to $\chi_\lambda$ via Mecke's formula then gives
\begin{equation}
    \norm*{\chi_\lambda}_\infty \leq 1 + \frac{1}{\varepsilon}\frac{1}{\lambda_T-\lambda}
\end{equation}
for all $\lambda\in\left[\lambda_T-\varepsilon',\lambda_T\right]$. Since $\norm*{\chi_\lambda}_\infty \geq \norm*{\chi_\lambda}_p$ for all $p\in\left[1,\infty\right]$, this completes the proof for $\lambda\in\left[\lambda_T-\varepsilon',\lambda_T\right]$. The result can be extended to all $\lambda\in\left[0,\lambda_T\right)$ by increasing the constant out the front because $\norm*{\chi_\lambda}_\infty$ is bounded on $\lambda\in\left[0,\lambda_T-\varepsilon'\right]$.
\end{proof}

\section{Magnetization}
\label{sec:MagnetizationBounds}

We introduce in this section a continuous and mark-dependent analogy of the magnetization first introduced for percolation models by Aizenman and Barsky \cite{AizBar87}. The name has its origins in ferromagnetic Ising models, and while the physical interpretation is now much more removed, it obeys the same differential inequalities. This has utility in describing both the susceptibility and percolation functions. 

For this magnetization function, we will introduce a parameter $\gamma\in\left(0,1\right)$. To be clear, this is not the susceptibility critical exponent from Definition~\ref{def_crit_exp}, but a parameter. We risk this potential confusion to maintain consistency with literature. The susceptibility critical exponent will not appear again in this paper, and hereafter every occurrence of $\gamma$ is this magnetization parameter.

\subsection{Magnetization preliminaries}

To arrive at our results for the percolation functions and the cluster trail behaviour, we make use of the magnetization function. For $\gamma\in (0,1)$ we enrich the MRCM by adding to each vertex a (Lebesgue) uniform $\left(0,1\right)$ label (independent of everything else), and let $\plg$ denote the resulting probability measure. A vertex $x \in\eta$ is called a \emph{ghost} vertex if its label is at most $\gamma$, and we write $x \in \G$. Similarly, we write $x \longleftrightarrow \G$ if $x$ is connected to a ghost vertex. We define then the magnetization as follows
\begin{align}
\label{eqn:define_magnetization}
    M(\lambda,\gamma,a) := \plg(\conn {\origin{a}}  {\mathcal{G}} {\xi^{\origin{a}}}),
\end{align}
and 
\begin{equation}
    M_{\sup}\left(\lambda,\gamma\right) := \esssup_{a\in\Ecal}M(\lambda,\gamma,a).
\end{equation}
We are going to also use the following susceptibility-type functions. For $\lambda\geq 0$ we define the finite susceptibility function $\chi^{\mathrm{f}}_\lambda\colon\Ecal\to \left[0,\infty\right]$ as
\begin{align}
    \chi^{\mathrm{f}}_\lambda(a) &:= \E_\lambda\left[\abs*{\C(\origin{a})} \mathds{1}_{ \abs*{\C(\origin{a})} < \infty}\right] = \sum_{k\in\N}k\pl(|\C(\origin{a})| = k).
\end{align}
For $\lambda\geq 0$ and $\gamma\in\left(0,1\right)$ we also define the ``ghost-free'' susceptibility function as
\begin{equation}
    \chi(\lambda,\gamma,a) = \elg\left[|\C(\origin{a})|\mathds{1}_{\C(\origin{a}) \cap \mathcal{G} = \emptyset}\right].
\end{equation}

Before we study the magnetization in more in detail, we relate it to the susceptibility and percolation functions.

\begin{lemma}
\label{gamma_to_0}
    For all intensities $\lambda \geq 0$ and all marks $a\in\Ecal$ we have
    \begin{align}
        &\lim_{\gamma \to 0}M(\lambda,\gamma,a) = \theta_\lambda(a)\\
        &\lim_{\gamma \to 0}\chi(\lambda,\gamma,a) = \chi^{\mathrm{f}}_{\lambda}(a).
    \end{align}
\end{lemma}

\begin{proof}
    The proof is quite direct and follows the same lines as in \cite{Gri99}. For the first one we write
    \begin{align}
         M(\lambda,\gamma,a) &= 1 - \plg(\nconn {\origin{a}}  {\mathcal{G}} {\xi^{\origin{a}}}) \nonumber \\
         &= 1 - \sum_{k \in \N}\plg(\C(\origin{a}) \cap \mathcal{G} = \emptyset \mid\abs*{\C\left(\origin{a}\right)} = k)\plg(\abs*{\C\left(\origin{a}\right)} = k) \nonumber \\
         &= 1 - \sum_{k \in \N}(1 - \gamma)^k\pl(\abs*{\C\left(\origin{a}\right)} = k). \label{eqn:magnetisation_expression}
    \end{align}
This sums contains non negative elements converges for every $\gamma > 0$ so we can take the limit $\gamma \to 0$ which gives us the first result.

For the second result we condition on $|\C(\origin{a})|$, which we know is almost-surely finite under the condition $\C(\origin{a}) \cap \mathcal{G} = \emptyset$ for all $\gamma > 0$ since otherwise that cluster would contain a ghost vertex almost-surely. We have then
\begin{multline}
    \chi(\lambda,\gamma,a) = \sum_{k\in\N} k \plg( \C(\origin{a}) \cap \mathcal{G} = \emptyset \mid |\C(\origin{a})| = k) \pl(|\C(\origin{a})| = k) \\
    =  \sum_{k \in\N} k (1 - \gamma)^k \pl(|\C(\origin{a})| = k).
\end{multline}
We take then the limit $\gamma \to 0$ with monotone convergence to conclude.
\end{proof}

Lemma~\ref{lem:analyticity} below uses the Weierstrass M-test to show analyticity of the magnetization in $\lambda$ and $\gamma$ for certain domains. In particular this proves that the $\lambda$ and $\gamma$ partial derivatives exist in these domains. The Weierstrass M-test shows analyticity through the following two theorems, which can be found in standard textbooks like \cite{ahlfors1953complex}.
\begin{theorem}[Weierstrass Theorem]
    Let $f_n$ be a sequence of analytic functions defined on an open subset $\Omega$ of the complex plane, which converges uniformly on the compact subsets of $\Omega$ to a function $f$. Then $f$ is analytic on $\Omega$.
\end{theorem}
\begin{theorem}[Weierstrass M-Test]
    Let $f_n$ be a sequence of complex-valued functions defined on a subset $\Omega$ of the plane and assume that there exist positive $M_n$ with $\abs*{f_n(z)}\leq M_n$ for all $z\in\Omega$, and $\sum_n M_n<\infty$. Then $\sum_n f_n$ converges uniformly on $\Omega$.
\end{theorem}

\begin{lemma}
\label{lem:analyticity}
    For all $\lambda\geq 0$ and $a\in\Ecal$, $\gamma\mapsto M(\lambda,\gamma,a)$ is analytic on $\left(0,1\right)$.
    For all $\gamma\in\left(0,1\right)$ and $a\in\Ecal$, $\lambda\mapsto M(\lambda,\gamma,a)$ is analytic on $\left(0,\infty\right)$.
\end{lemma}

\begin{proof}
    We start from the expression \eqref{eqn:magnetisation_expression} for the magnetization, and aim to use the Weierstrass M-test to show analyticity. 

    Showing the analyticity of $M(\lambda,\gamma,a)$ with respect to $\gamma$ is the easier argument and demonstrates the approach more cleanly. Our aim is to show that the complex function \newline $z\mapsto 1-\sum_{k\geq 0}\left(1-z\right)^k \pl\left(|\C(\origin{a})| = k\right)$ is complex analytic on some open subset of the complex plane that contains the open real interval $\left(0,1\right)$. Therefore the restriction to this interval (that is, our magnetization function) is real-analytic. 
    
    Since $\pl(|\C(\origin{a})| = k)$ are all $z$-independent, $z\mapsto \left(1-z\right)^k \pl(|\C(\origin{a})| = k)$ are clearly analytic on $\Complex$ for all $k\geq 0$. Furthermore, clearly $\abs*{\pl(|\C(\origin{a})| = k)}\leq 1$. Therefore for all $\varepsilon>0$ and $z\in\left\{\zeta\in\Complex\colon \abs*{1-\zeta}<1-\varepsilon\right\}$ we have
    \begin{equation}
        \abs*{\left(1-z\right)^k \pl(|\C(\origin{a})| = k)}\leq \left(1-\varepsilon\right)^k.
    \end{equation}
    The Weierstrass M-test then shows that $z\mapsto \sum_{k\geq 0}\left(1-z\right)^k \pl\left(|\C(\origin{a})| = k\right)$ is complex analytic on $\left\{\zeta\in\Complex\colon \abs*{1-\zeta}<1-\varepsilon\right\}$ for all $\varepsilon>0$, and is therefore analytic on $\left\{\zeta\in\Complex\colon \abs*{1-\zeta}<1\right\}$. Restricting to $\left(0,1\right)$ shows that $\gamma\mapsto 1-M(\lambda,\gamma,a)$ is real analytic on $\left(0,1\right)$, and therefore so is $\gamma\mapsto M(\lambda,\gamma,a)$.

    Showing the analyticity of $\lambda\mapsto M(\lambda,\gamma,a)$ is more complicated because of the $\lambda$-dependence of $\pl\left(|\C(\origin{a})| = k\right)$. Fortunately we can use Mecke's formula to get an explicit expression for this probability. Given $k\geq 0$, the cluster has size $k+1$ if there are $k$ vertices in $\eta$ (plus $\origin{a}$) that are all connected, and no other vertices in $\eta$ are connected to these. Furthermore, if we find such a candidate set of $k$ vertices, that set is unique. Then applying Mecke's formula gives
    \begin{align}
        &\pl(|\C(\origin{a})| = k+1) \nonumber\\
        &\hspace{1cm}= \frac{1}{k!}\el\left[\sum_{\Vec{x}\in\eta^{(k)}}\Id_{\left\{\forall y\in \left\{x_1,\ldots,x_k\right\}, \conn{y}{\origin{a}}{\xi^{\origin{a}}}\right\}\cap\left\{\forall y\in\eta\setminus\left\{x_1,\ldots,x_k\right\},\nconn{y}{\origin{a}}{\xi^{\origin{a}}}\right\}}\right]\nonumber\\
        &\hspace{1cm}= \frac{\lambda^k}{k!}\int_{\X^k} \left(\sum_{G\in \mathrm{Gr}_{k+1}}\left(\prod_{\left\{i,j\right\}\in E\left(G\right)}\connf\left(x_i,x_j\right)\right)\left(\prod_{\left\{i,j\right\}\not\in E\left(G\right)}\left(1-\connf\left(x_i,x_j\right)\right)\right)\right)\nonumber\\
        &\hspace{2.5cm}\times \exp\left(-\lambda \int_{\X} \left(1-\prod_{m\in\left\{0,1,\ldots,k\right\}}\left(1-\connf\left(y,x_m\right)\right)\right)\nu\left(\dd y\right)\right)\nu^{\otimes k}\left(\dd \Vec{x}_{\left[1,\ldots,k\right]}\right),\label{eqn:clustersize_individual}
    \end{align}
    where $\mathrm{Gr}_{k+1}$ denotes the set of simple connected graphs on $\left\{0,1,\ldots,k\right\}$, $E\left(G\right)$ denotes the set of edges of the graph $G$, and $x_0=\origin{a}$. The sum over $\mathrm{Gr}_{k+1}$ accounts for the various ways that that $\left\{\origin{a},x_1,\ldots,x_k\right\}$ can be connected, and the associated factors of $\connf$ and $1-\connf$ give the probability that they are indeed connected in that way. The elements of $\eta$ that are not $x_1,\ldots,x_k$ are distributed as a Poisson point process with intensity $\lambda\nu$, and so by Mecke's formula the expected number of these extra vertices that are adjacent to some vertex in our cluster is given by $\lambda \int_{\X} \left(1-\prod_{m\in\left\{0,1,\ldots,k\right\}}\left(1-\connf\left(y,x_m\right)\right)\right)\nu\left(\dd y\right)$. This number of extra vertices is also Poisson distributed, and so the probability that there are no such vertices is given by the exponential factor in \eqref{eqn:clustersize_individual}.

    For real integer $k\geq 0$ and $z\in\Complex$, we define
    \begin{multline}
        P_{k+1}\left(z\right) := \frac{z^k}{k!}\int_{\X^k} \left(\sum_{G\in\mathrm{Gr}_{k+1}}\left(\prod_{\left\{i,j\right\}\in E\left(G\right)}\connf\left(x_i,x_j\right)\right)\left(\prod_{\left\{i,j\right\}\not\in E\left(G\right)}\left(1-\connf\left(x_i,x_j\right)\right)\right)\right)\\
        \times \exp\left(-z \int_{\X} \left(1-\prod_{m\in\left\{0,1,\ldots,k\right\}}\left(1-\connf\left(y,x_m\right)\right)\right)\nu\left(\dd y\right)\right)\nu^{\otimes k}\left(\dd \Vec{x}_{\left[1,\ldots,k\right]}\right).\label{eqn:complex_extension}
    \end{multline}
    Clearly for $z\in\R_+$ this expression coincides with $\mathbb{P}_z(|\C(\origin{a})| = k+1)$. Now for real $s,\varepsilon>0$ we define
    \begin{equation}
        \Omega_{s,\varepsilon} := \left\{\zeta\in\Complex \colon \Re (\zeta) > s, \abs*{\zeta}<s+\varepsilon\right\}.
    \end{equation}
    For $z\in\Omega_{s,\varepsilon}$ we can bound $\abs*{z^k}\leq \left(s+\varepsilon\right)^k$ and $\abs*{\e^{-c z}}\leq \e^{-cs}$ (for real $c\geq 0$) to get
    \begin{equation}
        \abs*{P_{k+1}(z)} \leq \frac{\left(s+\varepsilon\right)^k}{s^k}P_{k+1}(s) \leq \left(1+\frac{\varepsilon}{s}\right)^k. \label{eqn:complex_bound}
    \end{equation}
    Here we have used $P_{k+1}(s)=\mathbb{P}_s(|\C(\origin{a})| = k+1)\leq 1$. Given $x_0,\ldots,x_k$, the integrand of \eqref{eqn:complex_extension} is a product of $z^k$ and an exponential of $z$, and is therefore analytic in $z$. We can then view $P_{k+1}(z)$ for $z\in\Omega_{s,\varepsilon}$ as a countably infinite sum of complex analytic functions. We have shown in \eqref{eqn:complex_bound} that this sum converges absolutely, and therefore we can use the Weierstrass M-test to show that each $z\mapsto P_{k+1}(z)$ is a complex analytic function on $\Omega_{s,\varepsilon}$.

    Furthermore, for $z\in\Omega_{s,\varepsilon}$ we have
    \begin{equation}
        \abs*{\left(1-\gamma\right)^{k+1}P_{k+1}(z)} \leq \left(1-\gamma\right)^{k+1}\left(1+\frac{\varepsilon}{s}\right)^k.
    \end{equation}
    Therefore for $\varepsilon < \frac{s\gamma}{1-\gamma}$ we can use the Weierstrass M-test to show that the infinite sum $\sum_{k\geq 1}\left(1-\gamma\right)^{k} P_k(z)$ is a complex analytic function on $\Omega_{s,\varepsilon}$. Stitching together these for suitable $s$ and $\varepsilon$ shows that $\sum_{k\geq 1}\left(1-\gamma\right)^{k} P_k(z)$ is complex analytic on
    \begin{equation}
        \Omega^{(\gamma)} := \bigcup_{\substack{s>0,\\ 0<\varepsilon < \frac{s\gamma}{1-\gamma}}}\Omega_{s,\varepsilon} = \left\{\zeta\in\Complex\colon \zeta\ne 0, \abs*{\mathrm{arg}(\zeta)}<\arccos\left(1-\gamma\right)\right\}.
    \end{equation}

    \begin{figure}
        \centering
        \begin{subfigure}[b]{0.49\linewidth}
        \centering
        \begin{tikzpicture}[scale=1.14]
            \draw[thick, fill=black!10, domain=-36.9:36.9] plot ({3.33*cos(\x)}, {3.33*sin(\x)});
            \draw[thick] (2.67,2) -- (2.67,-2);
            \draw[->] (-1,0) -- (4.2,0)node[right]{$\Re(\zeta)$};
            \draw[->] (0,-2) -- (0,2)node[above]{$\Im(\zeta)$};
            \draw[<->,dashed] (0,0) -- (2.67,2);
            \path (1.33,1) node[above, rotate=36.9]{$s+\varepsilon$};
            \draw[<->,dashed] (0,-1) -- (2.67,-1);
            \path (1.33,-1) node[above]{$s$};
        \end{tikzpicture}
        \caption{$\Omega_{s,\varepsilon}\subset \Complex$}
        \label{fig:small_domain}
        \end{subfigure}
        \hfill
        \begin{subfigure}[b]{0.49\linewidth}
        \centering
        \begin{tikzpicture}[scale=1.14]
            \draw[thick, fill=black!10] (4,-2) -- (0,0) -- (4,2);
            \draw[->] (-1,0) -- (4.2,0)node[right]{$\Re(\zeta)$};
            \draw[->] (0,-2) -- (0,2)node[above]{$\Im(\zeta)$};
            \draw [thick,domain=-26.6:0] plot ({1.5*cos(\x)}, {1.5*sin(\x)});
            \path (2.5,-0.5) node{\scriptsize $\arccos\left(1-\gamma\right)$};
        \end{tikzpicture}
        \caption{$\Omega^{(\gamma)}\subset\Complex$}
        \label{fig:big_domain}
        \end{subfigure}
        \caption{Diagrams of the open complex sets $\Omega_{s,\varepsilon}$ and $\Omega^{(\gamma)}$.}
    \end{figure}
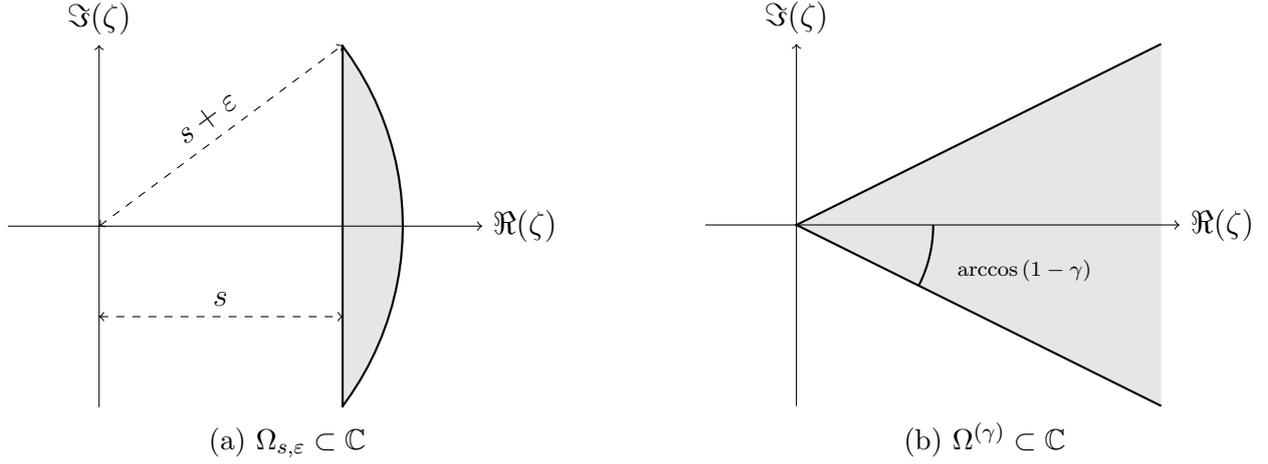
    Clearly $1-\sum_{k\geq 1}\left(1-\gamma\right)^{k} P_k(z)$ is then also complex analytic on $\Omega^{(\gamma)}$. Since $\left(0,\infty\right)\subset \Omega^{(\gamma)}\cap\R$, the restriction of this function to $\left(0,\infty\right)$ is real analytic. The restriction is exactly the magnetization $M(\lambda,\gamma,a)$.
\end{proof}

In order to get bounds for the magnetization we derive partial differential equation for it. However one crucial ingredient for doing so is the Margulis-Russo formula which only holds a priori in finite volume. Therefore we derive the differential inequalities first in the finite volume and then take the infinite volume limit, this is done in the next subsection. Before doing so we ensure that finite volume quantities converge to the appropriate infinite volume ones. 

Let $n \in \N$ and define $\Lambda_n := [-n,n]^d\times \Ecal$, then for $\lambda \geq 0$, $a\in\Ecal$, $\gamma \in (0,1)$ define the restricted magnetization as follows
\begin{align}
    M^{(n)}(\lambda,\gamma,a) :=  \plg(\conn {\origin{a}}  {\mathcal{G}} {\xi_{\Lambda_n}^{\origin{a}}}) = 1 - \sum_{k \in \N}(1 - \gamma)^k\pl\left(\abs*{\C^{(n)}\left(\origin{a}\right)} = k \right),
\end{align}
where $\C^{(n)}\left(\origin{a}\right)$ is the cluster of $\origin{a}$ in $\xi_{\Lambda_n}^{\origin{a}}$. By repeating the arguments of Lemma \ref{lem:analyticity} for $\Lambda_n$ instead of $\X$ we get the same result but on finite volume. For all $\lambda\geq 0$ and $a\in\Ecal$, $\gamma\mapsto M^{(n)}(\lambda,\gamma,a)$ is analytic on $\left(0,1\right)$. For all $\gamma\in\left(0,1\right)$ and $a\in\Ecal$, $\lambda\mapsto M^{(n)}(\lambda,\gamma,a)$ is analytic on $\left(0,\infty\right)$, which ensures the existence of the partial derivatives of the restricted magnetization. Furthermore for the convergence of the partial derivative with respect to $\lambda$ we are going to need the following assumption on our model:
\begin{enumerate}[label=\textbf{(D.\arabic* -)}]
\item \label{Assump:BoundExpectedDegreeWithAverage} \emph{``Every mark has bounded expected degree"}
    \begin{equation}
         \norm*{D}_{1,\infty} = \esssup_{a\in\Ecal}\int_\Ecal D(a,b)\Pcal(\dd b) < \infty. \label{eqn:supexpbound}
    \end{equation}
\end{enumerate}
Note that by Mecke's formula, $\esssup_{a\in\Ecal} \E_\lambda\left[\deg{\origin{a}}\right] = \lambda\esssup_{a\in\Ecal}\int_\Ecal D(a,b)\Pcal(\dd b)$ - hence the bounded expected degree description. The label \ref{Assump:BoundExpectedDegreeWithAverage} comes from the fact that it is implied by assumption \ref{Assump:BoundExpectedDegree}. We have then the following lemma.

\begin{lemma}
\label{fin_vol_lim}
Let $\gamma \in (0,1)$, $\lambda > 0$, and $a \in \Ecal$ such that $\int D(a,b)\Pcal(\dd b) < \infty$ then under the assumption \ref{Assump:BoundExpectedDegreeWithAverage} we have the following convergences as $n\to\infty$
\begin{align}
    M^{(n)}(\lambda,\gamma,a) &\to M(\lambda,\gamma,a),\label{convergence_mag}\\
    \frac{\partial M^{(n)}(\lambda,\gamma,a)}{\partial \gamma} &\to  \frac{\partial M(\lambda,\gamma,a)}{\partial \gamma}, \label{convergence_mag_partial_gamma}\\
    \frac{\partial M^{(n)}(\lambda,\gamma,a)}{\partial \lambda} &\to  \frac{\partial M(\lambda,\gamma,a)}{\partial \lambda}.\label{convergence_mag_partial_lambda}
\end{align}
\end{lemma}

\begin{remark}
    Whereas the assumption \ref{Assump:BoundExpectedDegreeWithAverage} appears somewhat naturally in the proof of \eqref{convergence_mag_partial_lambda}, it is actually a bit strong for the other two points. Indeed the reader might notice that when we use this assumption in the proof of \eqref{convergence_mag} and \eqref{convergence_mag_partial_gamma}, something weaker like having
    \begin{align}
    \label{Assumption:BoundedExpectedDegreeMinimal}
        \norm*{D}_{1,1}= \int_{\Ecal} \int_{\Ecal} D(a,b)\Pcal(\dd a)\Pcal(\dd b) < \infty,
    \end{align}
    would have been sufficient.
\end{remark}

\begin{proof}
    Let $\gamma \in (0,1)$, $\lambda > 0$, $a \in \Ecal$ such that $\int_\Ecal D(a,b)\Pcal(\dd b) < \infty$ and $n \in \N$. We assume without loss of generality that $\int D(a,b)\Pcal(\dd b) \leq \norm*{D}_{1,\infty}$. In the case that its not true, the proof still holds by replacing $\norm*{D}_{1,\infty}$ with $\max\left(\norm*{D}_{1,\infty},\int D(a,b)\Pcal(\dd b)\right)$.
    To show \eqref{convergence_mag} and \eqref{convergence_mag_partial_gamma} we follow the ideas from Lemma 4.3 in \cite{Mee95}. Notice that
    \begin{align}
        M(\lambda,\gamma,a) - M^{(n)}(\lambda,\gamma,a) = \sum_{k \in \N}(1 - \gamma)^k\left(\pl(\abs*{\C\left(\origin{a}\right)} = k) -\pl\left(\abs*{\C^{(n)}\left(\origin{a}\right)} = k \right) \right),
    \end{align}
    and,
    \begin{align}
        \sum_{k \in \N}\abs*{\left(1 - \gamma)^k\left(\pl(\abs{\C\left(\origin{a}\right)} = k \right) -\pl\left(\abs*{\C^{(n)}\left(\origin{a}\right)} = k \right)\right)} \leq \sum_{k \in \N}(1 - \gamma)^k < \infty.
    \end{align}
    By a similar reasoning than in Lemma \ref{lem:analyticity} we have for $n \in \N$, and also in the infinite volume limit, that
    \begin{align}
        \frac{\partial M^{(n)}(\lambda,\gamma,a)}{\partial \gamma} = \sum_{k \geq 1}k(1 - \gamma)^{k-1}\left(\pl\left(\abs*{\C^{(n)}\left(\origin{a}\right)} = k \right) \right)
        \label{derivative_mag_gamma}
    \end{align}
    In particular
    \begin{multline}
        \frac{\partial M(\lambda,\gamma,a)}{\partial \gamma} - \frac{\partial M^{(n)}(\lambda,\gamma,a)}{\partial \gamma}\\ = \sum_{k \geq 1}k (1 - \gamma)^{k - 1}\left(\pl(\abs*{\C\left(\origin{a}\right)} = k ) -\pl\left(\abs*{\C^{(n)}\left(\origin{a}\right)} = k \right) \right).
    \end{multline}
    Therefore, using dominated convergence, for \eqref{convergence_mag} and \eqref{convergence_mag_partial_gamma} it is enough to show for every $k \in \N$ that the following convergence holds
    \begin{align}
        \lim\limits_{n \to \infty}\pl(\abs{\C^{(n)}\left(\origin{a}\right)} = k) = \pl(\abs{\C\left(\origin{a}\right)} = k).
    \end{align}
    Let $k, l, n \in \N$ with $l < n$, decompose the event $\{\abs{\C\left(\origin{a}\right)} = k\}$ into $\{\abs{\C\left(\origin{a}\right)} = k, \C\left(\origin{a}\right) \subset \Lambda_l\}$ and $\{\abs{\C\left(\origin{a}\right)} = k, \C\left(\origin{a}\right) \not\subset \Lambda_l\}$ and decompose similarly for $\C^{(n)}\left(\origin{a}\right)$.
    Following a similar reasoning to that of Lemma \ref{lem:analyticity}, we get
    \begin{align}
    \label{restricted_cluster}
        &\pl(|\C(\origin{a})| = k+1, \C\left(\origin{a}\right) \subset \Lambda_l) \nonumber\\
        &\hspace{1cm}= \frac{\lambda^k}{k!}\int_{\Lambda_l^k} \left(\sum_{G\in \mathrm{Gr}_{k+1}}\left(\prod_{\left\{i,j\right\}\in E\left(G\right)}\connf\left(x_i,x_j\right)\right)\left(\prod_{\left\{i,j\right\}\not\in E\left(G\right)}\left(1-\connf\left(x_i,x_j\right)\right)\right)\right)\nonumber\\
        &\hspace{2.5cm}\times \exp\left(-\lambda \int_{\X} \left(1-\prod_{m\in\left\{0,1,\ldots,k\right\}}\left(1-\connf\left(y,x_m\right)\right)\right)\nu\left(\dd y\right)\right)\nu^{\otimes k}\left(\dd \Vec{x}_{\left[1,\ldots,k\right]}\right),
    \end{align}
    and also,
    \begin{align}
    \label{restricted_cluster_finvol}
        &\pl(|\C^{(n)}(\origin{a})| = k+1, \C\left(\origin{a}\right) \subset \Lambda_l) \nonumber\\
        &\hspace{0.5cm}= \frac{\lambda^k}{k!}\int_{\Lambda_l^k} \left(\sum_{G\in \mathrm{Gr}_{k+1}}\left(\prod_{\left\{i,j\right\}\in E\left(G\right)}\connf\left(x_i,x_j\right)\right)\left(\prod_{\left\{i,j\right\}\not\in E\left(G\right)}\left(1-\connf\left(x_i,x_j\right)\right)\right)\right)\nonumber\\
        &\hspace{1cm}\times \exp\left(-\lambda \int_\X \mathds{1}_{\Lambda_n}(y) \left(1-\prod_{m\in\left\{0,1,\ldots,k\right\}}\left(1-\connf\left(y,x_m\right)\right)\right)\nu\left(\dd y\right)\right)\nu^{\otimes k}\left(\dd \Vec{x}_{\left[1,\ldots,k\right]}\right).
    \end{align}
    Notice now that for fixed $x_1,x_2,\dots,x_k \in \Lambda_l^k$ and $y \in \X$ we have 
    \begin{align}
        \mathds{1}_{\Lambda_n}(y) \left(1-\prod_{m\in\left\{0,1,\ldots,k\right\}}\left(1-\connf\left(y,x_m\right)\right)\right) \leq \mathds{1}_{\Lambda_{n + 1}}(y) \left(1-\prod_{m\in\left\{0,1,\ldots,k\right\}}\left(1-\connf\left(y,x_m\right)\right)\right)
    \end{align}
    and
    \begin{align}
        \lim_{n \to \infty}\mathds{1}_{\Lambda_n}(y) \left(1-\prod_{m\in\left\{0,1,\ldots,k\right\}}\left(1-\connf\left(y,x_m\right)\right)\right) = 1-\prod_{m\in\left\{0,1,\ldots,k\right\}}\left(1-\connf\left(y,x_m\right)\right).
    \end{align}
    Therefore by monotone convergence theorem and continuity of the exponential function we get that 
    \begin{align}
        \lim_{n \to \infty}\exp\left(-\lambda \int_\X \mathds{1}_{\Lambda_n}(y) \left(1-\prod_{m\in\left\{0,1,\ldots,k\right\}}\left(1-\connf\left(y,x_m\right)\right)\right)\nu\left(\dd y\right)\right) \nonumber\\
        = \exp\left(-\lambda \int_\X \left(1-\prod_{m\in\left\{0,1,\ldots,k\right\}}\left(1-\connf\left(y,x_m\right)\right)\right)\nu\left(\dd y\right)\right).
    \end{align}
    We conclude now by another monotone convergence (decreasing in this case) that \eqref{restricted_cluster} is the limit as $n \to \infty$ of \eqref{restricted_cluster_finvol}.

    We show now that $\pl(|\C^{(n)}(\origin{a})| = k+1, \C\left(\origin{a}\right) \not\subset \Lambda_l)$ can be made arbitrarily small uniformly for all $n$ large, by taking $l$ fixed but  large enough. Fix $p, q, n \in \N$ with $0 < p < q < n$,  and let $E_n(p,q)$ be the event that in finite volume model in $\Lambda_n$, there is at least one point inside $\Lambda_p$ directly connected to a point outside of $\Lambda_q$. For every non-negative  integer-valued random variable $Y$ we have $\p(Y \geq 1) \leq \E(Y)$, combining this with Mecke formula it implies that
    \begin{align}
        \pl(E_n(p,q)) &\leq \lambda^2 \int_{\Lambda_p} \int_{\Lambda_l \setminus \Lambda_q} \varphi(x,y) \nu(\dd x)\nu(\dd y) \nonumber\\
        &\leq \lambda^2 \int_{\Lambda_p} \int_{\X \setminus \Lambda_q} \varphi(x,y) \nu(\dd x)\nu(\dd y).
    \end{align}
    This last estimate does not depend on $n$ and because of the assumption $\ref{Assump:BoundExpectedDegreeWithAverage}$ tends to $0$ as $q$ tends to infinity. Therefore $\pl(E_n(p,q))$ is small in $n$ for $q$ large enough.
    
    Now let $\epsilon > 0$ and take (hyper)-boxes $\Lambda_{p_1} \subset \cdots \subset \Lambda_{p_{k}}\subset \X$ such that the events $A_1, \dots, A_{k}$ have all probability at most $\epsilon$, uniformly in $n$, where
    \begin{align}
        A_1 &= \{\exists z\in\eta\setminus\Lambda_{p_1}\colon \origin{a} \sim {z} \text{ in } {\xi^{\origin{a}}_{\Lambda_n}}\}, \\
        A_m &= \{\exists w\in \eta\cap\Lambda_{p_{m-1}},z\in\eta\setminus\Lambda_{p_{m}}\colon  w \sim z \text{ in }{\xi^{\origin{a}}_{\Lambda_n}}\},
    \end{align}
    for $m = 2, \dots , k$. This is possible since the probability of $A_1$ tends to $0$ as $p_1 \to \infty$ because of our assumption over $a$. For $m \in \{2,\dots,k\}$ notice that $\pl(A_m) \leq \pl(E_n(p_{m-1},p_m))$ which tends to $0$ as $p_m \to \infty$. 
   
    Finally take $l = p_{k}$. If $\{\abs*{\C^{(n)}(\origin{a})} = k+1, \C\left(\origin{a}\right) \not\subset \Lambda_l \}$ occurs then there is a point outside $\Lambda_l$ connected to the origin in $k$ steps or fewer, and thus $\bigcup_{m = 1}^k A_m$ must occur. However, this last event has probability at most $ k \epsilon$ for all $n$ big enough, which finishes the proof of \eqref{convergence_mag} and \eqref{convergence_mag_partial_gamma}.

    For \eqref{convergence_mag_partial_lambda}, we first show that the following limit hold for all $k \in \N$
    \begin{align}
    \label{convergence_derivative_cluster_size}
        \lim_{n \to \infty}\frac{\dd \pl(\abs{\C^{(n)}\left(\origin{a}\right)} = k)}{\dd \lambda} = \frac{\dd \pl(\abs{\C\left(\origin{a}\right)} = k)}{\dd \lambda}.
    \end{align}
    In particular our proof shows (again) that these derivatives exist.
    
    Let $\lambda_{\min} > 0$, $I = (\lambda_{\min} , \infty)$ and define the function
    \begin{align}
        f \colon & \quad I \times \X^k \quad \longrightarrow \R_+ \nonumber\\
        &(\lambda,x_1,\dots,x_k) \longmapsto g(x_1,\dots,x_k) \exp(-\lambda h(x_1,\cdots,x_k)),
    \end{align}
    where
    \begin{align}
        g(x_1,\dots,x_k) = \left(\sum_{G\in \mathrm{Gr}_{k+1}}\left(\prod_{\left\{i,j\right\}\in E\left(G\right)}\connf\left(x_i,x_j\right)\right)\left(\prod_{\left\{i,j\right\}\not\in E\left(G\right)}\left(1-\connf\left(x_i,x_j\right)\right)\right)\right),
    \end{align}
    and
    \begin{align}
        h(x_1,\dots,x_k) = \int_{\X} \left(1-\prod_{m\in\left\{0,1,\ldots,k\right\}}\left(1-\connf\left(y,x_m\right)\right)\right)\nu\left(\dd y\right).
    \end{align}
    It is clear that $f$ is differentiable with respect to $\lambda$ and we have
    \begin{align}
        \frac{\partial f(\lambda,x_1,\dots,x_k)}{\partial \lambda} = -g(x_1,\dots,x_k)h(x_1,\dots,x_k)\exp(-\lambda h(x_1,\dots,x_k)).
    \end{align}
   Since we have $ 0 \leq \varphi \leq 1$, by the Weierstrass product inequality
    \begin{align}
         1-\prod_{m\in\left\{0,1,\ldots,k\right\}}\left(1-\connf\left(y,x_m\right)\right) \leq \sum_{m\in\left\{0,1,\ldots,k\right\}}\connf\left(y,x_m\right).
    \end{align}
    By assumption \ref{Assump:BoundExpectedDegreeWithAverage}, for $\nu$-almost all $x_1 = (\overline{x}_1, c_1) \in \X$ we have
    \begin{align}
       \int_{\X} \connf\left(y,x_1\right)\nu\left(\dd y\right) = \int_{\Ecal} D(c_1,b) \Pcal (\dd b) \leq \esssup_{c \in \Ecal} \int_{\Ecal} D(c,b) \Pcal (\dd b) = \norm{D}_{1,\infty}< \infty.
    \end{align}
    Therefore for almost all $(x_1,\dots,x_k) \in \X^k$,
    \begin{align}
        h(x_1,\dots,x_k) \leq (k + 1) \norm{D}_{1,\infty} < \infty.
    \end{align}
    This then implies that for all $\lambda \in I$ and almost all $(x_1,\cdots,x_k) \in \X^k$,
    \begin{align}
        \abs*{\frac{\partial f(\lambda,x_1,\dots,x_k)}{\partial \lambda}} \leq (k + 1)  \norm{D}_{1,\infty} g(x_1,\dots,x_k)\exp(-\lambda_{\min} h(x_1,\dots,x_k)).
    \end{align}
    This bound is an integrable function that doesn't depend on $\lambda$, and therefore by Leibniz integral rule we get that $\forall \lambda \in I$
    \begin{multline}
        \frac{\dd \pl(\abs{\C\left(\origin{a}\right)} = k + 1)}{\dd \lambda} \\
        = k \lambda^{- 1} \pl(\abs{\C\left(\origin{a}\right)} = k + 1) + \frac{\lambda^k}{k!}\int_{\X^k} \frac{\partial f(\lambda,x_1,\dots,x_k)}{\partial \lambda}\nu^{\otimes k}\left(\dd \Vec{x}_{\left[1,\ldots,k\right]}\right).
    \end{multline}
    Since $\lambda_{\min}$ was arbitrary this result holds for all $\lambda > 0$, in particular this shows (again) that $\pl(\abs{\C\left(\origin{a}\right)} = k + 1)$ is differentiable with respect to $\lambda$.

    Let $n \in \N$, and by doing the same procedure for $\C^{(n)}\left(\origin{a}\right)$ we get that  $\pl(\abs{\C^{(n)}\left(\origin{a}\right)} = k + 1)$  is differentiable with respect to $\lambda$. More precisely
    \begin{multline}
        \label{expression_derivative_cluter_size}
        \frac{\dd \pl(\abs{\C^{(n)}\left(\origin{a}\right)} = k + 1)}{\dd \lambda} \\
        = k \lambda^{- 1} \pl(\abs{\C^{(n)}\left(\origin{a}\right)} = k + 1) + \frac{\lambda^k}{k!}\int_{\X^k} \frac{\partial f_n (\lambda,x_1,\dots,x_k)}{\partial \lambda}\nu^{\otimes k}\left(\dd \Vec{x}_{\left[1,\ldots,k\right]}\right),
    \end{multline}
    where 
    \begin{align}
        f_n (\lambda,x_1,\dots,x_k) := g_n(x_1,\dots,x_k)\exp(-\lambda h_n(x_1,\cdots,x_k)),
    \end{align}
    with
    \begin{align}
        g_n(x_1,\dots,x_k) :=  \mathds{1}_{\Lambda_n^k}(x_1,\dots,x_k)g(x_1,\dots,x_k),
    \end{align}
    and
    \begin{align}
        h_n(x_1,\dots,x_k) := \int_{\X} \mathds{1}_{\Lambda_n}(y)\left(1-\prod_{m\in\left\{0,1,\ldots,k\right\}}\left(1-\connf\left(y,x_m\right)\right)\right)\nu\left(\dd y\right).
    \end{align}
    In particular
    \begin{align}
        \frac{\partial f_n(\lambda,x_1,\dots,x_k)}{\partial \lambda} = -g_n(x_1,\dots,x_k)h_n(x_1,\dots,x_k)\exp(-\lambda h_n(x_1,\dots,x_k)).
    \end{align}
    Therefore we get that
    \begin{align}
        &\abs*{\frac{\dd \pl(\abs{\C\left(\origin{a}\right)} = k + 1)}{\dd \lambda} - \frac{\dd \pl(\abs{\C^{(n)}\left(\origin{a}\right)} = k + 1)}{\dd \lambda}} \nonumber\\
        &\hspace{1cm} \leq k\lambda^{-1}\abs*{(\pl(\abs*{\C\left(\origin{a}\right)} = k + 1) -\pl\left(\abs*{\C^{(n)}\left(\origin{a}\right)} = k + 1 \right)}\nonumber\\
        &\hspace{2cm} + \frac{\lambda^k}{k!}\int_{\X^k}\abs*{\frac{\partial f(\lambda,x_1,\dots,x_k)}{\partial \lambda}- \frac{\partial f_n (\lambda,x_1,\dots,x_k)}{\partial \lambda}}\nu^{\otimes k}\left(\dd \Vec{x}_{\left[1,\ldots,k\right]}\right).
    \end{align}
    We already know that the first term on this bound goes to $0$ because of the preceding part of the proof. For the other one notice that $\frac{\partial f_n}{\partial \lambda}$ converge pointwise to $\frac{\partial f}{\partial \lambda}$ and
    \begin{align}
        \abs*{ \frac{\partial f_n (\lambda,x_1,\dots,x_k)}{\partial \lambda}}&\leq (k + 1) \norm{D}_{1,\infty}
        g_n(x_1,\dots,x_k)\exp(-\lambda h_n(x_1,\dots,x_k)) \label{eqn:derivativebound}\\
        & \leq (k + 1) \norm{D}_{1,\infty}
        g_n(x_1,\dots,x_k) \nonumber\\
        & \leq (k + 1) \norm{D}_{1,\infty}
        g(x_1,\dots,x_k).
    \end{align}
    This last bound is independent of $n$ and in order to get equation \eqref{convergence_derivative_cluster_size} we now show that this function $g$ is integrable which allows us to use dominated convergence theorem. For that purpose, given an element $G \in \mathrm{Gr}_{k+1}$ let $T(G)$ to be a spanning tree of $G$ rooted on $\origin{a}$. The choice of $T(G)$ is not unique in most cases, however the bound we get is independent of the choice and therefore any choice of spanning tree works for our purposes. We have then
    \begin{align}
        &\int_{\X^k} g(x_1,\dots,x_k) \nu^{\otimes k}\left(\dd \Vec{x}_{\left[1,\ldots,k\right]}\right)\nonumber\\
        &= \int_{\X^k} \left(\sum_{G\in \mathrm{Gr}_{k+1}}\left(\prod_{\left\{i,j\right\}\in E\left(G\right)}\connf\left(x_i,x_j\right)\right)\left(\prod_{\left\{i,j\right\}\not\in E\left(G\right)}\left(1-\connf\left(x_i,x_j\right)\right)\right)\right) \nu^{\otimes k}\left(\dd \Vec{x}_{\left[1,\ldots,k\right]}\right)\nonumber\\
        &\leq \sum_{G\in \mathrm{Gr}_{k+1}}\left(\int_{\X^k}\left(\prod_{\left\{i,j\right\}\in E\left(G\right)}\connf\left(x_i,x_j\right)\right) \nu^{\otimes k}\left(\dd \Vec{x}_{\left[1,\ldots,k\right]}\right) \right)\nonumber\\
        &\leq \sum_{G\in \mathrm{Gr}_{k+1}}\left(\int_{\X^k}\left(\prod_{\left\{i,j\right\}\in E\left(T(G)\right)}\connf\left(x_i,x_j\right)\right) \nu^{\otimes k}\left(\dd \Vec{x}_{\left[1,\ldots,k\right]}\right) \right).
    \end{align}
    The idea now is to bound this last quantity by starting the integration from the leaves of $T(G)$ and using the assumption \ref{Assump:BoundExpectedDegreeWithAverage}. Given a graph $G \in \mathrm{Gr}_{k+1}$ assume that in its associated spawning tree $T(G)$, $x_k$ is a leaf and $x_{k-1}$ is its parent (in the case that its not true picking any other leaf will work), and let's denote by $e$ the edge between them. Then
    \begin{align}
        &\int_{\X^k}\left(\prod_{\left\{i,j\right\}\in E\left(T(G)\right)}\connf\left(x_i,x_j\right)\right) \nu^{\otimes k}\left(\dd \Vec{x}_{\left[1,\ldots,k\right]}\right) \nonumber\\
        &= \int_{\X^{k-1}}\left(\prod_{\left\{i,j\right\}\in E\left(T(G)\right)\setminus \{e\}}\connf\left(x_i,x_j\right)\right)\left(\int_{\X}\connf\left(x_{k-1},x_k\right) \nu\left(\dd x_k \right)\right) \nu^{\otimes (k-1)}\left(\dd \Vec{x}_{\left[1,\ldots,(k-1)\right]}\right) \nonumber\\
        &\leq \norm{D}_{1,\infty}\int_{\X^{k-1}}\left(\prod_{\left\{i,j\right\}\in E\left(T(G)\right)\setminus \{e\}}\connf\left(x_i,x_j\right)\right) \nu^{\otimes (k-1)}\left(\dd \Vec{x}_{\left[1,\ldots,(k-1)\right]}\right).
    \end{align}
    We then iterate this procedure of isolating one edge (corresponding to a leaf) at a time and bounding its integral thanks to our integrability assumption to get the following bound
    \begin{align}
        \int_{\X^k}\left(\prod_{\left\{i,j\right\}\in E\left(T(G)\right)}\connf\left(x_i,x_j\right)\right) \nu^{\otimes k}\left(\dd \Vec{x}_{\left[1,\ldots,k\right]}\right) \leq \norm{D}_{1,\infty}^k.
    \end{align}
    Notice that this bound only depended on $G$ via the number of leaves in its spanning tree - or equivalently the number of vertices in $G$. For any $G \in \mathrm{G}_{k+1}$ and any spanning tree $T(G)$ chosen for it, we have $\abs*{E(T(G))} = k$.
    Therefore
    \begin{align}
        &\int_{\X^k} g(x_1,\dots,x_k) \nu^{\otimes k}\left(\dd \Vec{x}_{\left[1,\ldots,k\right]}\right) \leq \abs*{\mathrm{Gr}_{k + 1}} \norm{D}_{1,\infty}^k,
    \end{align}
    and we get the desired result \eqref{convergence_derivative_cluster_size}.

     We show now that for any $n \in \N$, and also for the infinite volume limit,
    \begin{align}
        \label{term_by_term_der}
        \frac{\partial M^{(n)}(\lambda,\gamma,a)}{\partial \lambda} = \sum_{k \in \N}(1 - \gamma)^{k + 1 }\frac{\dd \pl(\abs{\C^{(n)}\left(\origin{a}\right)} = k + 1)}{\dd \lambda}.
    \end{align}
    By using \eqref{eqn:derivativebound}, we find that 
    \begin{align}
        &\abs*{\frac{\lambda^k}{k!}\int_{\X^k} \frac{\partial f_n (\lambda,x_1,\dots,x_k)}{\partial \lambda}\nu^{\otimes k}\left(\dd \Vec{x}_{\left[1,\ldots,k\right]}\right)} \nonumber\\
        &\hspace{3cm}\leq \left(k+1\right)\norm{D}_{1,\infty}\mathbb{P}_{\lambda}\left(\abs{\C^{(n)}\left(\origin{a}\right)} = k + 1\right)\nonumber\\
        &\hspace{3cm}\leq \left(k+1\right)\norm{D}_{1,\infty}. \label{bound_integral_partial_fn}
    \end{align}
    Consider again $\lambda_{\min} > 0$ and $I = (\lambda_{\min},\infty)$, then fix $k_0 \in \N$ and $\lambda \in I$, using this last bound together with the expression \eqref{expression_derivative_cluter_size} we get
    \begin{align}
        &\abs*{\sum_{k = k_0}^\infty (1 - \gamma)^{k + 1}\frac{\dd \pl(\abs{\C^{(n)}\left(\origin{a}\right)} = k)}{\dd \lambda}} \nonumber\\
        &\leq \sum_{k = k_0}^\infty \abs*{\left(1 - \gamma\right)^{k + 1} \left( k \lambda^{- 1}\pl(\abs{\C^{(n)}\left(\origin{a}\right)} = k + 1) + \left(k+1\right)\norm{D}_{1,\infty} \right) } \nonumber\\
        & \leq \left(\lambda_{\min}^{-1} + \norm{D}_{1,\infty}\right)\sum_{k = k_0}^\infty (1 - \gamma)^{k + 1}(k + 1).
    \end{align}
    That last bound goes to $0$ as $k_0$ goes to $\infty$ uniformly on $\lambda$ (in $I$). Since $\lambda_{\min}$ was arbitrary we get that \eqref{term_by_term_der} holds for all $\lambda > 0$.
    
    To conclude notice that the bound \eqref{bound_integral_partial_fn} in independent of $n$, which allows us to get    
    \begin{multline}
        \abs*{\frac{\dd \pl(\abs{\C\left(\origin{a}\right)} = k + 1)}{\dd \lambda} - \frac{\dd \pl(\abs{\C^{(n)}\left(\origin{a}\right)} = k + 1)}{\dd \lambda}} \\\leq 2(k + 1)\left(\lambda^{-1} + \norm{D}_{1,\infty}\right),
    \end{multline}
    Since
    \begin{align}
        2\left(\lambda^{-1} + \norm{D}_{1,\infty}\right)\sum_{k = 0}^\infty (1 - \gamma)^{k + 1}(k + 1) < \infty,
    \end{align}
    applying dominated convergence then gives \eqref{convergence_mag_partial_lambda}.
     
\end{proof}

\begin{lemma}
    For $\lambda\geq 0$, $\gamma \in( 0,1)$, $a\in\Ecal$, $n \in \N$ and also in the infinite volume limit
\begin{equation}
    (1 - \gamma)\frac{\partial M^{(n)}}{\partial \gamma}(\lambda,\gamma,a) = \elg\left[\abs{\C^{(n)}\left(\origin{a}\right)}\mathds 1_{\{\nconn{\origin{a}}{\mathcal G}{\xi_{\Lambda_n}^{\origin{a}}}\}}\right],
    \label{connection_der_mag_chi}
\end{equation}
in particular,
\begin{equation}
     \lim_{\gamma \searrow 0}\frac{\partial M}{\partial \gamma}(\lambda,\gamma,a) = \sum_{k\in\N}k\pl\left(\abs*{\C(\origin{a})}=k\right) = \chi^{\mathrm{f}}_\lambda(a).
     \label{free_chi_der_mag}
 \end{equation}
\end{lemma}

\begin{proof}
    By the expression \eqref{derivative_mag_gamma} we have
    \begin{align}
        (1 - \gamma)\frac{\partial M^{(n)}}{\partial \gamma}(\lambda,\gamma,a) &= \sum_{k \geq 1}k(1 - \gamma)^{k}\left(\pl\left(\abs*{\C^{(n)}\left(\origin{a}\right)} = k \right) \right)\nonumber \\
        &= \elg\left[\abs{\C^{(n)}\left(\origin{a}\right)}\mathds 1_{\{\nconn{\origin{a}}{\mathcal G}{\xi_{\Lambda_n}^{\origin{a}}}\}}\right] ,
    \end{align}
    which gives \eqref{connection_der_mag_chi}. Then, again by \eqref{derivative_mag_gamma}
    \begin{equation}
        \frac{\partial M}{\partial \gamma}(\lambda,\gamma,a) = \sum_{k\geq 1}k\left(1-\gamma\right)^{k-1}\pl\left(\abs*{\C(\origin{a})}=k\right),
    \end{equation}
    and \eqref{free_chi_der_mag} follows by monotone convergence.
\end{proof}

\subsection{Magnetization Lower Bound}

We begin by deriving the following differential inequalities for the magnetization.

\begin{lemma}[Aizenman-Barsky differential inequalities on the magnetization] \label{lem:frst:diff_ineqs_magnetization} 
Let $\gamma \in (0,1)$, $\lambda>0$ and $a\in\Ecal$ such that $\int D(a,b)\Pcal(\dd b) < \infty$ then under assumption \ref{Assump:BoundExpectedDegreeWithAverage} we have
\begin{align}
 &(i) \label{first_dif_ineq} \ \frac{\partial M(\lambda,\gamma,a)}{\partial\lambda} \leq \frac{1-\gamma}{\lambda} M_{\sup}(\lambda,\gamma) \frac{\partial M(\lambda,\gamma,a)}{\partial\gamma}, \\
&(ii) \label{secon_dif_ineq} \  M(\lambda,\gamma,a) \leq \gamma \frac{\partial M(\lambda,\gamma,a)}{\partial\gamma} + M(\lambda,\gamma,a)^2 + \lambda M_{\sup}(\lambda,\gamma) \frac{\partial M(\lambda,\gamma,a)}{\partial\lambda}. 
\end{align}
\end{lemma}

\begin{proof}

These inequalities rely on the Russo formula which only holds in finite volume. Therefore we will derive the finite volume equivalent of the inequalities and then take limits thanks to Lemma \ref{fin_vol_lim}. One could also show that the Margulis-Russo formula holds, in this particular scenario in infinite volume, and do the same following proof directly on infinite volume.

Let $n \in \N$ and define for all $x \in \X$, $\xi'_n(x) := \xi[\mathscr{C}^{(n)}\left(\origin{a}\right) \cup \{x\}]$ where $\mathscr{C}^{(n)}\left(\origin{a}\right)$ is the cluster of $\origin{a}$ in $\xi_{\Lambda_n}^{\origin{a}}$. Then by Margulis-Russo formula and the lemma before we get that 
\begin{align}
\frac{\partial M^{(n)}(\lambda,\gamma,a)}{\partial \lambda} &= \int_{\Lambda_n} \el[\Id_{\{\conn {\origin{a}}  {\mathcal{G}} {\xi_{\Lambda_n}^{\origin{a},x}}\}} - \Id_{\{\conn {\origin{a}}  {\mathcal{G}} {\xi_{\Lambda_n}^{\origin{a}}}\}}] \nu(\mathrm{d}x) \nonumber \\
    &= \int_{\Lambda_n} \el[\Id_{\{\conn {\origin{a}}  {\mathcal{G}} {\xi_{\Lambda_n}^{\origin{a},x}}\} \cap \{\nconn {\origin{a}}  {\mathcal{G}} {\xi_{\Lambda_n}^{\origin{a}}}\}}] \nu(\mathrm{d}x) \nonumber \\
    &= \int_{\Lambda_n} \plg \left( x \in \piv{\conn{\origin{a}}{\mathcal G}{\xi^{\origin{a},x}}_{\Lambda_n}} \right)  \nu(\mathrm{d}x) \nonumber \\
    &= \int_{\Lambda_n} \plg\left( \conn{\origin{a}}{x}{\xi_{\Lambda_n}^{\origin{a},x}}, \nconn{\origin{a}}{\mathcal G}{\xi_{\Lambda_n}^{\origin{a}}}, \conn{x}{\mathcal G}{\xi_{\Lambda_n}^x} \right) \nu(\dd x) \nonumber \\
    & = \int_{\Lambda_n} \elg \Bigg[ \mathds 1_{\{\conn{\origin{a}}{x}{\xi_{\Lambda_n}^{\origin{a},x}}\}} \mathds 1_{\{\nconn{\origin{a}}{\mathcal G}{\xi_{\Lambda_n}^{\origin{a}}}\}}\nonumber\\
    & \hspace{3cm} \times\plg\left( \conn{x}{\mathcal G}{\xi_{\Lambda_n}^x \setminus \mathscr{C}^{(n)}\left(\origin{a}\right)} \mid \xi'_n(x) \right) \Bigg] \nu(\dd x)  \nonumber\\
    &\leq M^{(n)}_{\sup}(\lambda,\gamma) \int_{\Lambda_n} \plg \left( \conn{\origin{a}}{x}{\xi_{\Lambda_n}^{\origin{a},x}}, \nconn{\origin{a}}{\mathcal G}{\xi_{\Lambda_n}^{\origin{a}}} \right) \nu(\dd x) \nonumber \\
	&\leq  \frac{1}{\lambda}M_{\sup}(\lambda,\gamma)\elg\left[\abs{\C^{(n)}\left(\origin{a}\right)}\mathds 1_{\{\nconn{\origin{a}}{\mathcal G}{\xi_{\Lambda_n}^{\origin{a}}}\}}\right] \nonumber \\
	& = \frac{1-\gamma}{\lambda} M_{\sup}(\lambda,\gamma) \frac{\partial M^{(n)}(\lambda,\gamma,a)}{\partial\gamma}. 
\end{align}
Then we get equation \ref{first_dif_ineq} by applying Lemma \ref{fin_vol_lim}.

To prove the second inequality, let $n \in \N$ and note that 
\begin{equation}
    M^{(n)}(\lambda,\gamma,a) = \underbrace{\plg( |\C^{(n)}(\origin{a}) \cap \mathcal G| =1)}_{\text{(I)}} + \underbrace{\plg(|\C^{(n)}(\origin{a}) \cap \mathcal G| \geq 2)}_{\text{(II)}}.
\end{equation}
The first term is
\begin{equation}
    \text{(I)} = \sum_{k \geq 1} k \gamma(1-\gamma)^{k-1} \pl(|\C^{(n)}(\origin{a})|=k) = \gamma \frac{\partial M^{(n)}(\lambda,\gamma,a)}{\partial\gamma}.
\end{equation}
For an element $x \in \Lambda_n$, let $A_x$ be the event $\{x \in\mathcal G\} \cup\{ \conn{x}{\mathcal G}{\xi_{\Lambda_n}^x} \}$, the second term is then
\begin{equation}
    \text{(II)} = \underbrace{\plg(A_\origin{a} \circ A_\origin{a})}_{\text{(IIa)}} + \underbrace{\plg(|\C^{(n)}(\origin{a}) \cap \mathcal G| \geq 2, (A_\origin{a} \circ A_\origin{a})^c)}_{\text{(IIb)}},
\end{equation}
and we can bound $ \text{(IIa)} \leq \plg(A_\origin{a})^2 =  M^{(n)}(\lambda,\gamma,a)^2$ by the BK inequality. On the other hand,
\begin{equation}
\text{(IIb)} = \lambda \int_{\Lambda_n} \plg \Big( \conn{\origin{a}}{x}{\xi_{\Lambda_n}^{\origin{a},x}}, \nconn{\origin{a}}{\mathcal G}{\xi_{\Lambda_n}^\origin{a}}, A_x \circ A_x \text{ in } \xi_{\Lambda_n}^x[\eta \setminus \C^{(n)}(\origin{a})] \Big) \nu(\dd x).
\end{equation}
Since $\C^{(n)}(\origin{a})$ is a.s.~\emph{finite} we can partition via the size of it (notice this would also work in infintie volume since $\C(\origin{a})$ is a.s.~\emph{finite} if we assume that $\left\{\nconn{\origin{a}}{\mathcal G}{\xi^\origin{a}}\right\}$ holds):
\begin{align}
    \text{(IIb)} &= \lambda \sum_{m \in \N_0} \frac{1}{m!} \int_{\Lambda_n} \elg\left[ \sum_{\vec y \in \eta^{(m)}} \mathds 1_{\{ \C^{(n)}(\origin{a}) = \{\origin{a}, y_1, \ldots, y_m\}} \mathds 1_{\{\C^{(n)}(\origin{a}) \cap \mathcal G = \varnothing\}}\right. \nonumber\\
	&\hspace{5cm}\times\mathds 1_{\{\conn{\origin{a}}{x}{\xi_{\Lambda_n}^{\origin{a},x}}\}} \mathds 1_{\{A_x \circ A_x \text{ in } \xi_{\Lambda_n}^x[\eta \setminus \{y_1, \ldots, y_m\}] \}}  \Bigg] \nu(\dd x) \nonumber\\
	&= \sum_{m\in\N_0} \frac{\lambda^{m+1}}{m!} \int_{\Lambda_n \times \Lambda_n^m }  \plg \Big( \C\left(\origin{a},\xi_{\Lambda_n}^{\origin{a},\vec y_{[1,m]}}\right) = \{\origin{a}, y_1, \ldots, y_m\}, \C\left(\origin{a},\xi_{\Lambda_n}^{\origin{a},\vec y_{[1,m]}}\right) \cap \mathcal G = \varnothing, \nonumber\\
    &\hspace{5cm}	\conn{\origin{a}}{x}{\xi_{\Lambda_n}^{\origin{a}, \vec y, x}},  A_x\circ A_x\Big) \nu(\dd x)\nu^{\otimes m}( \dd \vec y_{[1,m]}) \nonumber\\
	&= \sum_{m\in\N_0} \frac{\lambda^{m+1}}{m!} \int_{\Lambda_n \times \Lambda_n^m }  \plg \Big(  \C\left(\origin{a},\xi_{\Lambda_n}^{\origin{a},\vec y_{[1,m]}}\right) = \{\origin{a}, y_1, \ldots, y_m\}, \C\left(\origin{a},\xi_{\Lambda_n}^{\origin{a},\vec y_{[1,m]}}\right) \cap \mathcal G = \varnothing, \nonumber\\
    &\hspace{4.5cm}\conn{\origin{a}}{x}{\xi_{\Lambda_n}^{\origin{a}, \vec y, x}}\Big) \times \plg\big(  A_x\circ A_x\big) \nu(\dd x) \nu^{\otimes m}(\dd \vec y_{[1,m]}). 
    \label{expression_IIb}
\end{align}
By the BK inequality, $\plg(A_x \circ A_x) \leq M_{\sup}(\lambda,\gamma) \plg(A_x)$, and so
\begin{align}
	\text{(IIb)} &\leq M_{\sup}(\lambda,\gamma) \sum_{m\in\N_0} \frac{\lambda^{m+1}}{m!} \int_{\Lambda_n \times \Lambda_n^m }  \plg \Big( \C\left(\origin{a},\xi_{\Lambda_n}^{\origin{a},\vec y_{[1,m]}}\right) = \{\origin{a}, y_1, \ldots, y_m\}, \nonumber\\
    &\hspace{2.5cm}\C\left(\origin{a},\xi_{\Lambda_n}^{\origin{a},\vec y_{[1,m]}}\right) \cap \mathcal G = \varnothing, \conn{\origin{a}}{x}{\xi_{\Lambda_n}^{\origin{a}, \vec y, x}},  A_x \Big) \nu(\dd x) \nu^{\otimes m}(\dd \vec y_{[1,m]}) \nonumber\\
		& = \lambda M_{\sup}(\lambda,\gamma) \int_{\Lambda_n} \plg \Big( \conn{\origin{a}}{x}{\xi_{\Lambda_n}^{\origin{a},x}}, \nconn{\origin{a}}{\mathcal G}{\xi_{\Lambda_n}^\origin{a}}, \conn{x}{\mathcal G}{\xi_{\Lambda_n}^{x}} \Big) \nu(\dd x) \nonumber\\
		&= \lambda M_{\sup}(\lambda,\gamma) \frac{\partial M^{(n)}(\lambda,\gamma,a)}{\partial \lambda},
\end{align}
then again we conclude by using Lemma \ref{fin_vol_lim} to get \ref{secon_dif_ineq} as required.

\end{proof}

\begin{remark}
    One thing that is interesting to note about the last proof is that, although it originated in the Bernoulli percolation on $\Z^d$ and we follow the ideas of \cite{AizBar87}. It seems that this strategy is better suited for continuous model since in our case the event $\big(  A_x\circ A_x\big)$ in \eqref{expression_IIb} is actually independent of the rest.
\end{remark}

The following lemma allows us to bound the the supremum norm $M_{\sup}(\lambda,\gamma)$ with the product of $M\left(\lambda,\gamma,a\right)$ and a given function. This complements the trivial reverse bound. Recall
\begin{equation}
    \Jcal_{\lambda,a} := \left(\sup_{k\geq 1}\left(\frac{\lambda}{2\left(1+\lambda\right)}\right)^k\essinf_{b\in\Ecal}D^{(k)}(a,b)\right)^{-1}.
\end{equation}
The similarity of $\Jcal_{\lambda,a}$ to $\Ical_{\lambda,a}$ means that many of the properties of $\Ical_{\lambda,a}$ also hold for $\Jcal_{\lambda,a}$. In particular, 
\begin{itemize}
    \item $\Jcal_{0,a}=\infty$ for all $a\in\Ecal$,
    \item $\lambda\mapsto\Jcal_{\lambda,a}$ (for all $a\in\Ecal$) and $\lambda\mapsto\esssup_{a\in\Ecal}\Jcal_{\lambda,a}$ are non-decreasing,
    \item \ref{Assump:AllReachablebySome} implies $\lambda\mapsto\esssup_{a\in\Ecal}\Jcal_{\lambda,a}$ is continuous and finite for $\lambda>0$.
\end{itemize}

\begin{lemma}
\label{lem:inf-pointwise_bound}
For all $\lambda>0$, $\gamma\in\left(0,1/2\right)$, and $a\in\Ecal$ such that $\int D\left(a,b\right)\Pcal\left(\dd b\right)<\infty$,
\begin{equation}
\label{eqn:inf-pointwise_bound}
    M_{\sup}(\lambda,\gamma) \leq \left(1+ \lambda \InfNorm{D}\Jcal_{\lambda,a}\right)M\left(\lambda,\gamma,a\right).
\end{equation}
\end{lemma}

\begin{proof}
The proof proceeds similarly to that of Lemma~\ref{thm:inf-pointwise bound susceptibility}. We begin by bounding $M\left(\lambda,\gamma,a\right)$ by considering the marks of the neighbours of $\origin{a}$. By a union bound over $\origin{a}\in\G$ and each neighbour being connected to $\G$, and then using Mecke's formula, we find
\begin{align}
    M\left(\lambda,\gamma,a\right) &\leq \gamma + \E_{\lambda,\gamma}\left[\sum_{\left(\xbar,b\right)\in\eta\colon \left(\xbar,b\right)\sim \origin{a}}\Id_{\left\{\conn{\left(\xbar,b\right)}{\G}{\xi^{\left(\xbar,b\right)}}\right\}}\right] \nonumber\\ &= \gamma + \lambda\int \plaga\left(\conn{\left(\xbar,b\right)}{\G}{\xi^{\left(\xbar,b\right)}}\right) \connf\left(\xbar;a,b\right) \dd \xbar \Pcal\left(\dd b\right) \nonumber\\
    &= \gamma + \lambda\int M\left(\lambda,\gamma,b\right) D\left(a,b\right) \Pcal\left(\dd b\right).
\end{align}
Then by taking the essential supremum over $a$ and splitting the integral using a supremum bound, we get
\begin{equation}
\label{eqn:MagnetInf_vs_Magnet1}
    M_{\sup}(\lambda,\gamma) \leq \gamma + \lambda\InfNorm{D}\norm*{M\left(\lambda,\gamma\right)}_1.
\end{equation}

Now we aim to get a lower bound on $M\left(\lambda,\gamma,a\right)$ in terms of $\norm*{M\left(\lambda,\gamma\right)}_1$. As in Lemma~\ref{thm:inf-pointwise bound susceptibility} we construct a new configuration $\left[\xi\right]^\origin{a}$ from an existing configuration $\xi$. In summary, we uniformly select a neighbour of $\origin{a}$ in $\xi^\origin{a}$ and let $\left[\xi\right]^\origin{a}$ have every vertex and edge present in $\xi$ and additionally the vertex $\origin{a}$ and the edge between $\origin{a}$ and $\left[x\right]$. If there are no such neighbours, we only add the vertex $\origin{a}$ and there are no further edges. Since $\int D\left(a,b\right)\Pcal\left(\dd b\right)<\infty$ this construction is $\pla$-almost surely defined.

Since the edges and vertices of $\left[\xi\right]^\origin{a}$ are a subset of the edges and vertices of $\xi^\origin{a}$, we have the bound
\begin{equation}
    M(\lambda,\gamma,a) \geq \pla\left(\conn{\origin{a}}{\G}{\left[\xi\right]^{\origin{a}}}\right) \geq \pla\left(\origin{a}\not\in\G,\conn{\origin{a}}{\G}{\left[\xi\right]^{\origin{a}}}\right).
\end{equation}
As in Lemma~\ref{thm:inf-pointwise bound susceptibility}, we can relate this probability to the magnetization of the vertex $\left[x\right]$. We once again find
\begin{align}
    \pla\left(\left[x\right]\text{ is defined}\right) &= 1-\exp\left(-\lambda\int D\left(a,b\right)\Pcal\left(\dd b\right)\right)\\
    \pla\left(\left[x\right]\in B\;\middle|\; \left[x\right]\text{ is defined} \right) &= \frac{\int_B\connf\left(u,\origin{a}\right)\nu\left(\dd u\right)}{\int_\X\connf\left(v,\origin{a}\right)\nu\left(\dd v\right)},
\end{align}
for any measurable set $B\subset \X$. If $\origin{a}$ is not a ghost vertex, then $\conn{\origin{a}}{\G}{\left[\xi\right]^{\origin{a}}}$ means that this connection must occur via $\left[x\right]$. The independence of $\left[x\right]$ from all the other vertices then means that for all $\xbar\in\Rd$ and $b\in\Ecal$
\begin{align}
    &\pla\left(\origin{a}\not\in\G,\conn{\origin{a}}{\G}{\left[\xi\right]^{\origin{a}}}\;\middle|\;\left[x\right]=\left(\xbar,b\right),\left[x\right]\text{ is defined}\right) \nonumber\\
    &\hspace{3cm}= \pla\left(\conn{\left(\xbar,b\right)}{\G}{\xi}\;\middle|\;\left[x\right]=\left(\xbar,b\right),\left[x\right]\text{ is defined}\right)\nonumber\\
    &\hspace{3cm} = \pla\left(\conn{\left(\xbar,b\right)}{\G}{\xi^{\left(\xbar,b\right)}}\right) \nonumber\\
    &\hspace{3cm}= M\left(\lambda,\gamma,b\right).
\end{align}

Therefore we can bound
\begin{align}
    M(\lambda,\gamma,a) &\geq \pla\left(\origin{a}\not\in\G,\conn{\origin{a}}{\G}{\left[\xi\right]^{\origin{a}}}\right)\nonumber\\
    &\geq \pla\left(\origin{a}\not\in\G,\conn{\origin{a}}{\G}{\left[\xi\right]^{\origin{a}}},\left[x\right]\text{ is defined}\right)\nonumber\\
    & = \left(1-\gamma\right)\int M(\lambda,\gamma,b)\frac{\connf\left(\left(\ubar,b\right),\origin{a}\right)}{\int\connf(v,\origin{a})\nu(\dd v)}\left(1- \e^{-\lambda\int\connf(v,\origin{a})\nu(\dd v)}\right)\dd \ubar\Pcal\left(\dd b\right)\nonumber\\
    & = \left(1-\gamma\right)\frac{\left(1- \e^{-\lambda\int D(a,b)\Pcal(\dd b)}\right)}{\int D(a,b)\Pcal(\dd b)}\int M(\lambda,\gamma,b) D(a,b)\Pcal\left(\dd b\right)\nonumber\\
    & \geq \left(1-\gamma\right)\frac{\lambda}{1+\lambda}\int M(\lambda,\gamma,b) D(a,b)\Pcal\left(\dd b\right). \label{eqn:Magnet_vs_Magnet1}
\end{align}

Now let $k\geq 1$ be a fixed integer. By iteratively applying \eqref{eqn:Magnet_vs_Magnet1} $k$ times, using an infimum bound to split an integral, and then optimising over $k$, we arrive at
\begin{equation}
    M\left(\lambda,\gamma,a\right) \geq \norm*{M\left(\lambda,\gamma\right)}_1\sup_{k\geq 1}\left(\left(1-\gamma\right)\frac{\lambda}{1+\lambda}\right)^k \essinf_{b\in\Ecal} D^{(k)}(a,b) \geq \frac{\norm*{M\left(\lambda,\gamma\right)}_1}{\Jcal_{\lambda,a}}.
\end{equation}
In this last inequality we have used $\gamma\leq 1/2$.

Then combining \eqref{eqn:MagnetInf_vs_Magnet1} with \eqref{eqn:Magnet_vs_Magnet1} gives
\begin{equation}
    M_{\sup}(\lambda,\gamma) \leq \gamma + \lambda \InfNorm{D}\Jcal_{\lambda,a} M\left(\lambda,\gamma,a\right).
\end{equation}
Finally using the bound $M\left(\lambda,\gamma,a\right) \geq \plaga\left(\origin{a}\in\G\right) = \gamma$ gives the result.
\end{proof}

\begin{corollary}[Lower bound on magnetization]
\label{lem:lowerboundmagnetization}
Let $\gamma\in\left(0,1/2\right)$, $a\in\Ecal$, and $\lambda>0$ be such that $\chi^{\mathrm{f}}_\lambda(a) = \infty$. Then then under the assumption \ref{Assump:BoundExpectedDegreeWithAverage} we have
\begin{equation}
    M(\lambda,\gamma,a) \geq \sqrt{\frac{\gamma}{1 + \left(1+ \lambda \InfNorm{D}\Jcal_{\lambda,a}\right)^2}}.
\end{equation}
\end{corollary}

\begin{proof}
Observe that if $\E_\lambda\left[\deg{\origin{a}}\right]=\infty$ then $\chi^{\mathrm{f}}_\lambda(a) < \infty$. Therefore the condition $\chi^{\mathrm{f}}_\lambda(a) = \infty$ and Mecke's formula implies that $\int D(a,b)\Pcal(\dd b) = \frac{1}{\lambda}\E_\lambda\left[\deg{\origin{a}}\right]<\infty$. We are therefore in a position to apply the Lemmas before. Inserting the first inequality  of Lemma~\ref{lem:frst:diff_ineqs_magnetization}, into the second yields
\begin{equation}
    \label{eq/combined_dif_ineq}
    M(\lambda,\gamma,a) \leq \gamma \frac{\partial M(\lambda,\gamma,a)}{\partial \gamma} + M(\lambda,\gamma,a)^2 + (1-\gamma) M_{\sup}(\lambda,\gamma)^2 \frac{\partial M(\lambda,\gamma,a)}{\partial\gamma}.
\end{equation}
Then by using Lemma~\ref{lem:inf-pointwise_bound} we get	
\begin{multline}
    M(\lambda,\gamma,a) \leq \gamma \frac{\partial M(\lambda,\gamma,a)}{\partial \gamma} + M(\lambda,\gamma,a)^2 \\+ (1-\gamma) \left(\left(1+ \lambda \InfNorm{D}\Jcal_{\lambda,a}\right)M\left(\lambda,\gamma,a\right)\right)^2 \frac{\partial M(\lambda,\gamma,a)}{\partial\gamma}.
    \label{bound_intermediary_cor_3.10}
\end{multline}

The fact that $\E_\lambda\left[\deg{\origin{a}}\right]<\infty$ implies that
\begin{align}
    \pl(|\C(\origin{a})| = 1) = \exp\left(-\lambda \int D(a,b)\Pcal(\dd b)\right) > 0,
    \label{eqn:no_neighbors}
\end{align}
in particular we have $\pl(|\C(\origin{a})| = \infty) < 1$, which makes the function $\gamma \mapsto M(\lambda,\gamma,a)$ strictly increasing by the expression \eqref{derivative_mag_gamma}. We denote the function $M_a \mapsto \gamma$ as the inverse to $\gamma \mapsto M(\lambda,\gamma,a)$ dividing then \eqref{bound_intermediary_cor_3.10} by $M(\lambda,\gamma,a)^2 \tfrac{\partial M(\lambda,\gamma,a)}{\partial\gamma}$ we obtain
\begin{multline}
    \frac{\partial}{\partial M_a} \left( \frac{\gamma}{M_a} \right) = \frac{1}{M_a} \frac{\partial\gamma}{\partial M_a} - \frac{1}{M_a^2} \gamma \leq (1-\gamma)\left(1+ \lambda \InfNorm{D}\Jcal_{\lambda,a}\right)^2 + \frac{\partial\gamma}{\partial M_a} \\
    \leq \left(1+ \lambda \InfNorm{D}\Jcal_{\lambda,a}\right)^2 + \frac{\partial\gamma}{\partial M_a}.
\end{multline}
Integrating from $0$ to $M_a$, we note that $\gamma(0)=0$, as well as $\lim_{M_a \to 0} \tfrac{\gamma(M_a)}{M_a} = \gamma'(0) = \tfrac{1}{M_a'(0+)} = \tfrac{1}{\chi_\lambda^{\mathrm{f}}(a)} =0$ (by Lemma~\ref{gamma_to_0} and $\chi^{\mathrm{f}}_\lambda(a) = \infty$) to get
\begin{equation}
    \frac{\gamma}{M_a} \leq \left(1+ \lambda \InfNorm{D}\Jcal_{\lambda,a}\right)^2 M_a + \gamma.
\end{equation}
Since $\gamma \leq M_a$ for all marks $a\in\Ecal$ we get
\begin{equation}
    \frac{\gamma}{M_a} \leq \left(1 + \left(1+ \lambda \InfNorm{D}\Jcal_{\lambda,a}\right)^2\right)M_a,
\end{equation}
which leads us to
\begin{equation}
    M(\lambda,\gamma,a) \geq \sqrt{\frac{\gamma}{1 + \left(1+ \lambda \InfNorm{D}\Jcal_{\lambda,a}\right)^2}}
\end{equation}
as required.
\end{proof}

\begin{remark}
    Although we are assuming \ref{Assump:BoundExpectedDegreeWithAverage} in order to get this result it would have been possible by just assuming the weaker assumption \eqref{Assumption:BoundedExpectedDegreeMinimal} mentioned in the remark after Lemma \ref{fin_vol_lim}. Indeed  to get the line \eqref{eq/combined_dif_ineq} we can combine the finite volume version of the differentials inequalities in the Lemma \ref{lem:frst:diff_ineqs_magnetization} and then take the infinite volume limit. This cancel the term of the derivative with respect to $\lambda$ and we can use the first two results of Lemma \ref{fin_vol_lim}, which, as mentioned in the previous remark only need the assumption \eqref{Assumption:BoundedExpectedDegreeMinimal}.

\end{remark}

\subsection{Magnetization Upper Bound}

Before proving our bound we need first to prove a differential inequality, which allows us to get an upper bound for the magnetization. It is adapted from a corresponding result for percolation on the lattice in \cite[Lemma A.4]{borgs2005random}. We define
\begin{equation}
    M_{\inf}(\lambda,\gamma) := \essinf_{a\in\Ecal}M(\lambda,\gamma,a).
\end{equation}

\begin{remark}
\label{obs_inf-sup}
Dividing \eqref{eqn:inf-pointwise_bound} by $M(\lambda,\gamma,a)$ on both sides and then taking the essential supremum over $a\in\Ecal$ on both sides gives
\begin{equation}
    \frac{M_{\sup}(\lambda,\gamma)}{M_{\inf}(\lambda,\gamma)} \leq 1+ \lambda \InfNorm{D}\esssup_{a\in\Ecal}\Jcal_{\lambda,a},
\end{equation}
for $\gamma\in\left(0,1/2\right)$.
\end{remark}

\begin{lemma}(Triangle differential inequality for the magnetization)
\label{tdim}
Under Assumptions \ref{Assump:BoundExpectedDegree} and \ref{TriangleCondition_Assumption}, if $\lambda \in [0,\lambda_T]$ and $\gamma \in (0,1/2)$, then for almost every $a \in \Ecal$
\begin{equation}
     M(\lambda,\gamma,a) \geq \frac{1}{2}M_{\inf}(\lambda,\gamma)^2\frac{\partial M(\lambda,\gamma,a)}{\partial \gamma}\kappa(\lambda,a) - M_{\sup}(\lambda,\gamma)^2\triangle_\lambda,
\end{equation}
where
\begin{equation}
    \kappa(\lambda,a) = \E_{\lambda}\left[\deg{\origin{a}}\right]^2 - (2\E_{\lambda}\left[\deg{\origin{a}}\right] + 1)\left(1+ \lambda \InfNorm{D}\esssup_{a\in\Ecal}\Jcal_{\lambda,a}\right)^2\triangle_\lambda.
\end{equation}
\end{lemma}

\begin{remark}
Notice that since $\frac{\partial M}{\partial \gamma}\geq 0$, this lemma is only useful if $\kappa(\lambda,a)>0$. Mecke's formula gives that $\E_{\lambda}\left[\deg{\origin{a}}\right] = \lambda\int\connf\left(\origin{a},x\right)\nu\left(\dd x\right) = \lambda\int D\left(a,b\right)\Pcal\left(\dd b\right)$, and so Assumption~\ref{TriangleCondition_Assumption} ensures that $\kappa\left(\lambda_T,a\right)>0$ for $\Pcal$-almost every $a\in\Ecal$.
\end{remark}

\begin{proof}
Recall $\origin{a} = (\zerobar,a)$, the point at the origin with fixed mark $a$. For $x \in \eta$ we define $\{\dconn x \G {\xi^x}\} :=  \{\conn{x}{\mathcal G}{\xi^x}\} \circ\{\conn{x}{\mathcal G}{\xi^x}\} $. Furthermore to avoid too heavy notations we denote $\mathscr{C}(x)$ for $\mathscr{C}(x,\xi^x)$. Let $\textrm{Piv}(\origin{a},\G,\xi)$ be the set of pivotal points for the connection of $\origin{a}$ to $\G$ in $\xi$. Now, still for $x \in \eta$, we define the following event $F_{a,x} = \{x \in \textrm{Piv} (\origin{a},\G,\xi^{\origin{a},x})\}\cap\{\dconn x \G {\xi^x} \}$. Let $F_a = \bigcup_{x \in \X} F_{a,x}$ and note that the union is disjoint. When $F_a$ occurs, $\origin{a}$ is connected to $\G$ which means that $M(\lambda,\gamma,a) \geq \plg(F_a)$, so in the rest of the proof we lower bound $\plg (F_a)$.

For $A \subset \eta$, we define the restricted ``ghost-free'' two-point function by
\begin{equation}
    \tlg^A(x,y) = \plg((\con x y, \ncon x \G) \textrm{ in } \xi^{x,y}[(\eta \setminus A)\cup \{x,y\}]).
\end{equation}
For $x$, $y \in \eta$ we say that $x$ and $y$ are connected in $\xi$ off $A$ and write $\{\conn x y \xi \text{ off } A\}$ for the event $\{\conn x y \xi[\eta \setminus A]\}$. To bound $\plg(F)$, we start by using the Mecke equation \eqref{eq:prelim:mecke_n}
\begin{multline}
    \plg(F_a) = \elg\left[\sum_{x \in \eta^\origin{a}} \mathds 1_{\{F_x \text{ in } \xi^{\origin{a}}\}}\right] = \elg\left[\sum_{x \in \eta^\origin{a}} \mathds 1_{\{x \in \text{Piv}(\origin{a},\G,\xi^{\origin{a}}\}}\mathds 1_{\{\dconn x \G \xi^{\origin{a}}\}}\right] \\
    = \plg(\dconn {\origin{a}} \G \xi^\origin{a}) + \lambda \int_{\X} \elg \left[\mathds 1_{\{F_{a,x} \text{ in } \xi^{\origin{a},x}\}}\right] \nu\left(\mathrm{d}x\right).
\end{multline}

For a point $x = (\Bar{x},b) \in \X$ we define the random set, the cluster of $x$ away from $\origin{a}$, by $\mathscr{C}^\origin{a}(x) := \mathscr{C}(x,\xi^{\origin{a},x}) \setminus (\mathscr{C}(\origin{a}) \cup \{x\})$. We condition on $\tilde \xi := \xi^{\origin{a},x}[\mathscr{C}^\origin{a}(x) \cup \{x\}]$, giving us information about the edges between $x$ and $\mathscr{C}^\origin{a}(x)$. We have then that
\begin{align}
    \elg \left[\mathds 1_{\{F_{a,x} \textrm{ in } \xi^{\origin{a},x}\}}\right] &= \elg \left[\mathds 1_{\{x \in \mathrm{Piv}(\origin{a},\G,\xi^{\origin{a},x})\}} \mathds 1_{\{\dconn x \G \xi^{\origin{a},x}\}}\right] \nonumber\\
    &= \elg\left[ \elg \left[\mathds 1_{\{x \in \mathrm{Piv}(\origin{a},\G,\xi^{\origin{a},x})\}} \mathds 1_{\{\dconn x \G \xi^{\origin{a},x}[\eta^x \setminus \mathscr{C}_\origin{a}]\}}| \tilde \xi\right]\right] \nonumber\\
    &=  \elg \left[ \mathds 1_{\{ \dconn x \G \xi^{\origin{a},x}\left[\eta^x \setminus \mathscr{C}_\origin{a}\right]\}} \elg \left[\mathds 1_{\{x \in \mathrm{Piv}(\origin{a},\G,\xi^{\origin{a},x})\}} | \tilde \xi\right]\right] \nonumber\\
    &= \elg \left[ \hat I (x) \tlg^{\C^\origin{a}( x)}(\origin{a},x) \right],
\end{align}
where $\hat I (x) :=  \mathds 1_{\{ \dconn x \G \xi^{\origin{a},x}[\eta^x \setminus \mathscr{C}(\origin{a})]\}}$ and by a consequence of the stopping-set lemma, see Corollary 3.4 in \cite{HeyHofLasMat19}, we have that 
\begin{multline}
    \tlg^{\cls x {\origin{a}}}(\origin{a},x) = \plg\left(\con {\origin{a}} x, \ncon {\origin{a}} \G \textrm{ in } \xi^{\origin{a},x}[(\eta^{\origin{a},x} \setminus \cls x {\origin{a}}) \cup \{x\}] | \tilde \xi\right) \\
    = \plg(\con {\origin{a}} x, \ncon {\origin{a}} \G \textrm{ in } \xi^{\origin{a},x}[\eta_{\langle \mathscr{C}^{\origin{a}}(x) \rangle} \cup \{\origin{a},x\}]).
\end{multline}

We use the following almost sure identities
\begin{align}
    \tlg^{\cls x {\origin{a}}}(\origin{a},x) &= \tlg(\origin{a},x) - \left(\tlg(\origin{a},x) - \tlg^{\cls x {\origin{a}}}(\origin{a},x)\right),\\
    \hat I (x) &= \mathds 1_{\{ \dconn x \G \xi^{x}\}} - \mathds 1_{\{ \dconn x \G \xi^{x}\}}\mathds 1_{\{\ndconn x \G  \xi^{\origin{a},x}[\eta^x \setminus \mathscr{C}_\origin{a}]\}}.
\end{align}
Therefore we get that
\begin{align}
    \label{bound:Fa}
    \plg(F_a) &= \underbrace{\plg (\dconn {\origin{a}} \G \xi^\origin{a})(\chi(\lambda,\gamma,a) +  M(\lambda,\gamma,a))}_{X_1} \nonumber\\
    &\hspace{2cm}- \underbrace{\lambda \int_{\X}\tlg(\origin{a},x) \elg \left[ \mathds 1_{\{ \dconn x \G \xi^{x}\}}\mathds 1_{\{\ndconn x \G  \xi^{\origin{a},x}[\eta^x \setminus \mathscr{C}_\origin{a}]\}}\right]\nu(\mathrm{d}x)}_{X_2} \nonumber\\
    &\hspace{4cm}- \underbrace{\lambda \int_{\X} \elg \left[\hat I (x) \left(\tlg(\origin{a},x) - \tlg^{\cls x {\origin{a}}}(\origin{a},x)\right)\right]\nu(\mathrm{d}x)}_{X_3}.
\end{align}
Having written this in the form $\plg(F_a) = X_1 - X_2 - X_3$, we now want to lower bound $X_1$ and upper bound $X_2$ and $X_3$.

\paragraph{Lower bound of $X_1$.}

Given two points $y$ and $z$ we define $E_{y,z}$ as the event that $y$ and $z$ are the only two neighbors of $\origin{a}$ in $\xi^{\origin{a}}$, and that $y$ is connected to $\G$, $z$ is connected to $\G$ but $y$ and $z$ are not connected in $\xi$. If $N(\origin{a})$ denotes the neighborhood of $\origin{a}$, we have
\begin{align}
    E_{y,z} = \{(N(\origin{a}) = \{y,z\}) \text{ in } \xi^\origin{a} \text{ and } (\con y \G, \con z \G, \ncon y z) \text{ in } \xi\}.
\end{align}

 It is clear that  $E_{y,z} \cap E_{y',z'} = \emptyset$ for $\{y,z\} \neq \{y',z'\}$. Furthermore if $E_{y,z}$ occurs for some couple of points then $\{ \dcon {\origin{a}} \G \}$ also occurs. So, again by the Mecke equation we get that  
\begin{align}
    &\plg(\dconn {\origin{a}} \G \xi^\origin{a}) \geq \elg\left[ \sum_{y,z \in \eta, \colon y \neq z}\mathds 1_{\{E_{y,z} \text{ in } \xi\}}\right] \nonumber\\
    &\hspace{0.5cm}= \lambda^2 \int_{\X^2} \plg \left( E_{y,z} \text{ in } \xi^{y,z}\right) \nu(\mathrm{d}y)\nu(\mathrm{d}z) \nonumber\\
    &\hspace{0.5cm}= \lambda^2 \exp\left(-\E_{\lambda}\left[\deg{\origin{a}}\right]\right) \int_{\X^2} \varphi(\origin{a},y) \varphi(\origin{a},z) \nonumber\\
    &\hspace{4cm}\times\plg \left((\con y \G, \con z \G, \ncon y z) \text{ in } \xi^{y,z} \right) \nu(\mathrm{d}y)\nu(\mathrm{d}z).
    \label{dcon}
\end{align}
Indeed we have that
\begin{align}
    &\plg((N(\origin{a}) = \{y,z\}) \text{ in } \xi^{\origin{a},y,z}) = \varphi(\origin{a},y)\varphi(\origin{a},z)\elg\left[\prod_{x \in \eta}(1 - \varphi(\origin{a},x))\right] \nonumber\\
    &\hspace{2cm}= \varphi(\origin{a},y)\varphi(\origin{a},z)\left(\sum_{m \geq 0} \frac{(-\lambda)^m}{m!} \int_{\X^m} \varphi(\origin{a},x_1)\ldots\varphi(\origin{a},x_m) \nu^{\otimes m}\left(\mathrm{d}x_{[1,m]}\right) \right)\nonumber\\
    &\hspace{2cm}= \varphi(\origin{a},y)\varphi(\origin{a},z)\exp\left({-\lambda\int_{\X}\varphi(\origin{a},x)\nu(\mathrm{d}x)}\right)\nonumber\\
    &\hspace{2cm}= \varphi(\origin{a},y)\varphi(\origin{a},z)\exp\left(-\E_{\lambda}\left[\deg{\origin{a}}\right]\right).
\end{align}

We then denote $W = W_{y,z}$ the event $\{(\con y \G, \con z \G, \ncon y z) \text{ in } \xi^{y,z} \}$ and we compute its probability by conditioning on $\mathscr{C}(y)$. Indeed,
\begin{align}
    \plg(W) &= \elg\left[\mathds 1_{\{(\con y \G, \con z \G, \ncon y z) \text{ in } \xi^{y,z} \}}\right] \nonumber\\
    &= \elg\left[\plg\left((\con y \G, \con z \G, \ncon y z) \text{ in } \xi^{y,z} \right)|\mathscr{C}(y) = A \right] \nonumber\\
    &= \elg\left[\plg\left((\conn y \G \xi^{y,z}, \conn z \G \xi^{y,z} \text{ off } \clx y {\xi^{y,z}}\right)|\mathscr{C}(y) = A\right].
    \label{wpro}
\end{align}
Now we can see why it was useful to condition on the cluster of $y$ because now we are computing the probability of two independent events, indeed we get
\begin{multline}
    \plg(W) \\= \elg[\plg((\conn y \G \xi^{y,z}|\mathscr{C}_y = A)\plg( \conn z \G \xi^{y,z} \text{ off } \clx y {\xi^{y,z}})|\mathscr{C}_y = A)] 
\end{multline}
Now we make use of the stopping set lemma to compute the second probability
\begin{align}
    &\plg( \conn z \G \xi^{y,z} \text{ off } \clx y {\xi^{y}})|\mathscr{C}(y) = A)\nonumber\\ 
    &\hspace{5cm}=  \plg( \conn z \G \xi^{z} \text{ off } \clx y {\xi^{y}})|\mathscr{C}(y) = A)\nonumber\\
    &\hspace{5cm}= \plg ( \conn z \G \xi^{z}[\eta_{\langle A \rangle } \cup \{z\}])\nonumber\\
    &\hspace{5cm}= \plg ( \conn z \G \xi^{z}) - \plg\left(\inconn z {A} \G {\xi^{z}}\right)\nonumber\\
    &\hspace{5cm}\geq M_{\inf}(\lambda,\gamma) - \plg\left(\inconn z {A} \G {\xi^{z}}\right),
    \label{offcon}
\end{align}
where in first  equality we are just removing the vertex $y$ from the configuration. This is permitted because we simultaneously know that $y \in \mathscr{C}(y)$ and that $z$ is connected to $\mathcal{G}$ without utilising vertices in that cluster. Then we used the stopping-set lemma, which leads us to
\begin{align}
    &\plg \left(\inconn z A \G {\xi^{z}}\right) \leq \elg \left[\sum_{w \in \eta} \mathds 1_{\{w \in \eta \setminus \eta_{\langle A \rangle}\}} \mathds 1_{\{\{\conn z w {\xi^z}\}\circ \{\conn w \G \xi\}\}} \right] \nonumber \\
    &\hspace{2cm}= \lambda \int_{\X}\plg \left(w \in \eta \setminus \eta_{\langle A \rangle}\right) \plg(\{\conn z w {\xi^{z,w}}\}\circ \{\conn w \G {\xi^w}\})\nu(\mathrm{d}w) \nonumber \\
    &\hspace{2cm}= \lambda \int_{\X}\Big(1 - \Bar{\varphi}(A,w)\Big) \plg(\{\conn z w {\xi^{z,w}}\}\circ \{\conn w \G {\xi^w}\})\nu(\mathrm{d}w) \nonumber\\
    &\hspace{2cm}\leq \lambda \int_{\X}\Big(1 - \Bar{\varphi}(A,w)\Big) \plg(\conn z w {\xi^z})\plg(\conn w \G \xi)\nu(\mathrm{d}w) \nonumber\\
    &\hspace{2cm}\leq \lambda M_{\sup}(\lambda,\gamma) \int_{\X}\Big(1 - \Bar{\varphi}(A,w)\Big) \tlam (w,z)\nu(\mathrm{d}w).
    \label{thcon}
\end{align}
We used the fact that for $w \in \eta$, it being in a $A$-thinning is independent of the events that $z$ is connected to $w$ and $w$ is connected to $\G$ in the second line. And then the BK inequality in the fourth line.

 Then by plugging the bound \eqref{thcon} into the equality \eqref{offcon} that we just computed, the resulting inequality into \eqref{wpro}, injecting everything into \eqref{dcon} leads to 2 terms. We then get
\begin{align}
    &\exp\Big(\E_{\lambda}\left[\deg{\origin{a}}\right]\Big)\plg(\dconn {\origin{a}} \G \xi^\origin{a}) \nonumber\\
    &\hspace{1cm}\geq \lambda^2 M_{\inf}(\lambda,\gamma)\int_{\X^2} \varphi(\origin{a},y)\varphi(\origin{a},z) \plg(\conn y \G {\xi^{y,z}})\nu(\mathrm{d}y)\nu(\mathrm{d}z) \nonumber\\
    &\hspace{2cm}- \lambda^3 M_{\sup}(\lambda,\gamma) \int_{\X^3} \varphi(\origin{a},y)\varphi(\origin{a},z)\tlam(w,z)\tlam(w,y)\nonumber\\
    &\hspace{7cm}\times\plg(\conn y \G {\xi^{y,z}})\nu(\mathrm{d}y)\nu(\mathrm{d}z)\nu(\mathrm{d}w) \label{lbdcon}.
\end{align}
Next notice that
\begin{equation}
    \plg(\conn y \G {\xi^{y,z}}) \geq \plg(\conn y \G {\xi^{y}}),
\end{equation}
and
\begin{align}
    &\plg(\conn y \G {\xi^{y,z}})\nonumber\\ 
    &\hspace{0.5cm}= \plg(\conn y \G {\xi^{y,z}}, z \in \text{Piv}(y,\G,\xi^{y,z})) + \plg(\conn y \G {\xi^{y,z}}, z \notin \text{Piv}(y,\G,\xi^{y,z})) \nonumber\\
    &\hspace{0.5cm}\leq \plg(\conn z \G {\xi^{z}}) + \plg(\conn y \G {\xi^{y}})\nonumber\\
    &\hspace{0.5cm}\leq 2M_{\sup}(\lambda,\gamma).
\end{align}
Plugging these into \eqref{lbdcon} this leads us to
\begin{align}
    &\exp\Big(\E_{\lambda}\left[\deg{\origin{a}}\right]\Big)\plg(\dconn {\origin{a}} \G \xi^\origin{a})\nonumber\\
    &\hspace{1cm}\geq \underbrace{\lambda^2 M_{\inf}(\lambda,\gamma)\int_{\X^2} \varphi(\origin{a},y)\varphi(\origin{a},z) \plg(\conn y \G {\xi^{y}})\nu(\mathrm{d}y)\nu(\mathrm{d}z)}_{Y_1}\nonumber \\
    &\hspace{2cm}- \underbrace{2\lambda^3 M_{\sup}(\lambda,\gamma)^2 \int_{\X^3} \varphi(\origin{a},y)\varphi(\origin{a},z)\tlam(w,z)\tlam(w,y)\nu(\mathrm{d}y)\nu(\mathrm{d}z)\nu(\mathrm{d}w)}_{Y_2}.
\end{align}
For $Y_1$, by applying an infimum bound on $\plg(\conn y \G {\xi^{y}})$ and factorising the integral we find
\begin{equation}
    Y_1 \geq M_{\inf}(\lambda,\gamma)^2\E_{\lambda}\left[\deg{\origin{a}}\right]^2.
\end{equation}
We now proceed to upper bound $Y_2$
\begin{multline}
    Y_2 \leq 2\lambda^3 M_{\sup}(\lambda,\gamma)^2\int\varphi(\origin{a},y)\varphi(\origin{a},z)\tlam(w,z)\tlam(w,y)\nu(\mathrm{d}y)\nu(\mathrm{d}z)\nu(\mathrm{d}w) \\
    \leq 2M_{\sup}(\lambda,\gamma)^2 \E_{\lambda}\left[\deg{\origin{a}}\right] \triangle_\lambda.
\end{multline}

So by neglecting the term $ M(\lambda,\gamma,a)$ in  $(\chi(\lambda,\gamma,a) +  M(\lambda,\gamma,a))$ from \eqref{bound:Fa} we finally get that
\begin{equation}
    X_1 \geq \chi(\lambda,\gamma,a)\E_{\lambda}\left[\deg{\origin{a}}\right]\left(\E_{\lambda}\left[\deg{\origin{a}}\right]M_{\inf}(\lambda,\gamma)^2 - 2 M_{\sup}(\lambda,\gamma)^2\triangle_\lambda\right).
\end{equation}

\paragraph{Upper bound of $X_2$.} We have by definition
\begin{equation}
    X_2 = \lambda \int_{\X}\tlg(\origin{a},x) \elg \left[ \mathds 1_{\{ \dconn x \G \xi^{x}\}}\mathds 1_{\{\ndconn x \G  \xi^{\origin{a},x} \text{ off } \mathscr{C}_\origin{a}\}}\right]\nu(\mathrm{d}x).
\end{equation}
If we want the double indicator not to vanish we need $x$ to be doubly connected to $\G$, but at least one one of these connections must happen through a point in $\mathscr{C}(\origin{a}) \setminus \{\origin{a}\}$ (this point cannot be $\origin{a}$ because the first indicator says that there is a double connection in $\xi^x$ and not $\xi^{\origin{a},x}$, so we have
\begin{align}
     &\elg\left[\mathds 1_{\{ \dconn x G \xi^{x}\}}\mathds 1_{\{\ndconn x G  \xi^{\origin{a},x} \text{ off } \mathscr{C}_\origin{a}\}}\right] \nonumber\\
     &\hspace{2cm}\leq \elg\left[\sum_{w \in \eta} \mathds 1_{\{\{\conn {\origin{a}} w \xi^\origin{a}\}\circ\{\conn x w {\xi^x}\}\circ\{\conn w \G \xi \}\circ\{\conn x \G {\xi^x[\eta \setminus \{w\}}]\}\}}\right] \nonumber\\
     &\hspace{2cm}\leq \lambda^2 M_{\sup}(\lambda,\gamma)^2 \int_\X \tlam(\origin{a},w) \tlam (x,w) \nu(\mathrm{d}w),
\end{align}
where in the last line we used the Mecke Formula and the BK inequality. Notice that in the second line, because we have disjoint occurrences, we can replace the event $\{\conn x \G {\xi^x}\}$ by $\{\conn x \G {\xi^x[\eta \setminus \{w\}}]\}$. After applying Mecke and BK we get 
\begin{equation}
\plg(\conn x \G {\xi^{x,w}[\eta \setminus \{w\}]}) = \plg(\conn x \G {\xi^x}) \leq M_{\sup}(\lambda,\gamma),
\end{equation}
and with the trivial bound $\tlg \leq \tlam$ we get
\begin{equation}
    X_2 \leq M_{\sup}(\lambda,\gamma)^2 \lambda^2\int_{\X^2} \tlam(\origin{a},x)\tlam(\origin{a},w) \tlam (x,w)\nu(\mathrm{d}x)\nu(\mathrm{d}w) \leq M_{\sup}(\lambda,\gamma)^2\triangle_\lambda .
\end{equation}

\paragraph{Upper bound of $X_3$.} Again, we start with the definition
\begin{equation}
    X_3 = \lambda \int_{\X} \elg \left[\hat I (x) \left(\tlg(\origin{a},x) - \tlg^{\cls x {\origin{a}}}(\origin{a},x)\right)\right]\nu(\mathrm{d}x).
\end{equation}
Now notice we have the easy bound
\begin{equation}
    \hat I (x) =  \mathds 1_{\{ \dconn x \G \xi^{\origin{a},x}[\eta^x \setminus \mathscr{C}_\origin{a}]\}}] \leq \mathds 1_{\{ \dconn x \G \xi^{x}\}},
\end{equation}
concerning the difference in the two-point functions we have
\begin{align}
    &\tlg(\origin{a},x) - \tlg^{\cls x {\origin{a}}}(\origin{a},x) \nonumber\\
    &\hspace{0.5cm}= \elg\left[\mathds 1_{\{(\con {\origin{a}} x, \ncon {\origin{a}} \G) \text{ in } \xi^{\origin{a},x}\}} - \mathds 1_{\{(\con {\origin{a}} x, \ncon {\origin{a}} \G) \text{ in } \xi^{\origin{a},x}[\eta_{\langle \cls x {\origin{a}} \rangle} \cup \{\origin{a},x\}]\}}\right]\nonumber\\
    &\hspace{0.5cm}\leq \plg\left((\con {\origin{a}} x, \ncon {\origin{a}} \G) \text{ in } \xi^{\origin{a},x}, \nconn {\origin{a}} x {\xi^{\origin{a},x}[\eta_{\langle \cls x {\origin{a}}}\cup \{\origin{a},x\}]}\right) \nonumber\\
    &\hspace{0.5cm}\leq \elg\left[\sum_{y \in \eta} \mathds 1_{\{y \in \eta_{\langle \cls x {\origin{a}} \rangle}\}}\mathds 1_{\{\nconn {\origin{a}} \G {\xi^{\origin{a},x}}\}}\mathds 1_{\{\conn {\origin{a}} y {\xi^\origin{a}}\} \circ \{\conn y x {\xi^x}\}}\right] \nonumber\\
    &\hspace{0.5cm}\leq \lambda \int_\X \plg(y \sim \cls x {\origin{a}}) \elg\left[\mathds 1_{\{\nconn {\origin{a}} \G {\xi^{\origin{a},x,y}}\}}\mathds 1_{\{\conn {\origin{a}} y {\xi^{\origin{a},y}}\} \circ \{\conn y x {\xi^{x,y}}\}}\right]\nu(\mathrm{d}x),
\end{align}
where in the first inequality we used the fact that if the indicator 
\begin{equation}
    \mathds 1_{\{(\con {\origin{a}} x, \ncon {\origin{a}} \G) \text{ in } \xi^{\origin{a},x}, \nconn {\origin{a}} x {\xi^{\origin{a},x}[\eta_{\langle \cls x {\origin{a}}}\cup \{\origin{a},x\}]}\}}
\end{equation} 
is $1$ then clearly the first indicator in the previous line is $1$ and the other is $0$. The next line we use an union bound and in the last one we used the Mecke formula and independence from $y$ being in the $\cls x {\origin{a}}$-thinning and the other events. In the following to be completely rigorous we should also do an approximation on finite boxes, not only for the application of the BK inequality but also because we will condition on sets that are almost-surely infinite.

Next we define $\mathscr{C}(\G,\xi)$ to be the points in $\xi$ that are connected to $\G$, notice that this random set is almost-surely infinite since $\G$ is almost-surely infinite and $\G \subset \mathscr{C}(\G,\xi)$. We then use again the Corollary 3.4 in \cite{HeyHofLasMat19} to obtain
\begin{align}
    &\elg\left[\mathds 1_{\{\nconn {\origin{a}} \G {\xi^{\origin{a},x,y}}\}}\mathds 1_{\{\conn {\origin{a}} y {\xi^{\origin{a},y}}\} \circ \{\conn y x {\xi^{x,y}}\}}\right] \nonumber\\
    &\leq \plg(\{\conn {\origin{a}} y {\xi^{\origin{a},y}} \text{ off } \mathscr{C}(\G,\xi^{\origin{a},x,y})\} \circ \{\conn y x {\xi^{x,y}} \text{ off } \mathscr{C}(\G,\xi^{\origin{a},x,y})\}) \nonumber\\
    &= \elg\left[\plg\left(\left\{\conn {\origin{a}} y {\xi^{\origin{a},y}\left[\eta^{\origin{a},y}_{\langle \mathscr{C}(\G,\xi^{\origin{a},x,y}) \rangle}\right]} \right\} \circ \left\{\conn y x {\xi^{x,y}\left[\eta^{x,y}_{\langle \mathscr{C}(\G,\xi^{\origin{a},x,y}) \rangle}\right]}\right\}\right)\right].
\end{align}
For a fixed $\mathscr{C}(\G,\xi^{\origin{a},x,y})$, the events considered are increasing so we can apply the BK inequality inside the expectation to get
\begin{align}
    &\elg\left[\mathds 1_{\{\nconn {\origin{a}} \G {\xi^{\origin{a},x,y}}\}}\mathds 1_{\{\conn {\origin{a}} y {\xi^{\origin{a},y}}\} \circ \{\conn y x {\xi^{x,y}}\}}\right] \nonumber\\
    &\hspace{1cm}\leq \elg\left[\plg\left(\left\{\conn {\origin{a}} y {\xi^{\origin{a},y}\left[\eta^{\origin{a},y}_{\langle \mathscr{C}(\G,\xi^{\origin{a},x,y}) \rangle}\right]} \right\}\right)\right.\nonumber\\
    &\hspace{6cm}\times\left.\plg\left(\left\{\conn y x {\xi^{x,y}\left[\eta^{x,y}_{\langle \mathscr{C}(\G,\xi^{\origin{a},x,y}) \rangle}\right]}\right\}\right)\right] \nonumber\\
    &\hspace{1cm}\leq \tlam(x,y)\tlg(\origin{a},y).
\end{align}
Plugging the bounds above in the definition of $X_3$ leads us to
\begin{equation}
    X_3 \leq \lambda^2 \int_{\X^2} \tlam(x,y)\tlg(\origin{a},y)\plg\left(\dconn x \G {\xi^x}, \conn y {\cls x {\origin{a}}} {\xi^{x,y}}\right) \nu(\mathrm{d}x)\nu(\mathrm{d}y).
\end{equation}
Moreover
\begin{align}
    &\plg(\dconn x \G {\xi^x}, \conn y {\cls x {\origin{a}}} {\xi^{x,y}}) \nonumber\\
    &\hspace{4cm}\leq \elg\left[\sum_{w \in \eta}\mathds 1_{\{\conn x \G {\xi^x} \} \circ \{\conn x w {\xi^x}\} \circ \{\conn w \G {\xi} \} \circ \{\conn y w {\xi^y   }\}}\right] \nonumber\\
    &\hspace{4cm}\leq \lambda^2 M_{\sup}(\lambda,\gamma)^2 \int_{\X^2} \tlam(x,w)\tlam(y,w)\nu(\mathrm{d}w),
\end{align}
and so
\begin{equation}
    X_3 \leq M_{\sup}(\lambda,\gamma)^2\chi(\lambda,\gamma,a)\triangle_\lambda.
\end{equation}

Using the fact that $\chi(\lambda,\gamma,a) = (1 - \gamma)\frac{\partial M(\lambda,\gamma,a)}{\partial \gamma}$ and putting together the bounds of $X_1, X_2$ and $X_3$ we get that
\begin{equation}
    M(\lambda,\gamma,a) \geq M_{\inf}(\lambda,\gamma)^2(1-\gamma)\frac{\partial M(\lambda,\gamma,a)}{\partial \gamma}\iota(\lambda,a) - M_{\sup}(\lambda,\gamma)^2\triangle_\lambda,
\end{equation}
where
\begin{multline}
    \iota(\lambda,a) := \inf_{\gamma \in (0,1/2)}\iota(\lambda,\gamma,a) \\:= \inf_{\gamma \in (0,1/2)}\left(\E_{\lambda}\left[\deg{\origin{a}}\right]^2 - (2\E_{\lambda}\left[\deg{\origin{a}}\right] + 1)\left(\frac{M_{\sup}(\lambda,\gamma)}{M_{\inf}(\lambda,\gamma)}\right)^2\triangle_\lambda\right).
\end{multline}
By then using Remark~\ref{obs_inf-sup}, we get that
\begin{multline}
    \iota(\lambda,a,\gamma) \geq \kappa(\lambda,a) \\:= \E_{\lambda}\left[\deg{\origin{a}}\right]^2 - (2\E_{\lambda}\left[\deg{\origin{a}}\right] + 1)\left(1+ \lambda \InfNorm{D}\esssup_{a\in\Ecal}\Jcal_{\lambda,a}\right)^2\triangle_\lambda.
\end{multline}
\end{proof}

\begin{corollary}(Upper bound for the magnetization)
\label{cor:square_root_bound_M}
Under Assumptions \ref{Assump:BoundExpectedDegree}, \ref{Assump:AllReachablebySome}, and \ref{TriangleCondition_Assumption}, there exists a constant $K$ such that for almost all marks $a \in \Ecal$ and all $\gamma\in\left(0,1/2\right)$, we have
\begin{equation}
    M(\lambda_T,\gamma,a) \leq \sqrt{K \gamma}.
\end{equation}
\end{corollary}

\begin{proof}

We start from the inequality from Lemma~\ref{tdim}. Then for $\lambda \in [0,\lambda_T]$ and $\gamma \in (0,1/2)$, and for almost all marks $a$, we have
\begin{equation}
\label{eq:lower-bound magn}
     M(\lambda,\gamma,a) \geq \frac{1}{2}M_{\inf}(\lambda,\gamma)^2\frac{\partial M(\lambda,\gamma,a)}{\partial \gamma}\kappa(\lambda,a) - M_{\sup}(\lambda,\gamma)^2\triangle_\lambda.
\end{equation}
We now use Lemma \ref{lem:inf-pointwise_bound} and the fact that $\gamma \leq 1/2$ to get that for almost all marks $a\in\Ecal$ we have
\begin{align}
    M_{\inf} (\lambda,\gamma) &\geq \essinf_{a \in \Ecal} \frac{M_{\sup}(\lambda,\gamma)}{1+ \lambda \InfNorm{D}\Jcal_{\lambda,a}} = \frac{M_{\sup}(\lambda,\gamma)}{1+ \lambda \InfNorm{D}\esssup_{a\in\Ecal}\Jcal_{\lambda,a}} = \frac{M(\lambda,\gamma,a)}{\cbar{\lambda}}, \\
    M_{\sup} (\lambda,\gamma) &\leq \cbar{\lambda}M(\lambda,\gamma,a),
\end{align}
where we denoted $\cbar{\lambda} := 1+ \lambda \InfNorm{D}\esssup_{a\in\Ecal}\Jcal_{\lambda,a}$. We know that $\cbar{\lambda} < \infty$ by Assumption \ref{Assump:AllReachablebySome}, and that $\cbar{\lambda}\geq 1$ by the non-negativity of $\connf$. We inject this now in \ref{eq:lower-bound magn} to get
\begin{equation}
     M(\lambda,\gamma,a) \geq \frac{1}{2}\left(\frac{M(\lambda,\gamma,a)}{\cbar{\lambda}}\right)^2\frac{\partial M(\lambda,\gamma,a)}{\partial \gamma}\kappa(\lambda,a) - \cbar{\lambda}^2M(\lambda,\gamma,a)^2\triangle_\lambda,
\end{equation}
By rearranging the terms using the fact that $ \frac{1}{2}\frac{\partial M(\lambda,\gamma,a)^2}{\partial \gamma} = M(\lambda,\gamma,a)\frac{\partial M(\lambda,\gamma,a)}{\partial \gamma} $  we get to the following differential inequality
\begin{equation}
    \frac{1}{2}\frac{\partial M(\lambda,\gamma,a)^2}{\partial \gamma} \leq \frac{4 \cbar{\lambda}^2}{\kappa(\lambda,a)}( 1 + \cbar{\lambda}^2 M(\lambda,\gamma,a)) \leq \frac{4 \cbar{\lambda}^2}{\kappa(\lambda,a)}( 1 + \cbar{\lambda}^2 )
\end{equation}
Now we integrate this inequality over the interval $[0,\gamma]$ and use Lemma~\ref{gamma_to_0} telling us that $M(\lambda,0+,a) = \theta_\lambda(a) = 0$  to get
\begin{align}
    M(\lambda,\gamma,a)^2 \leq \frac{8 \cbar{\lambda}^2(1 + \cbar{\lambda}^2)}{\kappa(\lambda,a)}\gamma \leq \frac{8 \cbar{\lambda}^2(1 + \cbar{\lambda}^2)}{\kappa(\lambda)}\gamma
\end{align}
where we denoted $\kappa(\lambda) := \essinf_{a \in \mathcal{E}}\kappa(\lambda,a)$. By Assumption \ref{TriangleCondition_Assumption} we know that $\kappa(\lambda_T) > 0$.
This gives us the result for $K := \frac{8 \cbar{\lambda_T}^2(1 + \cbar{\lambda_T}^2)}{\kappa(\lambda_T)}$.
Notices that we used $\gamma \leq 1/2$ in the proof but this result also holds for $\gamma \geq 1/2$ up to eventually increasing the constant $K$.
\end{proof}

\section{Percolation Proofs}
\label{sec:PercolationProofs}

In this section we use our magnetization results to derive bounds on the percolation function.

\subsection{Percolation Lower Bound}

For each $a\in\Ecal$ define
\begin{equation}
    \lambda_T(a) := \inf\left\{\lambda>0\colon \chi_\lambda(a)= \infty\right\}.
\end{equation}

\begin{lemma}
\label{lem:EqualCriticalIntensities}
    For $\Pcal$-almost every $a\in\Ecal$ such that the conditions $\int D(a,b)\Pcal(\dd b)<\infty$ and $\sup_{k\geq 1}\essinf_{b\in\Ecal} D^{(k)}(a,b)>0$ both hold,
    \begin{equation}
        \lambda_T(a) = \lambda_T^{(1)}.
    \end{equation}
\end{lemma}

\begin{proof}
    First observe that $\norm*{\chi_\lambda}_p<\infty$ implies that $\chi_\lambda(a)<\infty$ for $\Pcal$-almost every $a\in\Ecal$. Therefore 
    \begin{equation}
        \lambda_T(a) \geq \lambda_T^{(p)}
    \end{equation}
    for all $p\in\left[1,\infty\right]$ and $\Pcal$-almost every $a\in\Ecal$. In particular, this holds for $p=1$ which maximises $\lambda_T^{(p)}$.

    Recall we proved in \eqref{eqn:Susp_vs_Suspk} that
    \begin{equation}
        \chi_\lambda(a) \geq \norm*{\chi_\lambda}_1 \sup_{k\geq 1}\left(\frac{\lambda}{1+\lambda}\right)^k \essinf_{b\in\Ecal}D^{(k)}(a,b).
    \end{equation}
    Note that if $\sup_{k\geq 1}\essinf_{b\in\Ecal} D^{(k)}(a,b)>0$ holds, then for $\lambda>0$ we also have the inequality $\sup_{k\geq 1}\left(\frac{\lambda}{1+\lambda}\right)^k \essinf_{b\in\Ecal}D^{(k)}(a,b)>0$. This means that $\lambda>0$ and $\norm*{\chi_\lambda}_1=\infty$ implies $\chi_\lambda(a)=\infty$. On the other hand, if $\lambda =0$ then $\norm*{\chi_\lambda}_1 = \chi_\lambda(a)=1$ for all $a\in\Ecal$. Therefore $\lambda_T(a)\leq \lambda_T^{(1)}$ if $a$ satisfies $\sup_{k\geq 1}\essinf_{b\in\Ecal} D^{(k)}(a,b)>0$. This proves the result.
\end{proof}

\begin{lemma}
\label{lem:CriticalSusceptibility}
Suppose Assumption~\ref{Assump:BoundExpectedDegree} holds. Then for all $\Pcal$-almost every $a\in\Ecal$ such that $\sup_{k\geq 1}\essinf_{b\in\Ecal} D^{(k)}(a,b)>0$,
\begin{equation}
    \chi_{\lambda_T(a)}(a)=\infty.
\end{equation}
\end{lemma}

\begin{proof}
By a `method of generations' approach like we employed in Lemma~\ref{lem:nontriviallambdaO}, we can show
\begin{equation}
    \int T_\lambda(a,b)\Pcal(\dd b) \leq \esssup_{b\in\Ecal}D(a,b)\sum^\infty_{k=1}\lambda^{k-1}\InfNorm{D}^{k-1}.
\end{equation}
Therefore Assumption~\ref{Assump:BoundExpectedDegree} implies that for $\Pcal$-almost every $a\in\Ecal$ and $\lambda>0$ sufficiently small, $\chi_\lambda(a)<\infty$. Therefore $\lambda_T(a)>0$ for $\Pcal$-almost every $a\in\Ecal$.

We first address the case where $\lambda_T(a)=\infty$. If $D(a,b)=0$ for $\Pcal$-almost every $b\in\Ecal$, then $\essinf_{b\in\Ecal} D^{(k)}(a,b)=0$ for all $k\geq 1$. This contradicts this lemma's assumptions, and so there exist $\varepsilon>0$ and $\Pcal$-positive set $B\subset \Ecal$ such that $D(a,b)>\varepsilon$ for all $b\in B$. We can use this with Mecke's formula to bound the expected degree of $\origin{a}$. For $\lambda<\infty$,
\begin{equation}
    \E_\lambda\left[\deg \origin{a}\right] = \lambda\int D(a,b)\Pcal\left(\dd b\right) > \lambda \varepsilon \Pcal\left(B\right).
\end{equation}
Since $\chi_\lambda(a)>\E_\lambda\left[\deg \origin{a}\right]$, we then have $\chi_\lambda(a)\to \infty$ as $\lambda\to\infty$. In particular, if $\lambda_T(a)=\infty$, then $\chi_{\lambda_T(a)}(a)=\infty$.

Now let us suppose $\lambda_T(a)<\infty$. Recall $\chi_\lambda(a) = 1 + \lambda\int T_\lambda(a,b)\Pcal(\dd b)$. We therefore want to prove $\int T_{\lambda_T(a)}(a,b)\Pcal(\dd b) = \infty$, and we do so in a similar manner to how we proved $\OpNorm{\OptlamT }=\infty$ in the proof of Theorem~\ref{thm:OperatorLowerBound}. Since $\int T_\lambda(a,b)\Pcal(\dd b)=\infty$ for $\lambda>\lambda_T(a)$, we only need to show that the function $\lambda\mapsto 1/\int T_\lambda(a,b)\Pcal(\dd b)$ is continuous at $\lambda=\lambda_T(a)$. We do this by showing that a family of truncated versions (which converge pointwise to $1/\int T_\lambda(a,b)\Pcal(\dd b)$) is equicontinuous.

Recall $t^{(n)}_\lambda(a,b) := \esssup_{\xbar\in\left[-n,n\right]^d}\int \tau^{(n)}_\lambda\left(\left(\xbar,a\right),\left(\ybar,b\right)\right) \dd \ybar$, and consider the functions $\lambda\mapsto \int t^{(n)}_\lambda(a,b)\Pcal(\dd b)$. Since $h\mapsto \int\left(\esssup_{\xbar\in\left[-n,n\right]^d}\int \abs*{h\left(\left(\xbar,a\right),\left(\ybar,b\right)\right)} \dd \ybar\right)\Pcal(\dd b)$ satisfies the triangle inequality, we can use the Margulis-Russo formula and BK inequality to get
\begin{equation}
    \limsup_{\varepsilon\to0}\frac{1}{\varepsilon}\left(\int t^{(n)}_{\lambda+\varepsilon}(a,b)\Pcal(\dd b) - \int t^{(n)}_\lambda(a,b)\Pcal(\dd b)\right) \leq \int t^{(n)}_\lambda(a,b)t^{(n)}_\lambda(b,c)\Pcal^{\otimes2}(\dd b,\dd c).
\end{equation}
We can then split this integral in two by using a supremum bound on $b$, and then using the inequality \eqref{eqn:tauinftyboundpointwise} from Remark~\ref{obs:tauinftyboundpointwise} to get
\begin{multline}
    \limsup_{\varepsilon\to0}\frac{1}{\varepsilon}\left(\int t^{(n)}_{\lambda+\varepsilon}(a,b)\Pcal(\dd b) - \int t^{(n)}_\lambda(a,b)\Pcal(\dd b)\right) \leq \OneNorm{ T_\lambda}\int t^{(n)}_\lambda(a,b)\Pcal(\dd b) \\\leq \left(1+\lambda\InfNorm{D}\Ical_{\lambda,a}\right)\left(\int t^{(n)}_\lambda(a,b)\Pcal(\dd b)\right)^2.
\end{multline}
Therefore 
\begin{equation}
    \liminf_{\varepsilon\to 0} \frac{1}{\varepsilon}\left(\frac{1}{\int t^{(n)}_{\lambda+\varepsilon}(a,b)\Pcal(\dd b)} - \frac{1}{\int t^{(n)}_\lambda(a,b)\Pcal(\dd b)}\right) \geq -\left(1+\lambda\InfNorm{D}\Ical_{\lambda,a}\right).
\end{equation}
Note that $\InfNorm{D}<\infty$ by Assumption~\ref{Assump:BoundExpectedDegree}, and $\Ical_{\lambda,a}<\infty$ holds by having the inequality $\sup_{k\geq 1}\essinf_{b\in\Ecal} D^{(k)}(a,b)>0$. This bound then implies that the functions $\lambda\mapsto \int t^{(n)}_\lambda(a,b)\Pcal(\dd b)$ are equicontinuous on $\left[0,\lambda^*\right]$ for any finite $\lambda^*>0$. Therefore the non-increasing pointwise limit $\lambda\mapsto 1/\int T_\lambda(a,b)\Pcal(\dd b)$ is continuous everywhere and we have $\int T_{\lambda_T(a)}(0;a,b)\Pcal(\dd b) = \infty$ as required.
\end{proof}

\begin{proposition}
\label{thm:PercolationBoundPointwise}
For all $a\in\Ecal$ such that $\esssup_{b\in\Ecal} D(a,b)<\infty$, $\theta_{\lambda_T(a)}(a)=0$ implies that for all $\lambda\geq0$ we have
\begin{equation}
\label{eqn:PercolationBoundPointwise}
    \theta_\lambda(a) \geq \frac{1}{1+ \lambda \InfNorm{D}\Jcal_{\lambda_T(a),a}}\frac{1}{2\lambda}\left(\lambda - \lambda_T(a)\right)_+.
\end{equation}
\end{proposition}

\begin{proof}
First note that if $\lambda=\infty$ or $\InfNorm{D}=\infty$ or $\sup_{k\geq 1}\essinf_{b\in\Ecal} D^{(k)}(a,b)=0$ or $\lambda_T(a)=0$ then the result is trivial, so let us assume $\lambda<\infty$, $\InfNorm{D}<\infty$ (i.e. \ref{Assump:BoundExpectedDegree}), $\lambda_T(a)>0$, and $\sup_{k\geq 1}\essinf_{b\in\Ecal} D^{(k)}(a,b)>0$. In particular, this implies that $\Jcal_{\lambda_T(a),a}<\infty$.

Now note that there is nothing to prove for $\lambda \leq \lambda_T(a)$, so we aim to prove \eqref{eqn:PercolationBoundPointwise} for $\lambda=\lambda^*>\lambda_T(a)$. We start our reasoning from the second equation of Lemma~\ref{lem:frst:diff_ineqs_magnetization}, namely \begin{align}
      M(\lambda,\gamma,a) \leq \gamma \frac{\partial M(\lambda,\gamma,a)}{\partial\gamma} + M(\lambda,\gamma,a)^2 + \lambda M_{\sup}(\lambda,\gamma) \frac{\partial M(\lambda,\gamma,a)}{\partial\lambda},
\end{align}
which by multiplying by $\frac{1}{\gamma M}$ gives
\begin{align}
     \frac{1}{\gamma} \leq \frac{1}{M(\lambda,\gamma,a)}\frac{\partial M(\lambda,\gamma,a)}{\partial\gamma} + \frac{M(\lambda,\gamma,a)}{\gamma} + \frac{\lambda}{\gamma}\frac{M_{\sup}(\lambda,\gamma)}{M(\lambda,\gamma,a)}\frac{\partial M(\lambda,\gamma,a)}{\partial\lambda},
\end{align} 
Then using Lemma~\ref{lem:inf-pointwise_bound} leads us to
\begin{align}
\label{eqn:diffinequal1}
     \frac{1}{\gamma} \leq \frac{1}{M(\lambda,\gamma,a)}\frac{\partial M(\lambda,\gamma,a)}{\partial\gamma} + \frac{M(\lambda,\gamma,a)}{\gamma} + \frac{\lambda}{\gamma}\left(1+ \lambda \InfNorm{D}\Jcal_{\lambda,a}\right)\frac{\partial M(\lambda,\gamma,a)}{\partial\lambda},
\end{align} 
Since $\Jcal_{\lambda,a}$ is non-increasing in $\lambda$ for all $a\in\Ecal$, we have 
\begin{equation}
    c_a := 1+ \lambda^* \InfNorm{D}\Jcal_{\lambda_T(a),a} \geq 1+ \lambda \InfNorm{D}\Jcal_{\lambda,a}.
\end{equation}
Since $c_a\geq 1$, the differential inequality \eqref{eqn:diffinequal1} implies
\begin{equation}
\label{eqn:integrandTerms}
    0 \leq \frac{\partial \log M(\lambda,\gamma,a)}{\partial\gamma} + \frac{1}{\gamma}\frac{\partial }{\partial\lambda}\left(\lambda c_a M(\lambda,\gamma,a) - \lambda \right),
\end{equation}
for $\gamma\in\left[0,1/2\right]$ and $\lambda\in\left[0,\lambda^*\right]$.

We now integrate the two terms on the right hand side over $\left(\gamma,\lambda\right)\in\left[\gamma_1,\gamma_2\right]\times\left[\lambda_T(a),\lambda^*\right]$, where $0<\gamma_1\leq \gamma_2\leq 1/2$. First we write the whole expression as an integral over $\lambda$ and then $\gamma$ (we can make this choice because the integrand in non-negative and we use Tonelli's Theorem). Then we split the integral into two integrals, each over one term in \eqref{eqn:integrandTerms}. Then since the first term's integrand is clearly non-negative we can use Tonelli's Theorem to exchange the order of the integrals and perform the $\gamma$ integral first for this term.

Recall that for the first term first in \eqref{eqn:integrandTerms} we first integrate over $\gamma$. We then use that $M(\lambda,\gamma,a)$ is increasing in $\lambda$ for all $\gamma$ and $a$ to get
\begin{align}
    &\int^{\lambda^*}_{\lambda_T(a)}\int^{\gamma_2}_{\gamma_1}\frac{\partial \log M(\lambda,\gamma,a)}{\partial\gamma}\dd \gamma \dd \lambda = \int^{\lambda^*}_{\lambda_T(a)} \log  M(\lambda,\gamma_2,a) - \log M(\lambda,\gamma_1,a)\dd \lambda \nonumber\\
    &\hspace{5cm}\leq \int^{\lambda^*}_{\lambda_T(a)} \log  M(\lambda^*,\gamma_2,a) - \log M(\lambda_T(a),\gamma_1,a)\dd \lambda \nonumber\\
    &\hspace{5cm}=\left(\lambda^*-\lambda_T(a)\right)\left(\log  M(\lambda^*,\gamma_2,a) - \log M(\lambda_T(a),\gamma_1,a)\right).
\end{align}
For the second term we first integrate over $\lambda$ and use that $M(\lambda,\gamma,a)$ is non-negative and increasing in $\gamma$ to get
\begin{align}
    &\int^{\gamma_2}_{\gamma_1}\int^{\lambda^*}_{\lambda_T(a)}\frac{1}{\gamma}\frac{\partial}{\partial \lambda}\left(\lambda c_a M(\lambda,\gamma,a) - \lambda\right)\dd \lambda \dd \gamma \nonumber\\
    & \hspace{3cm}= \int^{\gamma_2}_{\gamma_1}\frac{1}{\gamma}\left(\lambda^* c_a M(\lambda^*,\gamma,a) - \lambda^* - \lambda_T(a) c_a M(\lambda_T(a),\gamma,a) + \lambda_T(a)\right) \dd \gamma \nonumber\\
    & \hspace{3cm}\leq \int^{\gamma_2}_{\gamma_1}\frac{1}{\gamma}\left(\lambda^* c_a M(\lambda^*,\gamma_2,a) - \lambda^* + \lambda_T(a)\right) \dd \gamma \nonumber\\
    & \hspace{3cm}= \left(\lambda^* c_a M(\lambda^*,\gamma_2,a) - \lambda^* + \lambda_T(a)\right)\log\frac{\gamma_2}{\gamma_1}.
\end{align}
In summary, with some rearranging we have
\begin{equation}
\label{eqn:RearrageMagnetization}
    0\leq \left(\lambda^*-\lambda_T(a)\right)\left(\frac{\log M(\lambda^*,\gamma_2,a)}{\log\gamma_2-\log\gamma_1} - \frac{\log M(\lambda_T(a),\gamma_1,a)}{\log\gamma_2-\log\gamma_1}\right) + \lambda^* c_a M(\lambda^*,\gamma_2,a) - \lambda^* + \lambda_T(a).
\end{equation}

Since $\theta_{\lambda_T(a)}(a)=0$, the cluster of $\origin{a}$ is finite almost surely. Therefore $\chi^{\mathrm{f}}_{\lambda_T(a)}(a) = \chi_{\lambda_T(a)}(a) = \infty$ from Lemma~\ref{lem:CriticalSusceptibility}. Therefore we can use Corollary~\ref{lem:lowerboundmagnetization} to show
\begin{equation}
    M(\lambda_T(a),\gamma_1,a) \geq \sqrt{\frac{\gamma_1}{1+2c_a^2}}
\end{equation}
for $\gamma_1\leq 1/2$. Therefore
\begin{equation}
    - \frac{\log M(\lambda_T(a),\gamma_1,a)}{\log\gamma_2-\log\gamma_1} \leq \frac{\frac{1}{2}\log\left(1+2c_a^2\right) - \frac{1}{2}\log \gamma_1}{\log\gamma_2-\log\gamma_1} \to \frac{1}{2}
\end{equation}
as $\gamma_1\downarrow0$. Therefore rearranging \eqref{eqn:RearrageMagnetization} gives
\begin{equation}
    M(\lambda^*,\gamma_2,a) \geq \frac{1}{2\lambda^*c_a}\left(\lambda^* - \lambda_T(a)\right).
\end{equation}
Then taking $\gamma_2\downarrow0$ gives
\begin{equation}
    \theta_{\lambda^*}(a) \geq \frac{1}{2\lambda^*c_a}\left(\lambda^* - \lambda_T(a)\right) = \frac{1}{2\lambda^*\left(1+ \lambda^* \InfNorm{D}\Jcal_{\lambda_T(a),a}\right)}\left(\lambda^* - \lambda_T(a)\right)
\end{equation}
as required.
\end{proof}

\begin{proof}[Proof of Theorem~\ref{thm:Percolation Mean-Field Bound}]
From Assumption~\ref{Assump:BoundExpectedDegree} we have $\lambda_T^{(1)}=\lambda_T$ and $\int D(a,b)\Pcal(\dd b)<\infty$ for $\Pcal$-almost every $a\in\Ecal$. Assumption~\ref{Assump:AllReachablebySome} implies that there exists a $\Pcal$-positive set (say, $E\subset\Ecal$) on which $\sup_{k\geq 1}\essinf_{b\in\Ecal} D^{(k)}(a,b)>0$. Therefore Lemma~\ref{lem:EqualCriticalIntensities} tells us that $\lambda_T(a)=\lambda_T$ on a $\Pcal$-positive set $E'\subset E$. Then Proposition~\ref{thm:PercolationBoundPointwise} with Assumption~\ref{Assump:BoundExpectedDegree} implies that there exists some $\Pcal$-positive set, $E''\subset E'$, on which $\theta_{\lambda_T}(a)>0$ or both $\theta_{\lambda_T}(a)=0$ and $\theta_\lambda(a) \geq \frac{1}{1+ \lambda \InfNorm{D}\Jcal_{\lambda_T(a),a}}\frac{1}{2\lambda}\left(\lambda - \lambda_T\right)_+$ for all $\lambda\geq0$.

If $\theta_{\lambda_T}(a)>0$ for any $\Pcal$-positive set of $a\in\Ecal$, then we have $\norm*{\theta_{\lambda_T}}_p>0$. Then since $\theta_{\lambda}(a)$ is non-decreasing in $\lambda$ for all $a\in\Ecal$, $\norm*{\theta_{\lambda}}_p$ is non-decreasing in $\lambda$ and our result is proven.

If $\theta_{\lambda_T}(a)=0$ for a $\Pcal$-almost every $a\in\Ecal$, then there exists a $\Pcal$-positive set $E'''\subset E''$ such that
\begin{equation}
    \norm*{\theta_\lambda}_p \geq \frac{1}{2\lambda}\left(\lambda -\lambda_T\right)_+\left(\int_{E'''}\left(\frac{1}{1+ \lambda \InfNorm{D}\Jcal_{\lambda_T,a}}\right)^{p}\Pcal(\dd a)\right)^\frac{1}{p}
\end{equation}
for all $\lambda\geq 0$. Assumption~\ref{Assump:BoundExpectedDegree} states that $\InfNorm{D}<\infty$, and so this last factor is strictly positive if and only if $\sup_{k\geq 1}\essinf_{b\in\Ecal} D^{(k)}(a,b)>0$ for a $\Pcal$-positive subset of $E'''$ and $\lambda_T>0$. However, \ref{Assump:AllReachablebySome} told us that this first inequality was true on $E\supset E'''$, and so this condition is true. For $\lambda_T>0$, note that \ref{Assump:BoundExpectedDegree} and Lemma~\ref{lem:lambdaO_lambdaT} imply that $\lambda_T=\lambda_O$. \ref{Assump:BoundExpectedDegree} also implies $\OpNorm{\Opconnf}<\infty$, and with Lemma~\ref{lem:nontriviallambdaO} this implies $\lambda_O>0$.

To get the non-dependence of the constants on $p$, note that $\norm*{\theta_\lambda}_p$ is minimised for $p=1$ by Jensen's inequality. We use the associated constant to get our result.
\end{proof}

\subsection{Percolation Upper Bound}
\label{sec:percolation_upper_bound}

We have the required lower bound for Theorem~\ref{thm:Percolation Mean-Field Behaviour} from Theorem~\ref{thm:Percolation Mean-Field Bound}, so we only need to prove the matching upper bound. We do this by taking the upper bound on the magnetization from Corollary~\ref{cor:square_root_bound_M}, and adapting the extrapolation principle used by Aizenman and Barsky in \cite{AizBar91}.

\begin{proof}[Proof of Theorem~\ref{thm:Percolation Mean-Field Behaviour}]
First we reparametrize the magnetization. Let $h = -\log(1 - \gamma)$ and for $\lambda>0$, $h\in\left(0,\log 2\right]$, and $a\in\Ecal$ define
\begin{equation}
    \widetilde{M}(\lambda,h,a) := M(\lambda,1 - \e^{-h},a).
\end{equation}
Then the first magnetization inequality from Lemma~\ref{lem:frst:diff_ineqs_magnetization} becomes
\begin{equation}
\label{eqn:tildeMagnetDiffInequal}
    \frac{\partial \widetilde{M}}{\partial \lambda} (\lambda,h,a) \leq \frac{\e^{-h}}{\lambda} \MagSupTilde(\lambda,h)\frac{1}{\e^{-h}}\frac{\partial \widetilde{M}}{\partial h} (\lambda,h,a) \leq \frac{1}{\lambda}\cbar{\lambda}\widetilde{M}\left(\lambda,h,a\right)\frac{\partial \widetilde{M}}{\partial h} (\lambda,h,a),
\end{equation}
where $\cbar{\lambda} = 1+ \lambda \InfNorm{D}\esssup_{a\in\Ecal}\Jcal_{\lambda,a}$ as before.

Fix $m \in [0,1]$, $\lambda \in \R_+$ and $a\in\Ecal$. Then by the continuity of $h\mapsto\widetilde{M}(\lambda,h,a)$ and the fact that $\widetilde{M}(\lambda,0,a) = 0$ and $\widetilde{M}(\lambda,\infty,a) = 1$ we know that we can solve the equation $\widetilde{M}(\lambda,h,a) = m$ for $h$, i.e, there exists $h := h(\lambda)$ so that $\widetilde{M}(\lambda,h(\lambda),a) = m$. Furthermore, since $h\mapsto\widetilde{M}(\lambda,h,a)$ is differentiable and strictly increasing, we know that $h(\lambda)$ is differentiable too. The mark $a$ is fixed for this part of the argument, and so we omit the dependence of $h(\lambda)$ upon $a$. By differentiating this equality we get 
\begin{equation}
    \frac{\partial \widetilde{M}}{\partial \lambda} + \frac{\partial \widetilde{M}}{\partial h}\frac{\partial h}{\partial \lambda}\bigg|_{\widetilde{M} = m} = 0.
\end{equation}
By combining this with \eqref{eqn:tildeMagnetDiffInequal} we get 
\begin{equation}
\label{ineq: slope magn}
0 \leq -\frac{\partial h}{\partial \lambda}\bigg|_{\widetilde{M} = m} \leq \frac{m\cbar{\lambda}}{\lambda}.
\end{equation}

We aim to describe $\theta_\lambda$ for $\lambda = \lambda_c + \varepsilon$ for some $\varepsilon>0$. In that case, we can further bound \eqref{ineq: slope magn} from above by $\frac{m\cbar{\lambda_c+\varepsilon}}{\lambda_c}$. Given some $h_0>0$, let us denote $\Lambda_1 = \left(\lambda_c + \varepsilon,h_0,a\right)$ and $m_1 = \widetilde{M}(\Lambda_1)$. Then define
\begin{equation}
    \varepsilon' = \varepsilon + \frac{\lambda_c h_0}{m_1 \cbar{\lambda_c+\varepsilon}},
\end{equation}
and further define $\Lambda_2 = \left(\lambda_c,\frac{m_1\cbar{\lambda_c+\varepsilon}\varepsilon'}{\lambda_c},a\right)$ and $m_2=\widetilde{M}(\Lambda_2)$. The key observation is then that $m_1 \leq m_2$. We prove this by considering the $\left(\lambda-h\right)$-plane and comparing $\Lambda_2$ to the contour line of $\widetilde{M}=m_1$. By construction, the straight line passing through $\Lambda_1$ and $\Lambda_2$ has gradient $-\frac{m_1\cbar{\lambda_c+\varepsilon}}{\lambda_c}$. Comparing this with \eqref{ineq: slope magn} shows that this straight line has a steeper (or equal) gradient than the contours of constant $\widetilde{M}$ - in particular that contour that passes through $\Lambda_1$. Therefore $\Lambda_2$ does not lie below the contour of $\widetilde{M}=m_1$. Since $\widetilde{M}$ is non-decreasing in $h$, this means that $m_2=\widetilde{M}\left(\Lambda_2\right) \geq m_1$.

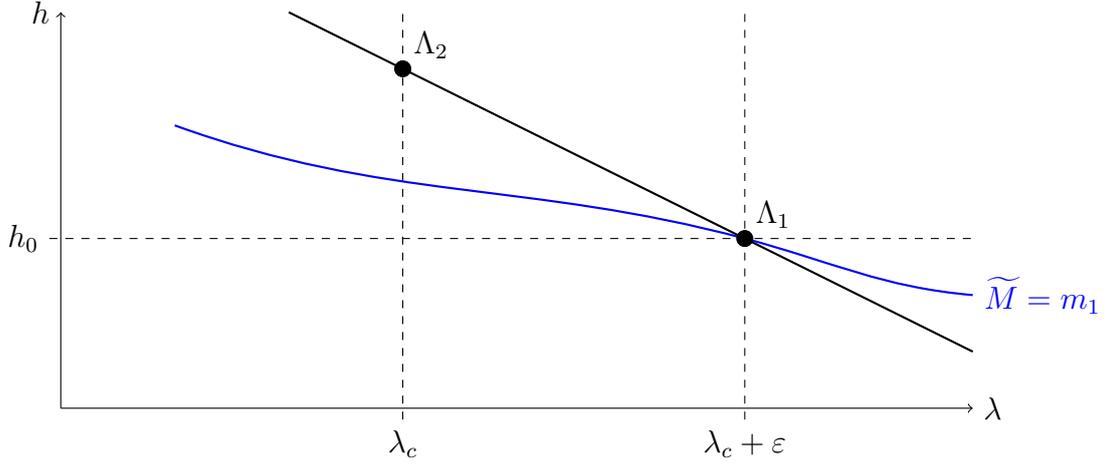
\begin{figure}
    \centering
    \begin{tikzpicture}[scale=1.4]
        \draw[->] (0,0.5) -- (0,4) node[left]{$h$};
        \draw[->] (0,0.5) -- (8,0.5) node[right]{$\lambda$};
        \draw[dashed] (3,0.4) node[below]{$\lambda_c$} -- (3,4);
        \draw[dashed] (6,0.4) node[below]{$\lambda_c+\varepsilon$} -- (6,4);
        \draw[dashed] (-0.1,2) node[left]{$h_0$} -- (8,2);
        \draw[thick] (8,1) -- (2,4);
        \draw[thick, blue] (8,1.5) node[right]{$\widetilde{M}=m_1$} to [out=175,in=345] (6,2) to [out=165,in=340] (1,3);
        \filldraw (6,2) circle (2pt) node[above right]{$\Lambda_1$};
        \filldraw (3,3.5) circle (2pt) node[above right]{$\Lambda_2$};
    \end{tikzpicture}
    \caption{Sketch showing $m_2\geq m_1$. The $\widetilde{M}=m_1$ contour has a shallower slope everywhere than the line connecting $\Lambda_1$ and $\Lambda_2$, and so $\Lambda_2$ lies above the $\widetilde{M}=m_1$ contour.}
    \label{fig:my_label}
\end{figure}

Since we have $\theta_{\lambda_c}(a)=0$ for $\Pcal$-almost every $a\in\Ecal$, we can apply Corollary~\ref{cor:square_root_bound_M} for these $a\in\Ecal$. This gives
\begin{align}
    \widetilde{M}\left(\lambda_c+\varepsilon,h_0,a\right) = m_1 &\leq m_2 = \widetilde{M}\left(\lambda_c,\frac{m_1\cbar{\lambda_c+\varepsilon}\varepsilon'}{\lambda_c},a\right)\nonumber\\
    &\leq \sqrt{K}\left(1 - \exp\left(-\frac{m_1\cbar{\lambda_c+\varepsilon}\varepsilon'}{\lambda_c}\right)\right)^{\frac{1}{2}}\nonumber\\
    & = \sqrt{K}\left(1 - \exp\left(-\frac{m_1\cbar{\lambda_c+\varepsilon}\varepsilon}{\lambda_c}\right) + \exp\left(-\frac{m_1\cbar{\lambda_c+\varepsilon}\varepsilon}{\lambda_c}\right)\left(1-\e^{-h_0}\right)\right)^{\frac{1}{2}}\nonumber\\
    &\leq \sqrt{K}\left(\frac{m_1\cbar{\lambda_c+\varepsilon}\varepsilon}{\lambda_c} + \gamma_0\right)^{\frac{1}{2}},
\end{align}
where in the last inequality we have used $1-\exp(x) \leq x$ and $\exp(x)\leq 1$ for all $x\geq 0$. We therefore want to solve the quadratic inequality
\begin{equation}
    M\left(\lambda_c+\varepsilon,\gamma_0,a\right)^2 - KM\left(\lambda_c+\varepsilon,\gamma_0,a\right)\frac{\cbar{\lambda_c+\varepsilon}\varepsilon}{\lambda_c} - K\gamma_0 \leq 0
\end{equation}
for $M\left(\lambda_c+\varepsilon,\gamma_0,a\right)$. Since $M\left(\lambda_c+\varepsilon,\gamma_0,a\right)\geq 0$, this inequality becomes
\begin{equation}
    M\left(\lambda_c+\varepsilon,\gamma_0,a\right) \leq \frac{K\varepsilon\cbar{\lambda_c+\varepsilon}}{\lambda_c}+\sqrt{K\gamma_0}.
\end{equation}
Taking $\gamma_0\downarrow 0$ then gives
\begin{equation}
    \theta_{\lambda_c + \epsilon}(a) \leq \frac{K\cbar{\lambda_c+\varepsilon}}{\lambda_c}\varepsilon \leq 2\frac{K\cbar{\lambda_c + \epsilon}}{\lambda_c}\varepsilon
\end{equation}
where this last inequality holds for $\varepsilon$ sufficiently small. Taking the essential supremum over $a\in\Ecal$ then produces the result.
\end{proof}

\section{Cluster Tail Proof}
\label{sec:ClusterTailProof}

Here we prove both the upper and lower bounds in Theorem~\ref{thm:Cluster Tail Mean-Field Behaviour} using the upper bound for the magnetization from Corollary~\ref{cor:square_root_bound_M}. Therefore our proofs for the upper and lower bounds both require Assumption~\ref{TriangleCondition_Assumption}.

\begin{proof}[Proof of Theorem \ref{thm:Cluster Tail Mean-Field Behaviour}]
We first establish that $\chi^{\mathrm{f}}_{\lambda_c}(a)=\infty$ for $\Pcal$-almost every $a\in\Ecal$. Via Lemma~\ref{thm:equality of lambda_T} and Corollary~\ref{thm:Sharpness of Phase Transition}, we know that the Assumptions \ref{Assump:BoundExpectedDegree} and \ref{Assump:AllReachablebySome} imply that $\lambda_O=\lambda_T=\lambda_c$. It then follows from \cite[Proposition~7.2]{DicHey2022triangle} that \ref{TriangleCondition_Assumption} implies $\theta_{\lambda_c}(a)=0$ for $\Pcal$-almost every $a\in\Ecal$. This - with Lemma~\ref{lem:CriticalSusceptibility} - then implies that $\chi^{\mathrm{f}}_{\lambda_c}(a)=\chi_{\lambda_c}(a)=\infty$ for $\Pcal$-almost every $a\in\Ecal$. We are therefore able to use Corollary~\ref{lem:lowerboundmagnetization} and Corollary~\ref{cor:square_root_bound_M} to get lower and upper bounds on the magnetization for $\Pcal$-almost every $a\in\Ecal$.

Note that
\begin{align}
    \p_{\lambda_c} \left(\left|\C\left(\origin{a}\right)\right| \geq n\right) &=  \frac{\e}{\e-1} \sum_{l \geq n} \left(1- \tfrac 1e\right) \p_{\lambda_c} \left(\left|\C\left(\origin{a}\right)\right| = l\right) \nonumber\\
	& \leq  \frac{\e}{\e-1} \sum_{l \geq n} \left(1- \left(1-\tfrac 1n\right)^l\right) \p_{\lambda_c} \left(\left|\C\left(\origin{a}\right)\right| = l\right) \nonumber\\ 
	& \leq \frac{\e}{\e-1} M\left(\lambda_c,\tfrac 1n,a\right).
\end{align}
With Corollary~\ref{cor:square_root_bound_M} we can bound the magnetization above to get
\begin{equation}
    \p_{\lambda_c} (|\C(\origin{a})| \geq n) \leq \frac{\e}{\e-1}\sqrt{K}n^{-1/2}.
\end{equation}

To prove the lower bound, let $0 \leq \gamma < \tilde\gamma < 1$. Then
\begin{align}
    M(\lambda_c,\gamma,a) & \leq \gamma \sum_{l < n} l \p_{\lambda_c} \left(\left|\C\left(\origin{a}\right)\right| = l\right) + \sum_{l \geq n} \p_{\lambda_c} \left(\left|\C\left(\origin{a}\right)\right| = l\right) \nonumber\\
	& \leq \frac{\gamma}{\tilde\gamma} \sum_{l < n} \e^{l \tilde\gamma} \left(1-\left(1-\tilde\gamma\right)^l\right) \p_{\lambda_c} \left(\left|\C\left(\origin{a}\right)\right| = l\right) + \p_{\lambda_c} \left(\left|\C\left(\origin{a}\right)\right| \geq n\right) \nonumber\\
		& \leq \frac{\gamma}{\tilde\gamma} \e^{\tilde\gamma n} M(\lambda_c, \tilde{\gamma},a) + \p_{\lambda_c} \left(\left|\C\left(\origin{a}\right)\right| \geq n\right).
\end{align}
In the above, we used that $\left(1-\left(1-\gamma\right)^l\right) \leq l\gamma$ in the first bound, as well as
\begin{equation}
    l \tilde\gamma \leq \e^{l \tilde\gamma} - 1 = \e^{l \tilde\gamma} \left( 1- \left(\e^{- \tilde\gamma}\right)^l\right) \leq \e^{l \tilde\gamma} \left(1- \left(1-\tilde\gamma\right)^l\right)
\end{equation}
in the second bound. Plugging in $\tilde\gamma = 1/n$ and using Corollary \ref{cor:square_root_bound_M}, we obtain
\begin{equation}
    \p_{\lambda_c} \left(\left|\C\left(\origin{a}\right)\right| \geq n\right) \geq M(\lambda_c,\gamma,a) - n\gamma\e\sqrt{K}n^{-\frac{1}{2}}.
\end{equation}
We now use Corollary~\ref{lem:lowerboundmagnetization} to lower bound the first term, so that
\begin{equation}
    \p_{\lambda_c} \left(\left|\C\left(\origin{a}\right)\right| \geq n\right) \geq \sqrt{\gamma} \left(\sqrt{\left(1+\cbar{\lambda_c}^2\right)^{-1}} - \e\sqrt{K n\gamma } \right).
\end{equation}
Choosing $\gamma = \frac{1}{4K\e^2\left(1+\cbar{\lambda_c}^2\right)} n^{-1}$ shows that
\begin{equation}
    \p_{\lambda_c} \left(\left|\C\left(\origin{a}\right)\right| \geq n\right) \geq C \sqrt{1/n}
\end{equation}
for $C := \frac{1}{4\e\sqrt{K}}\left(1+\cbar{\lambda_c}^2\right)^{-1}$ .
\end{proof}





\newcommand{\etalchar}[1]{$^{#1}$}

\begin{acks}
This work is supported by  \textit{Deutsche Forschungsgemeinschaft} (project number 443880457) through priority program ``Random Geometric Systems'' (SPP 2265). The authors would also like to thank Markus Heydenreich for his advice and many interesting discussions.
\end{acks}

\end{document}